\documentclass[11pt]{article}
\usepackage{amsmath,amssymb,amsthm,mathtools,bm}
\usepackage[margin=1in]{geometry}
\usepackage{enumitem}
\usepackage[utf8]{inputenc}
\usepackage{verbatim}
\usepackage{graphicx}
\usepackage{subcaption}

\usepackage[colorlinks,citecolor=blue,urlcolor=blue,linkcolor=blue]{hyperref}

\numberwithin{equation}{section}
\newtheorem{Definition}{Definition}[section]
\newtheorem{Remark}{Remark}[section]
\newtheorem{Theorem}{Theorem}[section]
\newtheorem{Lemma}{Lemma}[section]
\newtheorem{Proposition}{Proposition}[section]
\newtheorem{Corollary}{Corollary}[section]
\newtheorem{Assumption}{Assumption}[section]

\newcommand{\be}{\begin{equation}}
 \newcommand{\ee}{\end{equation}}
\newcommand{\bee}{\begin{equation*}}
 \newcommand{\eee}{\end{equation*}}
\newcommand{\bi}{\begin{itemize}}
 \newcommand{\ei}{\end{itemize}}

\DeclareMathOperator*{\argmax}{arg max}

\def \E{\mathbb{E}}

\def \R{\mathbb{R}}

\def \Lc{{\mathcal L}}
\def \Cc{{\mathcal C}}
\def \Sc{{\mathcal S}}
\def \Ac{{\mathcal A}}
\def \Bc{{\mathcal B}}
\def \Ec{{\mathcal E}}
\def \Pc{{\mathcal P}}

\def \Hc{\mathcal{H}}

\def \eps{\varepsilon}
\def \rtild{\widetilde{r}}
\def \Leb{\operatorname{Leb}}

\newcommand{\thl}{\underline\eta}
\newcommand{\thu}{\overline\eta}
\newcommand{\norm}[1]{\left\|#1\right\|}

\newcommand{\vpi}{\theta} 

\title{Policy Iteration Convergence for Exploratory Two-State Mean-Field Control Problems}
\author{
 Zhenhua Wang\thanks{Zhongtai Securities Institute for Financial Studies, Shandong University; Shandong Province Key Laboratory of Financial Risk, Jinan, Shandong, China. \url{zhenhuaw@sdu.edu.cn}}
\and Xiaoli Wei \thanks{Institute for Advance Study in Mathematics, Harbin Institute of Technology, Harbin, China. \url{ xiaoli.wei@hit.edu.cn,}}
\and Xiang Yu\thanks{Department of Applied Mathematics,  The Hong Kong Polytechnic University, Kowloon, Hong Kong. \url{xiang.yu@polyu.edu.hk}}}

\date{\vspace{-0.5in}}

\begin{document}
\maketitle

\begin{abstract}
This paper studies the theoretical policy iteration convergence for two-state mean-field control (MFC) problems with Wright--Fisher common noise under entropy regularization. The population dynamics lead to a degenerate Hamilton--Jacobi--Bellman equation of Kimura type. After an arcsine transformation, the degeneracy becomes a singular Bessel drift whose coefficients and  weights depend on the feedback stochastic policy. This policy dependence prevents a direct application of estimates for a fixed operator. We establish weighted Sobolev estimates and H\"older embeddings uniformly over the admissible policy class. Combining these estimates with logarithmic entropy bounds and boundary estimates that remain valid when the policy changes, we prove the well-posedness, monotone improvement and convergence of policy iteration to the optimal solution. 
We also establish the convergence to the original MFC problem when the entropy parameter $\lambda$ tends to zero. Both the regularization error and the performance loss of the regularized optimal policy are bounded by $O(\lambda(1+|\ln\lambda|))$ in a weighted H\"older norm. An epidemic containment example of MFC is presented to illustrate the policy iteration convergence and the vanishing entropy convergence.
\end{abstract}


\section{Introduction}\label{sec:introduction}

Recently, continuous-time reinforcement learning (RL) theories and algorithms have attracted considerable attention in the communities of stochastic control and quantitative finance. The dynamic exploitation in RL to achieve the maximum reward naturally connects the learning problem to continuous-time stochastic control, in which some mathematical tools such as relaxed control, dynamic programming and martingale characterization can be invoked again in the new context to understand the learning mechanism. The pioneer work Wang, Zariphopoulou and Zhou \cite{WangZariphopoulouZhou2020} proposed the entropy regularization to the so-called exploratory control formulation to facilitate the characterization of the optimal stochastic policy. Later, Jia and Zhou \cite{jiazhou2022a, jiazhou2022b, jiazhou2023} laid theoretical foundations for policy evaluation, policy gradient and  policy iteration RL algorithms in learning the optimal policy in classical diffusion models. The fast growth in continuous-time RL  also motivated interesting theoretical studies such as the model-based policy iteration convergence. In the classical stochastic control setting, Kerimkulov, \v{S}i\v{s}ka and Szpruch \cite{KerimkulovSiskaSzpruch2020} established exponential convergence and stability of of policy iteration for drift-controlled diffusions on a finite horizon, using BSDE techniques. In the exploratory formulation, Huang, Wang and Zhou \cite{HuangWangZhou2025} proved the policy iteration convergence over an infinite horizon by combining Sobolev and H\"older estimates with a sublinear entropy bound. Tran, Wang and Zhang \cite{TranWangZhang2025} studied the convergence of policy iteration with controlled diffusion or unbounded coefficients under a large-discount condition. Ma, Wang and Zhang \cite{MaWangZhang2026} developed a probabilistic approach based on representations of solutions and their derivatives and obtained quantitative convergence under suitable assumptions. Huang, Yu and Zhang \cite{HYZ26} employed this probabilistic approach in the exploratory time-inconsistent control problem and proved the convergence of policy iteration towards the equilibrium policy.

Meanwhile, continuous-time RL theories and algorithms have been rapidly generalized in various applications beyond the classical models, among which the extension of continuous-time RL to mean-field control (MFC) has been an appealing research direction with many exclusive new features. MFC concerns the cooperative population where all agents share the same goal to attain the social optimum, which is often interpreted as an optimal control problem by a social planner on behalf of interacting agents. Just to name a few, Wei and Yu \cite{weiyu2025} firstly extended the q-learning theory to MFC in the McKean-Vlasov diffusion model by developing the integrated q-function and model-free policy iteration algorithm; Ren, Wei, Yu and Zhou \cite{RWYZ26a, RWYZ26b} further established the integrated q-function, two-layer fixed point characterization of the optimal policy and model-free policy iteration algorithms in MFC with controlled common noise; Frikha, Germain, Lauri\`ere, Pham and Song \cite{FGLPS25} studied the stochastic policy gradient algorithm in MFC problems; Bayraktar, Hennandez, Yan, Zhu \cite{BHYZ26} investigated the stochastic policy gradient for ergodic MFC problems; Cheng, Guo, Pham and Zhang \cite{ChengGPZ} proposed a deterministic policy gradient algorithm for extended MFC problems. In the continuous-time framework, the mean-field HJB equation involves the derivatives of the value function in Wasserstein space with respect to the population law. Consequently, the standard PDE estimates and probabilistic arguments in the aforementioned studies on the policy iteration convergence in the single agent's RL problem cannot be applied to the mean-field model, thereby leaving the convergence of policy iteration for MFC in McKean-Vlasov diffusion model still an open problem. This gap between the well-performed mean-field RL algorithms and the under-developed convergence analysis of model-based policy iteration motivates the present study. To circumvent the profound challenge due to the lack of proper estimates for the derivatives of the value function with respect to probability measures, we consider the MFC problem with finite states under entropy regularization. To facilitate the arguments and presentation, we focus on the case of two-state model in this work. The general cases with multiple states will be left for future study.

Increasing research interests have been seen in finite-state MFC and mean-field game (MFG) problems. Cecchin \cite{Cecchin2021} studied finite-horizon finite-state MFC, characterizes the value function as a viscosity solution of an HJB equation on the simplex, and establishes convergence from centralized finite-agent control problems. Bayraktar, Cecchin, Cohen and Delarue \cite{BayraktarCecchinCohenDelarue2021} introduced Wright--Fisher common noise to establish uniqueness in finite-state MFGs through a nonlinear Kimura master equation. Their subsequent work \cite{BayraktarCecchinCohenDelarue2022} derived this formulation as a limit of weighted finite-player games. In the potential setting, Cecchin and Delarue \cite{CecchinDelarue2022} examined a stochastic MFC problem with Wright--Fisher common noise, establish classical solvability of its parabolic HJB equation under sufficiently strong boundary repulsion, and use vanishing common noise to select MFG equilibria. These studies demonstrate the role of Kimura operators in mean-field models. 

Notably, even in the mean-field model with finite states, the social planner's state variable is the population distribution, which takes values in a continuous probability simplex. In the presence of common noise, this distribution becomes stochastic, and the associated Hamilton--Jacobi--Bellman (HJB) equation contains second-order terms. Wright--Fisher common noise has a covariance structure adapted to the simplex and degenerates at its boundary, leading to differential operators of Kimura type. This boundary degeneracy creates substantial difficulties in studying the regularity of value functions and the convergence of the policy iteration. In particular, the previous studies in single agent's model show the importance of uniform derivative estimates in policy iteration. However, their results do not directly apply to the boundary structure considered in the finite-state MFC problem: the diffusion is degenerate in the original population variable, while the transformed equation has a singular drift that may be merely measurable.

The present paper aims to establish convergence of a policy iteration algorithm (PIA) for an infinite-horizon discounted two-state MFC problem with Wright--Fisher common noise. The social planner chooses a distribution of actions for each individual state, and Shannon entropy is added to the population-weighted reward. This exploratory formulation yields Gibbs policies and a PIA that alternates linear policy evaluation with pointwise policy improvement. We study convergence to the regularized optimum and quantify the performance of the resulting policy in the original problem as the entropy parameter vanishes. To the best of our knowledge, this is the first work to establish convergence of policy iteration for MFC with Wright--Fisher common noise. A central difficulty is that the boundary behavior of the policy-evaluation operator depends on the feedback policy. Consequently, both the operator and the weighted spaces used to estimate its solutions vary along the iteration, requiring estimates that are uniform over the admissible policy class.

Apart from PIA convergence, the present work is also related to the vanishing entropy regularization approach in the literature. It has been observed in some recent studies that the Gibbs
measure characterization of the optimal relaxed control under entropy regularization brings the improvement of solution regularity and simplifications in fixed point arguments. By further establishing the convergence
as entropy regularization vanishes, it offers a new route to confirm the existence of solution in the 
original continuous-time control or game problems. Tang, Zhang and Zhou \cite{TangZhangZhou2022} proved the convergence of the exploratory HJB equation to the classical control problem using the viscosity solution argument. In linear-quadratic MFG, Guo, Xu and Zariphopoulou \cite{GuoXuZariphopoulou2022} showed the convergence of equilibrium policy to the solution in the original MFG as temperature parameter tends to zero. In the time-inconsistent setting, Wang, Yu, Zhang and Zhou \cite{WangYuZhangZhou2026} established the existence of equilibria in the original problem by the vanishing entropy regularization approach. Bayraktar, Wang, Yu and Zhang \cite{BayraktarWangYuZhang2026} extended this approach to time-inconsistent MFGs with uniformly elliptic dynamics.  It is worth noting that Cohen and Priya \cite{CohenPriya2026} recently studied a finite-state MFC problem with deterministic population dynamics and a relative-entropy penalty on unbounded transition intensities, obtaining the convergence result as entropy tends to zero. In contrast, the entropy regularization in our work is for the purpose of randomization of stochastic policy, while Wright--Fisher common noise leads to a second-order Kimura equation. In our setting, how to show the convergence as entropy vanishes deserves some careful investigations.

Our methodology and main theoretical contributions are outlined as follows.

Firstly, we establish uniform weighted Sobolev estimates and H\"older embeddings for a family of Bessel operators with measurable coefficients; see Theorem~\ref{thm:main.besselest}. The analytical starting point is the one-dimensional diffusion theory of Feller \cite{Feller1951} and the Kimura theory developed by Epstein and Mazzeo \cite{EpsteinMazzeo2010,EpsteinMazzeo2013}. We use the global arcsine transformation $s=\Phi(p)=2\arcsin\sqrt p$, which maps $[0,1]$ onto $[0,\pi]$ and converts the Kimura coefficient $p(1-p)$ into a constant second-order coefficient. This transformation realizes the standard one-dimensional Wright--Fisher geometry globally, and the norm in \eqref{eq:CK.original} measures H\"older regularity of the value and its intrinsic first derivative in this geometry. The transformed drift contains singular terms proportional to $\cot(s/2)$ and $\tan(s/2)$, with coefficients determined by averages of the action-dependent boundary transition rates under the feedback policy. The associated weights $m^\theta$ and $\widehat m^\theta$ therefore depend on the policy, and estimates in the weighted space for a fixed policy are insufficient for the subsequent compactness argument. To address this difficulty, we establish uniform weight-ratio and logarithmic-moment estimates in Lemma~\ref{lm:ratio}, a weighted Hardy inequality in Lemma~\ref{lm:hardy}, and a uniform embedding in Proposition~\ref{prop:sobolev.embed}. The constants depend only on common bounds on the boundary parameters and the bounded drift remainders. These estimates place the policy values in a common H\"older space without requiring continuity of the feedback policy. This is relevant even when the primitive transition rates are smooth, since averaging them against a Borel feedback policy may produce a discontinuous effective drift.

Secondly, we establish the well-posedness of policy evaluation and the policy improvement on the entire state interval. For a Borel feedback stochastic policy, the interior differential expression alone need not determine the process at an accessible boundary. Using scale--speed theory and Dirichlet forms, Lemma~\ref{lm:closed.diffusion} specifies a conservative diffusion with zero boundary occupation time. Proposition~\ref{prop:PDE.wellposed} then identifies its stochastic payoff with the zero-flux resolvent solution, which yields a well-defined PIA for all initial states in the closed interval (see Corollary~\ref{cor:PIA.welldefined}). To compare consecutive policy values, we must also ensure that the value obtained under the current policy satisfies the zero-flux condition associated with the improved policy. This does not follow directly from the zero-flux condition for the current weight. We resolve this issue in Lemma~\ref{lm:endpoint.crossflux} by deriving endpoint derivative estimates of order $s(1+|\ln s|)$ near $0$, with the analogous estimate near $\pi$. These estimates imply zero flux for every weight in the admissible class and provide the boundary control needed for the policy improvement argument in Proposition~\ref{prop:policy.improvement} and the verification arguments.

Thirdly, we prove that the entire value sequence generated by PIA converges to the entropy-regularized optimum. Although monotonicity and boundedness give a pointwise limit, passing to the limit in the HJB equation requires compactness of the derivatives. An additional difficulty comes from the entropy term, which depends on the preceding value gradient and may become unbounded near the endpoints. Lemma~\ref{lm:entropy.est} bounds its growth logarithmically in the gradient, with a boundary contribution of order $|\ln(\sin s)|$. The required action-space assumption is a uniform local volume lower bound, which covers the uniform cone condition used in \cite{HuangWangZhou2025} without requiring a prescribed cone or a convex action set (see Remark~\ref{rem:action.volume}). The logarithmic moments in Lemma~\ref{lm:ratio} control the boundary contribution uniformly over the policy class. Combining these bounds with Theorem~\ref{thm:main.besselest}, we obtain a recursive H\"older estimate whose logarithmic growth yields a bound uniform in the iteration index and, moreover, in $\lambda\in(0,1]$ (see Theorem~\ref{thm:compact.un}). Starting from the uniform action density, the values thus converge in $\Cc^{1,\beta}([0,\pi])$ for every $0<\beta<\alpha$, and Theorem~\ref{thm:convergence} verifies that the Gibbs feedback associated with the limit is optimal. By the norm correspondence in \eqref{eq:CK.original}, the original values converge in $\Cc_K^{1,\beta}([0,1])$. This intrinsic norm controls the value and $\sqrt{p(1-p)}v_p$ up to the boundary; ordinary $\Cc^1$ regularity in $p$ at the endpoints is not required. The convergence argument applies without a large-discount condition or a smallness restriction on the entropy parameter.

Finally, we analyze the convergence as entropy regularization vanishes. Theorem~\ref{thm:vanishing.entropy} constructs and verifies the unregularized value function and establishes uniqueness of the HJB solution within the class satisfying the stated endpoint estimates. Writing $v^{*,\lambda}$ and $v^{*,0}$ for the regularized and unregularized optimal values, respectively, Theorem~\ref{thm:lambda.rate} and Corollary~\ref{cor:entropy.optimal.regret} give
\bee
\begin{aligned}
	\|v^{*,\lambda}-v^{*,0}\|_{\Cc_K^{1,\alpha}([0,1])}, 
	\|v^{*,0}-V_0^{\theta^{*,\lambda}}\|_{\Cc_K^{1,\alpha}([0,1])}
	&\leq C\lambda(1+|\ln\lambda|),\qquad 0<\lambda\leq1,
\end{aligned}
\eee
where $\theta^{*,\lambda}$ is an optimal regularized feedback and $C$ is independent of $\lambda$. The second estimate measures the performance loss when this feedback is applied to the problem without the entropy reward. To establish these bounds, we first estimate the difference between the regularized and unregularized Hamiltonians in weighted $L^k$ norms (see Lemma~\ref{lm:regularize.Hconvergence}), since the singular drift obstructs a direct uniform estimate in space. We then linearize the difference of the HJB equations. The secant slopes of the unregularized Hamiltonians lie between the corresponding action-dependent drift coefficients, so the resulting drift remains in the same Bessel class. The uniform resolvent and Sobolev estimates therefore yield the stated bounds for the values and their intrinsic first derivatives. These results justify the use of small positive entropy parameters to approximate the original control problem while retaining control up to the degenerate boundary.

The rest of the paper is organized as follows. Section~\ref{sec:model} formulates the two-state MFC problem. Section~\ref{sec:arcsin.transform} introduces the arcsine transformation and the resulting singular drift. Section~\ref{sec:Bessel.estimates} establishes the uniform Bessel estimates and embeddings. Section~\ref{sec:PIA} proves the well-posedness of PIA and the policy improvement property, and Section~\ref{sec:PIA.convergence} establishes convergence and optimality. Section~\ref{sec:convergence.entropy} studies the vanishing-entropy limit and derives the convergence rate. Finally, Section~\ref{sec:epidemic} presents an example of epidemic containment with numerical illustrations.

\section{Exploratory Two-State MFC with Kimura Operator}\label{sec:model}

 For some $\ell\in\mathbb N$, let $A\subset\mathbb R^\ell$ be an action space with $0<\operatorname{Leb}(A)<\infty$. We consider the two-state MFC with common noise such that the state process $X_t \in \{1, 2\}$. Given an action $a \in A$ and a population fraction $p \in [0, 1]$ of agents in state $1$, the agent's transition rate matrix is described by
 \bee
\begin{pmatrix}
-q^{1,2}(p,a) & q^{1,2}(p,a)\\
q^{2,1}(p,a) & -q^{2,1}(p,a)
\end{pmatrix},
\eee
where $q^{i,j} \geq 0$, $i \neq j$ is the transition rate from state $i$ to state $j$.

A relaxed Markov feedback policy is a pair $\theta=(\theta_1,\theta_2)$ with $\theta_i:[0,1]\to\Pc(A)$, $i=1,2$, where $\theta_i(p,da)$ is the action distribution of an agent in state $i$ when the population fraction in state $1$ is $p$. Let $P_t$ be the fraction of the population in state $1$ at time $t$ under $\theta$. Then $P_t$ evolves as the controlled Wright--Fisher diffusion
\be\label{eq:model.X}
\begin{aligned}
dP_t^{\theta}={}&2\Bigg[(1-P_t^{\theta})\int_A q^{2,1}(P_t^{\theta},a) \theta_2(P_t^{\theta},da)
-P_t^{\theta}\int_A q^{1,2}(P_t^{\theta},a) \theta_1(P_t^{\theta},da)\Bigg]dt+\sqrt{2P_t^{\theta}(1-P_t^{\theta})} dB_t,
\end{aligned}
\ee
where $B$ is the common Brownian noise. For each feedback policy, the diffusion on $[0,1]$ is specified by the interior generator in \eqref{eq:model.X}, together with the following boundary conditions: its lifetime is infinite, and it spends zero Lebesgue time at the boundary $\{0,1\}$. The existence and uniqueness in law of this diffusion, including for initial values $0$ and $1$, follow from Lemma \ref{lm:closed.diffusion} and the transformation in Remark \ref{rem:original.diffusion}. 

For a probability density $\vartheta$ on $A$, define its entropy by
$
\Hc(\vartheta):=-\int_A \vartheta(a)\ln\vartheta(a) da.
$
Denote $\Ac$ as the set of Borel relaxed feedback policies. For each $\vpi\in\Ac$, let $P_t^\theta$ denote the original-coordinate
diffusion obtained from Lemma \ref{lm:closed.diffusion} through Remark \ref{rem:original.diffusion}. In view that the rewards are bounded under Assumption~\ref{ass:action} (iv), the following expectation is finite for every initial point $p\in[0,1]$:
\be\label{eq:value.unregularized.x}
\begin{aligned}
V_0^{\vpi}(p):=\E_p^{\vpi}\Bigg[\int_0^\infty e^{-\rho t}\Bigg\{
&P_t^{\vpi}\int_A r(1,P_t^{\vpi},a) \vpi_1(P_t^{\vpi},da)+(1-P_t^{\vpi})\int_A r(2,P_t^{\vpi},a) \vpi_2(P_t^{\vpi},da)
\Bigg\}dt\Bigg].
\end{aligned}
\ee
Through the paper, $q^{i,j}(p,a)$, $i\neq j$, and $r(i,p,a)$, $i=1,2$ are assumed to be jointly Borel measurable in $(p,a)\in[0,1]\times A$. For fixed $\lambda>0$, let $\Ac_\lambda\subset\Ac$ be the set of feedback policies with finite entropy for a.e. $p\in(0,1)$ such that the following 
expected payoff is finite for $p\in[0,1]$:
\be\label{eq:value.regularized.x}
\begin{aligned}
V_\lambda^{\vpi}(p):=\E_p^{\vpi}\Bigg[\int_0^\infty e^{-\rho t}\Bigg\{
&P_t^{\vpi}\left[\int_A r(1,P_t^{\vpi},a)\vpi_1(P_t^{\vpi},a) da+\lambda\Hc(\vpi_1(P_t^{\vpi}))\right]\\
&+(1-P_t^{\vpi})\left[\int_A r(2,P_t^{\vpi},a)\vpi_2(P_t^{\vpi},a) da+\lambda\Hc(\vpi_2(P_t^{\vpi}))\right]
\Bigg\}dt\Bigg].
\end{aligned}
\ee
The optimal values are defined by
$
v^{*,0}(p):=\sup_{\vpi\in\Ac}V_0^{\vpi}(p)$ and
$v^{*,\lambda}(p):=\sup_{\vpi\in\Ac_\lambda}V_\lambda^{\vpi}(p)
$, respectively. For $\lambda>0$, the exploratory HJB equation on $(0,1)$ is written by
\be\label{eq:HJB.x}
\begin{aligned}
0=&-\rho v(p)+p(1-p)v_{pp}(p)\\
&+\sup_{\vartheta_1,\vartheta_2\in\Pc(A)}\Bigg\{p\int_A\Big[-2q^{1,2}(p,a)v_p(p)+r(1,p,a)-\lambda\ln\vartheta_1(a)\Big]\vartheta_1(a) da\\
&\qquad \qquad \qquad \quad \qquad+(1-p)\int_A\Big[2q^{2,1}(p,a)v_p(p)+r(2,p,a)-\lambda\ln\vartheta_2(a)\Big]\vartheta_2(a) da\Bigg\}.
\end{aligned}
\ee

\begin{Assumption}[Variable boundary weights with uniformly bounded remainder]\label{ass:b}
There exist some finite constants $\thl,\thu,M$ such that
\be\label{eq:main.boundaryrelaxed}
0<\thl\leq 2q^{2,1}(0, a), 2q^{1,2}(1, a)\leq\thu<\infty,\quad a\in A.
\ee
Moreover,
\be\label{eq:assum.remainder}
\begin{aligned}
&|2(1-p)[q^{2,1}(p,a)-q^{2,1}(0, a)]|\leq M\sqrt{p(1-p)},\\
&|2p[q^{1,2}(p,a)-q^{1,2}(1, a)]|\leq M\sqrt{p(1-p)},
\end{aligned}
\quad p\in(0,1),\ a\in A.
\ee
\end{Assumption}

\begin{Remark}
The endpoint conditions are one-sided: the left Kimura boundary weight is determined only by $q^{2,1}(0, a)$ and the right Kimura boundary weight is determined only by $q^{1,2}(1, a)$. The square-root conditions imply that, for every policy pair, the effective drift has a uniformly bounded remainder after subtraction of its state-dependent Jacobi profile. No regularity of the policy pair in $p$ is required for this conclusion.
\end{Remark}


\section{Global Arcsine Transformation}\label{sec:arcsin.transform}

Starting from this section, we employ the following change of variables, which corresponds to the intrinsic one-dimensional Wright--Fisher coordinate whose derivative is the reciprocal of the square root of the Kimura diffusion coefficient. Before introducing the transformation, let us briefly recall some related studies.

As noted in \cite{EpsteinMazzeo2010}, the one-dimensional normal-form transformation dates back to Feller \cite{Feller1951}. In \cite[Section 4]{EpsteinMazzeo2010}, Epstein and Mazzeo study the local behavior of Wright--Fisher operators near degenerate boundaries by decomposing the operator into a boundary-frozen model $L_{b_0}=x\partial_{xx}+b_0\partial_x$ and a lower-order perturbation. Under the further transformation $s=2\sqrt{x}$, a direct calculation shows that this model operator becomes
$$
\partial_{ss}+\frac{2b_0-1}{s}\partial_s,
$$
which is of Bessel type. This approach is developed extensively for generalized Kimura diffusions and their intrinsic square-root geometry in \cite{EpsteinMazzeo2013}. The global arcsine transformation $s=2\arcsin\sqrt{x}$, adopted in the present work, also appears in the literature on one-dimensional Wright--Fisher/Jacobi diffusions; see, for instance, \cite{NeuenkirchSzpruch2014,DucTranJost2018}.

The main distinction in our setting is: the singular drift coefficients depend on the action and, after averaging, on the feedback policy. Therefore, both the associated Bessel-type operators and their intrinsic weights vary with the policy and hence across policy iterations, while the intrinsic Wright--Fisher distance remains fixed. The convergence analysis therefore requires estimates uniform over this family of operators and weights.

A key contribution of this paper is to establish weighted Sobolev estimates and embeddings with constants uniform over the entire admissible policy class, despite the policy dependence of the singular coefficients and weights. These uniform estimates provide the essential regularity and compactness results crucially needed in establishing the well-posedness and convergence of the policy iteration algorithm.

Set $s=\Phi(p):=2\arcsin\sqrt{p}$ for $p\in[0,1]$. Then $\Phi$ is a homeomorphism from $[0,1]$ onto $[0,\pi]$, with inverse $p(s)=\sin^2(s/2)$.
Given $v(p)$, define $u(s):=v(p(s))$. For $p\in(0,1)$, equivalently $s\in(0,\pi)$, we have
$
u_s=\sqrt{p(1-p)}v_p$ and $u_{ss}=p(1-p)v_{pp}+\frac{1-2p}{2}v_p$. Then, it holds that 
$$
p(1-p)v_{pp} = u_{ss} - \frac{1-2p}{2\sqrt{p(1-p)}} u_s = u_{ss} + \left(\frac{p}{2} -\frac{1-p}{2} \right)\frac{u_s }{\sqrt{p(1-p)}}. 
$$
This, together with the weighted sum of the two Hamiltanian in \eqref{eq:HJB.x}, gives
\bee
\begin{aligned}
&-p2q^{1,2}(p, a)\frac{u_s}{\sqrt{p(1-p)}}+(1-p)2q^{2,1}(p, a)\frac{u_s}{\sqrt{p(1-p)}}+ \left(\frac{p}{2} -\frac{1-p}{2} \right)\frac{u_s}{\sqrt{p(1-p)}}\\
&= p\frac{-2q^{1,2}(p, a)+\frac12}{\sqrt{p(1-p)}} u_s+ (1-p) \frac{2q^{2,1}(p,a)-\frac12}{\sqrt{p(1-p)}} u_s
\end{aligned}
\eee 
For $i=1,2$, define
\be\label{eq:def.Q}
\begin{cases}
Q_1(s,a):=\frac{-2q^{1,2}(p(s), a)+\frac12}{\sqrt{p(s)(1-p(s))}},\\
Q_2(s,a):=\frac{2q^{2,1}(p(s), a)-\frac12}{\sqrt{p(s)(1-p(s))}}, \\
\widetilde r_i(s,a):=r(i,p(s), a)\quad \text{ for $i=1,2.$}
\end{cases}
\ee
The HJB equation \eqref{eq:HJB.x} becomes
\bee
\begin{aligned}
 0=&-\rho u+u_{ss}+p(s)\sup_{\vartheta_1\in\Pc(A)}\left\{\int_A[Q_1(s,a)u_s+\widetilde r_1(s,a)-\lambda\ln\vartheta_1(a)]\vartheta_1(a)da\right\}\\
 &+(1-p(s))\sup_{\vartheta_2\in\Pc(A)}\left\{\int_A[Q_2(s,a)u_s+\widetilde r_2(s,a)-\lambda\ln\vartheta_2(a)]\vartheta_2(a)da\right\},\quad s\in(0,\pi).
\end{aligned}
\eee
The two positive factors $p(s)$ and $1-p(s)$ do not affect the corresponding maximizers. The sum-weighted drift can be computed as 
$$
\begin{aligned}
	Q(p(s),a)=&pQ_1(s,a)+(1-p)Q_2(s,a)
	=  \frac{-p2q^{1,2}+2(1-p)q^{2,1} +(\frac{p}{2}-\frac{1-p}{2})}{\sqrt{p(1-p)}}  \\
	=& \frac{p \left[ -2q^{1,2}(1, a)+\frac{1}{2} \right] +p\left[ -2q^{1,2}(p, a) +2 q^{1,2}(1, a) \right]}{ \sqrt{p(1-p)}} \\
	&+ \frac{(1-p) \left[  2q^{2,1}(0, a)-\frac{1}{2} \right] +(1-p)\left[ 2q^{2,1}(p, a) - 2q^{2,1}(0, a) \right]}{ \sqrt{p(1-p)}}\\
	=& \tan\frac{s}{2} \left[ -2q^{1,2}(1, a)+\frac{1}{2} \right] + \cot\frac{s}{2} \left[  2q^{2,1}(0, a)-\frac{1}{2} \right]\\
	&+ \frac{p}{ \sqrt{p(1-p)}}\left[ -2q^{1,2}(p, a) +2 q^{1,2}(1, a) \right]+\frac{(1-p)}{ \sqrt{p(1-p)}}\left[ 2q^{2,1}(p, a) - 2q^{2,1}(0, a) \right].
\end{aligned}
$$
This calculation identifies where the singularities come from. The original drift in the $p$-variable is bounded at the boundary, but the Jacobian of the intrinsic change of variables is of order $p^{-1/2}$ at $0$ and $(1-p)^{-1/2}$ at $1$, and the left singularity is determined only by $q^{2,1}(0, a)$, while the right singularity is determined only by $q^{1,2}(1, a)$. That is,
\begin{align*}
	&\frac{2q^{2,1}}{\sqrt{p}} = \text{$O\left( \frac{1}{\sqrt{p}} \right)$ near 0 (i.e., after the arcsine transformation, $\cot{\frac s2}$ near 0 );}\\
	&\frac{2q^{1,2}}{\sqrt{1 - p}} = \text{$O\left( \frac{1}{\sqrt{1-p}} \right)$ near $1$ (i.e., after the arcsine transformation, $\tan\left( \frac s2 \right)$ near $\pi$)}.
\end{align*}
The remaining terms measure the deviation of the transition rates from their corresponding boundary values and are uniformly bounded by Assumption~\ref{ass:b}. 

We now make the decomposition precise. Set
\be\label{eq:def.q}
\begin{cases}
q_1(s,a):=&-2\sqrt{\frac{p}{1-p}} [q^{1,2}(p, a)-q^{1,2}(1, a)] =-2\tan\frac{s}{2}[q^{1,2}(p, a)-q^{1,2}(1, a)],\\
 q_2(s,a):=&2\sqrt{\frac{1-p}{p}} [q^{2,1}(p, a)-q^{2,1}(0, a)] =2\cot \frac{s}{2} [q^{2,1}(p, a)-q^{2,1}(0, a)].
\end{cases}
\ee
Then
$$
\begin{aligned}
 p(s)Q_1(s,a)&
 = -\left(2q^{1,2}(1, a)-\frac12\right)\tan\frac{s}{2}+q_1(s,a),\\
 (1-p(s))Q_2(s,a)&
 =\left(2q^{2,1}(0, a)-\frac12\right)\cot\frac{s}{2}+q_2(s,a).
\end{aligned}
$$
Thus the drift is decomposited as a singular part and a bounded part:
$$
\begin{aligned}
 Q(p(s),a):= &pQ_1(s,a)+(1-p)Q_2(s,a)=P(s,a) +q(s,a),
\end{aligned}
$$
with
$$
\begin{aligned}
&P(s, a):= -\tan\frac{s}{2} \left(2q^{1,2}(1, a)-\frac12\right)+ \cot\frac{s}{2} \left(2q^{2,1}(0, a)-\frac12\right), \\
& q(s,a):= q_1(s,a)+q_2(s,a),
\end{aligned}
$$

For a policy $\vpi\in\Ac$, we use the same notation for its representation in the transformed coordinate, namely $\vpi_i(s,\cdot):=\vpi_i(p(s),\cdot)$. For $\vpi\in\Ac_\lambda$, the same convention is used for its density $\vpi_i(s,a)$. Note that $\Phi$ is a homeomorphism, this identifies Borel feedback stochastic policies in the two coordinates bijectively; the distribution and the entropy of the stochastic policies are unchanged. We set
$$
W_0^{\vpi}(s):=V_0^{\vpi}(p(s)),\qquad W_\lambda^{\vpi}(s):=V_\lambda^{\vpi}(p(s)).
$$

For an admissible policy pair $\vpi=(\vpi_1,\vpi_2)$, define
$$
\begin{aligned}
&\Lambda_1^{\vpi}(s):=2\int_A q^{1,2}(1, a) \vpi_1(s,da),\quad
\Lambda_2^{\vpi}(s):=2\int_A q^{2,1}(0, a) \vpi_2(s,da),\\
&P^{\vpi}(s):=\left(\Lambda_2^{\vpi}(s)-\frac12\right)\cot\frac{s}{2}-\left(\Lambda_1^{\vpi}(s)-\frac12\right)\tan\frac{s}{2},\\
&q^{\vpi}(s):=\int_Aq_1(s,a) \vpi_1(s,da)+\int_Aq_2(s,a) \vpi_2(s,da).
\end{aligned}
$$
That is, $Q^{\vpi}(s)=P^{\vpi}(s)+q^{\vpi}(s)$. By Assumption~\ref{ass:b},
\be\label{eq:q.bound}
\sup_{\vpi\in\Ac}\|q^{\vpi}\|_{L^\infty((0,\pi))}\leq 2 M.
\ee
Then, $Q^{\vpi}(s)=\frac{2\Lambda_2^{\vpi}(s)-1}{s}+O(1)$ near $0$ and $Q^{\vpi}(s)=-\frac{2\Lambda_1^{\vpi}(s)-1}{\pi-s}+O(1)$ near $\pi$. As a result, the unbounded part $Q^{\vpi}$ is of Bessel type with policy-dependent exponents.

\begin{Assumption}\label{ass:action}
 Assume that:
 \begin{enumerate}[label=(\roman*)]
  \item The endpoint maps are uniformly Lipschitz in the action: there exists $L_b<\infty$ such that
  $$
  2|q^{2,1}(0, a)-q^{2,1}(0, a')|+2|q^{1,2}(1, a)-q^{1,2}(1, a')|\leq L_b|a-a'|.
  $$
  \item
  There exists $L_q<\infty$ such that
  $$
  |q_i(s,a)-q_i(s,a')|\leq L_q|a-a'|,\quad i=1,2,
  $$
  uniformly in $s\in(0,\pi)$ and $a,a'\in A$.
  \item The rewards are uniformly Lipschitz in the action: there exists $L_r<\infty$ such that
  $$
  \sup_{p\in(0,1), a,a'\in A}|r(i,p, a)-r(i,p, a')|\leq L_r|a-a'|,\quad i=1,2.
  $$
  \item The rewards are uniformly bounded: $\sup_{i,p, a}|r(i,p, a)|<\infty$.
  \item The action space $A$ is compact, and there exist constants $\gamma>0$ and $\iota\in(0,1)$ such that
  \bee
  \Leb(A\cap B_r(a))\geq\iota\Leb(B_r(a))=\iota\Leb(B_1(0)) r^\ell \quad \forall a\in A,\;\forall  0<r< \gamma.
  \eee
 \end{enumerate}
\end{Assumption}

\begin{Remark}\label{rem:action.volume}
Assumption~\ref{ass:action}(v) only requires that the action space $A$ occupies a fixed positive proportion of every sufficiently small ball centered at a point of $A$. It imposes no smoothness or convexity requirement on $A$, and does not require these intersections to contain a prescribed geometric shape. In particular, it results from the uniform cone condition in \cite[Assumption 4.10)]{HuangWangZhou2025} (see also \cite[Assumption 3]{BayraktarHuangWangZhou2025}).
In the one-dimensional case, the one-sided interval condition in \cite[Assumption 4.10]{HuangWangZhou2025} 
gives the same conclusion with $\iota=1/2$. Therefore, this assumption allows full-dimensional bounded polytopes and ellipsoids, as well as more general sets satisfying the stated volume bound without a uniform cone condition.
\end{Remark}

\section{The Bessel Operator and Two Uniform Estimates}\label{sec:Bessel.estimates}
This section devotes to provide \textit{uniform} Sobolev embedding and Sobolev estimate on singular operator of Bessel type particularly tailed for the PIA, which is non-trivial and new, and may have their own interests.

Within this section, we consider the PDE:
\be\label{eq:pde.bessel}
-\rho u+\Bc u+q(s)u_s=f,\quad \Bc u:=u_{ss}+P(s)u_s,
\ee
where $\norm q_{L^\infty((0,\pi))}\leq M$ and
\be\label{eq:Pdef}
 P(s)=\left(\Lambda_2(s)-\frac12\right)\cot\frac s2-\left(\Lambda_1(s)-\frac12\right)\tan\frac s2.
\ee
Define
$
m(s):=\exp\left(\int_{\pi/2}^{s}P(t) dt\right)$ and
$m^q(s):=m(s)\exp\left(\int_{\pi/2}^{s}q(t) dt\right)$.
Then
\be\label{eq:P.transform}
\Bc u=\frac1m(mu_s)_s \text{ and } u_{ss}+(P+q)u_s=\frac1{m^q}(m^q u_s)_s \text{ in the distributional sense.  }
\ee
Given $\norm q_{L^\infty}\leq M$, then
\begin{equation}\label{eq:measure-equivalence}
 e^{-M\pi}m(s)\leq m^q(s)\leq e^{M\pi}m(s),\quad s\in(0,\pi).
\end{equation}

Here, let us collect the function spaces and norms used throughout the paper. For a positive locally integrable weight $w$ on an interval $I$ and $1\leq k<\infty$, let $L^k(I,w\,ds)$ be the space of measurable functions, identified almost everywhere with finite norm
\bee
\|f\|_{L^k(I,w\,ds)}:=\left(\int_I |f(s)|^k w(s)\,ds\right)^{1/k}.
\eee
We abbreviate $L^k((0,\pi),w\,ds)$ as $L^k(w\,ds)$, while omitting the weight means integration with respect to Lebesgue measure. The space $L^\infty(I)$ carries the essential-supremum norm, which agrees with the supremum norm for continuous functions. For $1\leq k\leq\infty$, the unweighted spaces $W^{j,k}(I)$, $j=1,2$ carry the norm $\sum_{\ell=0}^j\|\partial_s^\ell u\|_{L^k(I)}$, with weak derivatives. The subscript $\mathit{loc}$ indicates that this regularity holds on every compact subinterval.
For $1<k<\infty$, the weighted Sobolev space $W^{2,k}(w\,ds)$ consists of $u\in W^{2,k}_{\mathrm{loc}}((0,\pi))$ such that $u,u_s,u_{ss}\in L^k(w\,ds)$, with the norm given by the sum of these three $L^k$ norms. In particular,
\begin{equation}\label{eq:W2pdef}
\norm{u}_{W^{2,k}(m ds)}:=\norm{u}_{L^k(m ds)}+\norm{u_s}_{L^k(m ds)}+\norm{u_{ss}}_{L^k(m ds)}.
\end{equation}
The same convention applies to $m^q$, $m^\vpi$ and $\widehat m^\vpi$, with the policy-dependent weights defined below. By \eqref{eq:measure-equivalence}, the norms with weights $m$ and $m^q$ are uniformly equivalent.

We also apply the following H\"older seminorm and norms $[0,\pi]$:
\bee
\begin{aligned}
 [u]_{\Cc^\alpha([0,\pi])}&:=\sup_{s\neq t, s,t\in [0,\pi]}\frac{|u(s)-u(t)|}{|s-t|^\alpha},  \quad  \norm{u}_{\Cc^{l}([0,\pi])}:=\sum_{j=0}^l\norm{u^{(j)}}_{L^\infty([0,\pi])}\quad \; \text{for } l=0,1,\\
 \|u\|_{\Cc^\alpha([0,\pi])}&:=\|u\|_{L^\infty([0,\pi])}+[u]_{\Cc^\alpha([0,\pi])},\quad 
\|u\|_{\Cc^{1,\alpha}([0,\pi])}:=\|u\|_{\Cc^\alpha([0,\pi])}+\|u_s\|_{\Cc^\alpha([0,\pi])}.
\end{aligned}
\eee
and $\Cc([0,\pi]), \Cc^{1}([0,\pi]), \Cc^{\alpha}([0,\pi]), \Cc^{1,\alpha}([0,\pi])$ consist of the continuous functions with the corresponding finite seminorms and norms respectively.
The same conventions apply on compact subintervals.
For $p,p'\in[0,1]$, define the intrinsic Wright--Fisher distance by
$
d_{\mathrm{WF}}(p,p'):=\left|\int_{p'}^p\frac{dr}{\sqrt{r(1-r)}}\right|.
$
Then $\Phi$ introduced in  Section \ref{sec:arcsin.transform} is an isometry from $[0,1]$ onto the Euclidean interval $[0,\pi]$.
The integrals are understood as convergent improper integrals when an endpoint is $0$ or $1$. For $0<\alpha<1$, $u\in\Cc([0,1])$ and $v\in\Cc([0,1])\cap\Cc^1((0,1))$, set
\bee
[u]_{\alpha;\mathrm{WF}}:=
\sup_{\substack{p,p'\in[0,1], p\neq p'}}\frac{|f(p)-f(p')|}{d_{\mathrm{WF}}(p,p')^\alpha},\qquad 
D_Kv(p):=\sqrt{p(1-p)}\,v_p(p),\quad \text{$p\in(0,1)$}.
\eee
We denote by $\Cc_K^{1,\alpha}([0,1])$ the space of such functions for which $D_Kv$ extends continuously to $[0,1]$ and the following norm is finite:
\begin{equation}\label{eq:CK.original}
\begin{aligned}
\|v\|_{\Cc_K^{1,\alpha}([0,1])}
:=\|v\|_{L^\infty([0,1])}+[v]_{\alpha;\mathrm{WF}}+\|D_Kv\|_{L^\infty([0,1])}+[D_Kv]_{\alpha;\mathrm{WF}}
=\|v(p(s))\|_{\Cc^{1,\alpha}([0,\pi])}.
\end{aligned}
\end{equation}
The last equality follows from
$d_{\mathrm{WF}}(p,p')=|\Phi(p)-\Phi(p')|$ and $(v(p(s)))_s=(D_Kv)(p(s))$. Hence,
\bee
\Cc_K^{1,\alpha}([0,1])=\left\{v\in\Cc([0,1]):v\circ p\in\Cc^{1,\alpha}([0,\pi])\right\}.
\eee
The distance $d_{\mathrm{WF}}$ gives the standard one-dimensional Wright--Fisher geometry; see \cite[Section 5.2.1]{EpsteinMazzeo2013}. Our norm measures H\"older continuity of both the value and its intrinsic derivative $D_Kv$ in this geometry. We retain the notation $\Cc_K^{1,\alpha}$ to distinguish this space from the hybrid space $\Cc_{\mathrm{WF}}^{1,\alpha}$ in \cite[Section 5.2.3]{EpsteinMazzeo2013}, which controls the ordinary first derivative. 
For $\vpi\in\Ac$, recall the drifts $P^\vpi,Q^\vpi$ from Section~\ref{sec:arcsin.transform}. Let us define the policy-dependent singular and full weights by
\be\label{eq:pi.weights}
m^\vpi(s):=\exp\left(\int_{\pi/2}^sP^\vpi(t)dt\right),\quad \widehat m^\vpi(s):=\exp\left(\int_{\pi/2}^sQ^\vpi(t)dt\right).
\ee
Under the norm
$
\|u\|_{\Sc(\vpi)}:=\big(\int_0^\pi (u^2+u_s^2)\widehat m^\vpi(s)\,ds\big)^{1/2},
$
$\Sc(\vpi)$ is defined as the completion of $\Cc^1([0,\pi])$.
We use the completion in view that the weight may be singular at the endpoints. Assumption~\ref{ass:b} and Lemma~\ref{lm:ratio} below give $\int_0^\pi\widehat m^\vpi(s)\,ds<\infty$, so this norm is finite on $\Cc^1([0,\pi])$. Throughout the paper, $0<\alpha<1$ denotes a fixed H\"older exponent.

The PDE \eqref{eq:pde.bessel} may not admit a unique solution, we will first define a family of solution satisfying a natural zero-flux boundary behavior and then focus our discussion within this regular branch.
\begin{Definition}\label{def:regular.linear}
We call $u$ a \emph{regular solution} to \eqref{eq:pde.bessel} if
$
u\in \Cc^1((0,\pi))\cap W^{2,k}_{\mathrm{loc}}((0,\pi))\cap L^k(m ds),
$
\eqref{eq:pde.bessel} holds a.e. on $(0,\pi)$, $m^q u_s$ has a locally absolutely continuous representative, and
\begin{equation}\label{eq:zero-flux0}
\lim_{s\to0+}m^q(s)u_s(s)=0,\qquad
\lim_{s\to\pi-}m^q(s)u_s(s)=0.
\end{equation}
\end{Definition}


\begin{Theorem}\label{thm:main.besselest}
 Assume
 \be\label{eq:theta.bounds}
 0<\thl\leq \Lambda_i(s) \leq \thu<\infty,\quad \text{a.e. on $(0,\pi)$ for }i=1,2,
 \ee
 where $\thl, \thu$ are the constants in \eqref{eq:main.boundaryrelaxed}.  Choose $k$ such that
 $k>\max\{1,2\thu\}$ and $0<\alpha<1-\frac{\max\{1,2\thu\}}{k}.$
 Then the following constants can be chosen uniformly over all measurable pairs $(\Lambda_2,\Lambda_1)$ satisfying the above bounds:
 \begin{enumerate}[label=(\roman*)]
  \item {\bf Uniform Bessel type embedding:} There exists a finite constant $C_E=C_E(k,\alpha,\thl,\thu)>0$ such that for every $u\in W^{2,k}(m ds)$,
  $$
  \norm{u}_{\Cc^{1,\alpha}([0,\pi])}\leq C_E\norm{u}_{W^{2,k}(m ds)}.
  $$
  \item {\bf Uniform Bessel type Sobolev estimate:} If $\norm q_{L^\infty((0,\pi))}\leq M$ and $u$ is a regular solution of \eqref{eq:pde.bessel}, then
  $$
  \norm{u}_{W^{2,k}(m ds)}\leq C_S\left(\norm u_{L^k(m ds)}+\norm f_{L^k(m ds)}\right),
  $$
  where $C_S=C_S(k,\rho,M,\thl,\thu)>0$ is a finite constant.
 \end{enumerate}
\end{Theorem}
The proof of Theorem \ref{thm:main.besselest} is divided into the following subsections.

\subsection{Uniform structure of the Bessel weights}


We first provide a uniform endpoint ratio estimates on the $m, m^q$ measures.

\begin{Lemma}\label{lm:ratio}
There exists some $C\geq1$, depending only on $\thl,\thu,M$, such that for every measurable pair $\Lambda_1,\Lambda_2$ satisfying \eqref{eq:theta.bounds} and every $q$ with $\|q\|_{L^\infty((0,\pi))}\leq M$,
\begin{align}\label{eq:lm.bessel.ratio}
\begin{cases}
C^{-1}\left(\frac ts\right)^{2\thu-1}\leq \frac{m^q(t)}{m^q(s)}\leq C\left(\frac ts\right)^{2\thl-1},\qquad\; 0<t\leq s\leq \frac{\pi}{8}, \\
C^{-1}\left(\frac ts\right)^{2\thu-1}\leq \frac{m^q(\pi-t)}{m^q(\pi-s)}\leq C\left(\frac ts\right)^{2\thl-1},\;\quad 0<t\leq s\leq \frac{\pi}{8}, \\
C^{-1}\leq m^q(s)\leq C,\qquad\qquad\qquad \qquad \quad  s\in[\frac{\pi}{8},\frac{7\pi}{8}]. 
\end{cases}
\end{align}
It also holds that
\begin{equation}\label{eq:mass-ratio}
 \begin{cases}
C^{-1}sm^q(s)\leq\int_0^s m^q(t) dt\leq Csm^q(s),\qquad \qquad \qquad \quad \; 0<s\leq \frac{\pi}{8},\\
C^{-1}sm^q(\pi-s)\leq\int_0^s m^q(\pi-t) dt\leq Csm^q(\pi-s),\qquad 0<s\leq \frac{\pi}{8}.
 \end{cases}
\end{equation}
After possibly changing $C$,
\begin{align}\label{eq:general-envelope}
 \begin{cases}
C^{-1}s^{2\thu-1}\leq  m^q(s)\leq Cs^{2\thl-1}, \qquad \qquad \qquad 0<s\leq \frac{\pi}{8},\\
C^{-1}(\pi-s)^{2\thu-1}\leq  m^q(s)\leq C(\pi-s)^{2\thl-1}, \quad \frac{7\pi}{8}\leq s<\pi,
 \end{cases}
\end{align}
Then, for every finite $k\geq1$,
\begin{equation}\label{eq:general-masslog}
\sup_{(\Lambda_2,\Lambda_1,q)}\int_0^\pi m^q(s) ds<\infty,
\qquad
\sup_{(\Lambda_2,\Lambda_1,q)}\left\|\ln\left(2+\frac1{\sin(\cdot)}\right)\right\|_{L^k(m^qds)}<\infty,
\end{equation}
where the suprema are over the uniformly controlled class. All above conclusions hold with $m^q$ replaced by $m$.
\end{Lemma}

\begin{proof}
We prove the assertions near $0$ first. For $0<s\leq \pi/8$,
$
\cot\frac{s}{2}=\frac2s+O(s)$ and $\tan\frac{s}{2}=O(s).
$
Hence \eqref{eq:Pdef}, \eqref{eq:theta.bounds}, and $\|q\|_\infty\leq M$ give
$$
P(s)+q(s)=\frac{2\Lambda_2(s)-1}{s}+R_0(s)\quad  \text{with}\;  |R_0(s)|\leq C_0,
$$
where $C_0$ is a finite constant that depends only on $\thl,\thu,M$. Thus, for $0<t\leq s\leq \pi/8$,
$$
(2\thl-1)\ln\frac{s}{t}-C_0(s-t)
\leq\int_t^s[P(r)+q(r)] dr
\leq(2\thu-1)\ln\frac{s}{t}+C_0(s-t).
$$
This, together with
$
\frac{m^q(t)}{m^q(s)}=\exp\left(-\int_t^s[P(r)+q(r)] dr\right),
$
proves the first line of \eqref{eq:lm.bessel.ratio}. Near the right endpoint, by setting $r=\pi-s$, we have
$\tan\frac{s}{2}=\frac2r+O(r)$ and $\cot\frac{s}{2}=O(r)$, then
$$
P(\pi-r)+q(\pi-r)=-\frac{2\Lambda_1(\pi-r)-1}{r}+R_1(r) \quad \text{with}\; |R_1(r)|\leq C_0,
$$
and the same integration argument yields the second line in \eqref{eq:lm.bessel.ratio}. On $[\pi/8,7\pi/8]$, $P+q$ is uniformly bounded. Since $m^q(\pi/2)=1$, integration from $\pi/2$ gives the third line in \eqref{eq:lm.bessel.ratio}. Taking $q=0$ proves the corresponding assertions for $m$. The estimate \eqref{eq:mass-ratio} follows from \eqref{eq:lm.bessel.ratio}. Indeed, since both powers $2\thl-1, 2\thu-1$ are integrable at $0$ by $2\thl-1>-1$, for $0<t\leq s\leq \pi/8$,
integrating the first line in $t$ of \eqref{eq:lm.bessel.ratio} gives the first line of \eqref{eq:mass-ratio}, the right-endpoint version follows from the second line of \eqref{eq:lm.bessel.ratio}. By the last line in \eqref{eq:lm.bessel.ratio}, $m^q(\pi/8)$ and $m^q(7\pi/8)$ are uniformly bounded above and below. Then \eqref{eq:general-envelope} follows by fixing $s=\pi/8$ in the first two lines of \eqref{eq:lm.bessel.ratio}.

As $2\thl-1>-1$, the upper bounds in \eqref{eq:general-envelope} imply the first formula in \eqref{eq:general-masslog}. Moreover, $\sin s=s+o(|s|)$ near $0$,
it then follows that
$$
\ln\left(2+\frac1{\sin s}\right)\leq C_1(1+|\ln s|)\quad  0<s\leq \frac{\pi}{8}
$$
Thus, by the first line in \eqref{eq:general-envelope}, for every $k\in[1,\infty)$,
$$
\int_0^{\frac{\pi}{8}}  \left| \ln(2+\frac{1}{\sin s}) \right|^k m^q(s)ds\leq C_2  \int_0^{\frac{\pi}{8}}(1+|\ln s|^k)s^{2\thl-1} ds<\infty.
$$
The right endpoint can be treated in the same fashion, while the logarithm is bounded on the middle interval. This concludes \eqref{eq:general-masslog}.

At last, following the same arguments by taking $q\equiv 0$ gives the result for $m$.
\end{proof}

\subsection{Proof of Theorem \ref{thm:main.besselest} Part (i)}

\begin{Lemma}\label{lm:inverse-weight}
 Let $k>\max\{1,2\thu\}$, set $\alpha_*:=1-\max\{1,2\thu\}/k>0$ and $k'=\frac{k}{k-1}$. There exists a finite constant $C>0$ depending only on $k,\thl,\thu$, such that for every interval $I\subset(0,\pi)$,
 \begin{equation}\label{eq:inverse-weight}
  \left(\int_I m(s)^{-1/(k-1)} ds\right)^{1/k'}\leq C|I|^{\alpha_*}.
 \end{equation}
\end{Lemma}

\begin{proof}
Utilizing \eqref{eq:general-envelope} with $q\equiv 0$ in Lemma~\ref{lm:ratio}, we get that
 $$
 m(s)\geq C^{-1}s^{2\thu-1} \text{ on  $(0, \frac{\pi}{8}]$ },\quad m(s)\geq C^{-1}(\pi-s)^{2\thu-1} \text{ on $[\frac{7\pi}{8}, \pi)$}.
 $$
 On $[\pi/8, 7\pi/8]$, $m$ is uniformly bounded below. If $2\thu-1\le0$, $\max\{1,2\thu\}=1$, and then the above inequalities tells that $m^{-1/(k-1)}$ is bounded near each endpoint, so \eqref{eq:inverse-weight} holds with exponent $1/k'=1-1/k=\alpha_*$. Suppose $2\thu-1>0$, it follows that on $(0, \pi/8]$,
 $$
 m(s)^{-1/(k-1)}\leq Cs^{-(2\thu-1)/(k-1)}.
 $$
 By $k>2\thu$, we have $0<(2\thu-1)/(k-1)<1$. Therefore, for every interval of length $h$ contained in $(0, \pi/8]$,
 $$
 \int_I m(s)^{-1/(k-1)} ds\leq \int_I Cs^{-(2\thu-1)/(k-1)}  ds \leq Ch^{1-(2\thu-1)/(k-1)}.
 $$
 Taking the $1/k'$ power gives
 $$
 h^{\frac{1-(2\thu-1)/(k-1)}{k'}}=h^{1-2\thu/k}=h^{\alpha_*}.
 $$
 In sum, the result is verified for the left-endpoint. The part for the right endpoint follows in a similar way. In the middle, by $\max\{1,2\thu\}\geq1$, we have the stronger bound $Ch^{1/k'}$ and $1/k'\geq\alpha_*$. Splitting a general interval into at most three pieces completes the proof.
\end{proof}

\begin{Proposition}\label{prop:sobolev.embed}
Suppose $\Lambda_1,\Lambda_2$ be arbitrary measurable functions satisfying \eqref{eq:theta.bounds}, and let $m$ be defined as in \eqref{eq:Pdef}. Let $k>\max\{1,2\thu\}$, set $\alpha_*:=1-\max\{1,2\thu\}/k>0$ and $k'=\frac{k}{k-1}$. For $0<\alpha<\alpha_*$, every $u\in W^{2,k}(m ds)$ has a unique representative in $\Cc^{1,\alpha}([0,\pi])$, and
 \begin{equation}\label{eq:embedding}
  \norm{u}_{\Cc^{1,\alpha}([0,\pi])}\leq C_E\norm{u}_{W^{2,k}(m ds)},
 \end{equation} 
 where $C_E$ depends only on $k,\alpha,\thl,\thu$.
\end{Proposition}

\begin{proof}
 By the definition of $W^{2,{k}}(m ds)$ in \eqref{eq:W2pdef}, $u_s$ admits a representative that is locally absolutely continuous on $(0,\pi)$, whose a.e. derivative is $u_{ss}$. We use this representative below. For $0<s<t<\pi$, the fundamental theorem of calculus and H\"older's inequality give
 \begin{align*}
  |u_s(t)-u_s(s)|\leq \int_s^t|u_{ss}(r)| dr
  \leq \left(\int_s^t|u_{ss}(r)|^{k}m(r) dr\right)^{1/{k}}\left(\int_s^t m(r)^{-1/({k}-1)} dr\right)^{1/{k'}}.
 \end{align*}
 This, together with Lemma~\ref{lm:inverse-weight}, gives
 \be\label{eq:embed.interior}
 |u_s(t)-u_s(s)|\leq C\|u_{ss}\|_{L^{k}(m ds)}|t-s|^{\alpha_*},\quad \text{for }0<s<t<\pi.
 \ee

For any sequence $s_j\to 0+$, \eqref{eq:embed.interior} shows that $(u_s(s_j))_{j\geq1}$ is Cauchy. The same argument applies to sequences tending to $\pi$. In particular, the limits $g(0):=\lim_{s\to0+}u_s(s)$ and $g(\pi):=\lim_{s\to\pi-}u_s(s)$ exist. Define $g(s):=u_s(s)$ for $s\in(0,\pi)$. Letting $s\to0+$ and $t\to\pi-$ in \eqref{eq:embed.interior} respectively shows that $g\in\Cc^{\alpha_*}([0,\pi])$ and
 \be\label{eq:embed.0}
 [g]_{\Cc^{\alpha_*}([0,\pi])}\leq C\|u_{ss}\|_{L^{k}(m ds)}.
 \ee
In view that the locally absolutely continuous representative of $u$ satisfies $u(s)-u(\pi/8)=\int_{\pi/8}^s u_s(t)\,dt$ on $(0,\pi)$, the function $\widetilde u(s):=u(\pi/8)+\int_{\pi/8}^s g(t) dt $ defined on $[0,\pi]$ agrees with $u$ almost everywhere on $(0,\pi)$, and belongs to $\Cc^1([0,\pi])$. It thus has derivative $\widetilde u_s=g$.
 For $0<\alpha<\alpha_*$, it follows from \eqref{eq:embed.0} that $[\widetilde u_s]_{\Cc^\alpha([0,\pi])} \leq\pi^{\alpha_*-\alpha}C\|u_{ss}\|_{L^{k}(m ds)}$.

We now focus on the middle interval $\left[\frac{\pi}{3},\frac{2\pi}{3}\right].$ By the last line of \eqref{eq:lm.bessel.ratio} in Lemma~\ref{lm:ratio} with $q\equiv 0$, there exist constants $0<c_0<c_1<\infty$, independent of the choices of functions $\Lambda_1,\Lambda_2$, such that
 $
 c_0\leq m(s)\leq c_1$ for  all $s\in \left[\frac{\pi}{3},\frac{2\pi}{3}\right].
 $
 Hence,
 $$
 \int_{[\frac{\pi}{3},\frac{2\pi}{3}]} |u_s(s)|^{k} ds \leq c_0^{-1}\int_{[\frac{\pi}{3},\frac{2\pi}{3}]} |u_s(s)|^{k} m(s)ds \leq c_0^{-1}\norm{u_s}_{L^{k}(m ds)}^{k}.
 $$
As a result, there exists a point $s_1\in [\pi/3,2\pi/3]$ such that
 $$
 |u_s(s_1)|^{k} \leq \frac{3}{\pi}\int_{[\frac{\pi}{3}, \frac{2\pi}{3}]} |u_s(s)|^{k} ds.
 $$
 Suppose not, $|u_s(s)|^{k}$ would be strictly larger than the right-hand side for every $s\in [\pi/3,2\pi/3]$, which yields a contradiction after integration over $[\pi/3,2\pi/3]$. Hence,
 $|u_s(s_1)|\leq C\norm{u_s}_{L^{k}(m ds)}.$
 For an arbitrary $s\in(0,\pi)$, we have
 $$
 |u_s(s)|\leq |u_s(s_1)|+|u_s(s)-u_s(s_1)|
 \leq C\norm{u_s}_{L^{k}(m ds)}+[u_s]_{C^{\alpha_*}((0,\pi))}|s-s_1|^{\alpha_*}.
 $$
 This, together with \eqref{eq:embed.0}, implies that
\be\label{eq:embed.1}
 \norm{u_s}_{L^\infty([0,\pi])}\leq C\left(\norm{u_s}_{L^{k}(m ds)}+\norm{u_{ss}}_{L^{k}(m ds)}\right).
\ee
 The same argument gives a pointwise bound for $u$. Indeed, as $m$ is uniformly bounded below on $[\pi/3,2\pi/3]$,
 $
 \int_{[\frac{\pi}{3},\frac{2\pi}{3}]} |u(s)|^{k} ds\leq C\norm{u}_{L^{k}(m ds)}^{k},
 $
 and hence there exists $s_2\in [\pi/3,2\pi/3]$ such that
 $
 |u(s_2)|\leq C\norm{u}_{L^{k}(m ds)}.
 $
 For any $s\in(0,\pi)$, the fundamental theorem of calculus gives
 $$
 |u(s)|\leq |u(s_2)|+\left|\int_{s_2}^s u_r(r) dr\right| \leq C\norm{u}_{L^{k}(m ds)}+|s-s_2|\norm{u_s}_{L^\infty([0,\pi])},
 $$
 which leads to
\be\label{eq:embed.2}
 \norm{u}_{L^\infty([0,\pi])} \leq C\left(\norm{u}_{L^{k}(m ds)}+\norm{u_s}_{L^\infty([0,\pi])}\right).
\ee
Moreover, for arbitrary $s,t\in(0,\pi)$
and $\alpha\in(0,1)$, it holds that
 $$
 \frac{|u(t)-u(s)|}{|t-s|^\alpha} \leq\norm{u_s}_{L^\infty([0,\pi])}|t-s|^{1-\alpha} \leq\pi^{1-\alpha}\norm{u_s}_{L^\infty([0,\pi])},
 $$
 and thus
 $
 [u]_{C^\alpha([0,\pi])}\leq\pi^{1-\alpha}\norm{u_s}_{L^\infty([0,\pi])}.
 $
 Combining this with \eqref{eq:embed.0} -- \eqref{eq:embed.2} concludes \eqref{eq:embedding}.
\end{proof}

\subsection{Proof of Theorem \ref{thm:main.besselest} Part (ii)}

Through the rest of this subsection, we take $1<{k}<\infty$ and
$
{k'}=\frac{{k}}{{k}-1}.
$
Let us first recall the classical two-weight Hardy criterion in the following lemma (see, e.g., \cite[Section~1.1]{KufnerPersson2003}). 

\begin{Lemma}\label{lm:hardy.prepare}
 Let $-\infty<a<b<\infty$, and  take measurable functions $V,W:(a,b)\to(0,\infty)$.

 \begin{enumerate}
  \item[(i)] The inequality
  \begin{equation}\label{eq:lmhardy1}
   \int_a^b \left|\int_a^s F(t),dt \right|^{k} V(s)ds
   \leq C\int_a^b |F(s)|^{k}W(s) ds
  \end{equation}
  holds for every measurable $F$ for which the right-hand side is finite if
  \begin{equation}\label{eq:lmhardy0}
   \sup_{a<r<b}\left(\int_r^b V(s)ds\right)^{1/{k}}\left(\int_a^r W(s)^{-1/({k}-1)}ds\right)^{1/{k'}} <\infty.
  \end{equation}
  Moreover, the constant $C$ in \eqref{eq:lmhardy1} depends only on ${k}$ and the quantity in \eqref{eq:lmhardy0}.

  \item[(ii)] Similarly, the inequality
  \begin{equation}\label{eq:lmhardy3}
   \int_a^b\left| \int_s^b F(t),dt \right|^{k}V(s)ds \leq C \int_a^b |F(s)|^{k}W(sds
  \end{equation}
  holds for every measurable $F$ for which the right-hand side is finite if
  \begin{equation}\label{eq:lmhardy2}
   \sup_{a<r<b}\left(\int_a^r V(s)ds\right)^{1/{k}} \left(\int_r^b W(s)^{-1/({k}-1)}ds\right)^{1/{k'}} <\infty.
  \end{equation}
  Again, the constant $C$ in \eqref{eq:lmhardy3} depends only on ${k}$ and the quantity in \eqref{eq:lmhardy2}.
 \end{enumerate}
\end{Lemma}

\begin{Lemma}\label{lm:hardy}
There exists $C_H<\infty$, depending only on ${k},\thl,\thu,M$, such that if $u\in \Cc^1((0,\pi))$, $m^q u_s\in W^{1,1}_{\mathrm{loc}}((0,\pi))$, \eqref{eq:zero-flux0} holds, and 
$g(s):=\frac{1}{m^q(s)}(m^q(s)u_s(s))_s\in L^{k}(m^q ds)$, then
\begin{equation}\label{eq:lm.hardy}
\norm{\frac{u_s}{\min\{s,\pi-s\}}}_{L^{k}(m^q ds)}\leq C_H\norm{g}_{L^{k}(m^q ds)}.
\end{equation}
As a result, with $C=\frac{\pi}{2}C_H$,
\begin{equation}\label{eq:first-derivative}
\norm{u_s}_{L^{k}(m^q ds)}\leq C\norm{g}_{L^{k}(m^q ds)}.
\end{equation}
\end{Lemma}

\begin{proof}
 Throughout the proof, $C$ is a generic positive constant that only depends on ${k},\thl,\thu,M$ and may change from line to line. We first treat the left endpoint $0$. By Lemma \ref{lm:ratio}, the total mass of $\it_0^\pi m^q(s) ds$ is uniformly finite. Indeed, \eqref{eq:mass-ratio} and the third line in \eqref{eq:lm.bessel.ratio} give
 $ \int_0^\pi m^q(s) ds\leq C.$
Then $g\in L^{k}(m^q ds)$ and H\"older's inequality imply that
 $$
 \int_0^\pi |g(s)|m^q(s) ds \leq\norm{g}_{L^{k}(m^q ds)} \left(\int_0^\pi m^q(s) ds\right)^{1/{k'}} <\infty.
 $$
 Hence $(m^q u_s)_s=m^q g\in L^1((0,\pi))$, and $m^q u_s$ has an absolutely continuous representative. It follows from integration from $\varepsilon$ to $s$ that
 $$
 m^q(s)u_s(s)-m^q(\varepsilon)u_s(\eps)
 = \int_\varepsilon^s m^q(t)g(t) dt.
 $$
 Sending $\varepsilon\to0+$ and applying the left zero-flux condition in \eqref{eq:zero-flux0}, we obtain
 \begin{equation}\label{eq:flux-representation}
 u_s(s) =\frac1{m^q(s)}\int_0^s m^q(t)g(t) dt.
 \end{equation}
  Thus, to show  \eqref{eq:lm.hardy} near $0$, it suffices to prove
 $$
 \int_0^{\frac{\pi}{8}}\left| \frac1{s m^q(s)}\int_0^s m^q(t)g(t) dt \right|^{k} m^q(s) ds
 \leq C\int_0^{\frac{\pi}{8}}|g(s)|^{k}m^q(s) ds.
 $$
 We next apply the Hardy criterion in Lemma \ref{lm:hardy.prepare} Part (i) with $a=0$, $b=\pi/8$, $r\in(0,\pi/8)$ and
 $$
 F(t)=m^q(t)g(t),\quad  V(s)=s^{-{k}}(m^q(s))^{1-{k}},\quad W(s)=(m^q(s))^{1-{k}}.
 $$
 As $(W(s))^{-1/({k}-1)}=m^q(s),$ the relevant Hardy quantity is
 \begin{equation}\label{eq:hardy-quantity}
  \sup_{0<r<\frac{\pi}{8}}\left(\int_r^{\frac{\pi}{8}}s^{-{k}}(m^q(s))^{1-{k}} ds \right)^{1/{k}}  \left(\int_0^r m^q(s) ds\right)^{1/{k'}}.
 \end{equation}
 For $0<r\leq s\leq \pi/8$, the first line of \eqref{eq:lm.bessel.ratio} in Lemma~\ref{lm:ratio} gives
 $
 m^q(s) \geq C^{-1}m^q(r) \left(\frac{s}{r}\right)^{2\thl-1}.
 $
 As $1-{k}<0$, it holds that
 $
 (m^q(s))^{1-{k}} \leq C (m^q(r))^{1-{k}} \left(\frac{s}{r}\right)^{-(2\thl-1)({k}-1)}.
 $
 Consequently,
 \begin{align*}
  \int_r^{\frac{\pi}{8}}s^{-{k}}(m^q(s))^{1-{k}} ds&\leq C (m^q(r))^{1-{k}}r^{(2\thl-1)({k}-1)}\int_r^{\frac{\pi}{8}}s^{-{k}-(2\thl-1)({k}-1)} ds.
 \end{align*}
By ${k}+(2\thl-1)({k}-1) =1+2\thl({k}-1)>1,$ we have
 \begin{align*}
  \int_r^{\frac{\pi}{8}}s^{-{k}-(2\thl-1)({k}-1)}ds&=\frac{r^{1-{k}-(2\thl-1)({k}-1)}-\left(\pi/8\right)^{1-{k}-(2\thl-1)({k}-1)}}{{k}+(2\thl-1)({k}-1)-1 }
\leq C r^{1-{k}-(2\thl-1)({k}-1)},
 \end{align*}
which gives
 \begin{equation}\label{eq:lm.hardy.1}
  \int_r^{\frac{\pi}{8}}s^{-{k}}(m^q(s))^{1-{k}} ds \leq C (m^q(r))^{1-{k}}r^{1-{k}}.
 \end{equation}
 Meanwhile, \eqref{eq:mass-ratio} in Lemma~\ref{lm:ratio} gives
$
  \int_0^r m^q(s) ds \leq C r m^q(r).
$
Substituting this and \eqref{eq:lm.hardy.1} into \eqref{eq:hardy-quantity} and noting that  $\frac{1-{k}}{{k}}+\frac1{{k'}}=0$, we obtain
 \begin{align*}
 & \left(  \int_r^{\frac{\pi}{8}}s^{-{k}}m^q(s)^{1-{k}} ds  \right)^{1/{k}}  \left(  \int_0^r m^q(s) ds \right)^{1/{k'}} \\
  & \leq  C  \left[  m^q(r)^{1-{k}}r^{1-{k}} \right]^{1/{k}} \left[ r m^q(r) \right]^{1/{k'}}
  =  C(m^q(r))^{\frac{1-{k}}{{k}}+\frac1{{k'}}} r^{\frac{1-{k}}{{k}}+\frac1{{k'}}}=C.
 \end{align*}
 Invoking Lemma \ref{lm:hardy.prepare} (i) with the help of \eqref{eq:flux-representation} yields
 \begin{equation}\label{eq:hardy-left}
  \norm{\frac{u_s}{s}}_{L^{k}((0,\frac{\pi}{8}),m^q ds)}  \leq  C\norm{g}_{L^{k}((0, \frac{\pi}{8}),m^q ds)}.
 \end{equation}
 The argument for the right endpoint follows a similar fashion and
 we obtain
 \begin{equation}\label{eq:hardy-right}
  \norm{\frac{u_s}{\pi-s}}_{L^{k}((\frac{\pi}{8},\pi),m^q ds)}  \leq C\norm{g}_{L^{k}((\frac{7\pi}{8},\pi),m^q ds)}.
 \end{equation}

 We next consider the middle interval $[\pi/8,7\pi/8]$. Thanks to the third line of \eqref{eq:lm.bessel.ratio} in Lemma~\ref{lm:ratio},
 $C^{-1}\leq m^q(s)\leq C$ for all $s\in[\pi/8,7\pi/8].$
Then by \eqref{eq:flux-representation}, H\"older's inequality and $\int_0^\pi m^q(s) ds\leq C$, for every $s\in[\pi/8,7\pi/8]$, it holds that
 \begin{align*}
  |u_s(s)|  &\leq \frac1{m^q(s)}  \int_0^s m^q(t)|g(t)| dt
  \leq  C \int_0^\pi m^q(t)|g(t)| dt\\
  &\leq  C\norm{g}_{L^{k}(m^q ds)}  \left(  \int_0^\pi m^q(t) dt  \right)^{1/{k'}} 
   \leq  C\norm{g}_{L^{k}(m^q ds)}.
 \end{align*}
 Consequently, using the uniform boundedness of $m^q$ on the middle interval again, we have
 \begin{equation}\label{eq:hardy-middle}
  \norm{u_s}_{L^{k}((\frac{\pi}{8},\frac{7\pi}{8}),m^q ds)} \leq C\norm{g}_{L^{k}(m^q ds)}.
 \end{equation}
In view of $\min\{s,\pi-s\}\geq \frac{\pi}{8}$ on $[\pi/8,7\pi/8]$, \eqref{eq:hardy-middle} implies $\norm{u_s/\min\{s,\pi-s\}}_{L^{k}((\pi/8, 7\pi/8),m^q ds)}\leq C\norm{g}_{L^{k}(m^q ds)}$. It results from a combination of \eqref{eq:hardy-left}, \eqref{eq:hardy-right} and \eqref{eq:hardy-middle} that \eqref{eq:lm.hardy} holds. Then, \eqref{eq:first-derivative} readily follows from \eqref {eq:lm.hardy} and the fact that $\min\{s,\pi-s\}\leq\frac{\pi}{2}$.
\end{proof}

\begin{Proposition}\label{prop.sobolev.est}
 Let $\rho>0$ and $P,m$ be generated by arbitrary measurable $\Lambda_1,\Lambda_2$ satisfying \eqref{eq:theta.bounds}. Suppose $\norm q_{L^\infty((0,\pi))}\leq M$ and $f\in L^{k}(m ds)$. Let $u$ be a regular solution of
 \begin{equation}\label{eq:resolvent}
  -\rho u+\Bc u+q(s)u_s=f
 \end{equation}
 satisfying the zero-flux condition \eqref{eq:zero-flux0}. Then
$$
  \norm{u}_{W^{2,{k}}(m ds)}\leq C_S\left(\norm{u}_{L^{k}(m ds)}+\norm{f}_{L^{k}(m ds)}\right),
$$
 where $C_S>0$ is a finite constant depending only on ${k},\rho,M,\thl,\thu$.
\end{Proposition}

\begin{proof}
 By \eqref{eq:P.transform}, \eqref{eq:resolvent} becomes
$
  \frac1{m^q}(m^q u_s)_s=\rho u+f.
$
 Then Lemma~\ref{lm:hardy} with $g= \rho u+f$ gives
 \begin{equation}\label{eq:lm.est.us}
  \norm{u_s}_{L^{k}(m^q ds)}+\norm{\frac{u_s}{\min\{\cdot, \pi-\cdot \}}}_{L^{k}(m^q ds)}\leq C \|\rho u+f \|_{L^{k}(m^q ds)}.
 \end{equation}
 Noting $\norm q_{L^\infty}\leq M$ and \eqref{eq:Pdef}, we have that
$|P(s)+q(s)|\leq \frac{C}{\min\{s,\pi-s\}}$ for all $s\in(0,\pi)$
 with $C$ depending only on $\thl,\thu,M$. Then, by $u_{ss}=\rho u+f-(P+q)u_s$, we have
 $$
 \norm{u_{ss}}_{L^{k}(m^q ds)}\leq \|\rho u+f \|_{L^{k}(m^q ds)}+C\norm{\frac{u_s}{\min\{\cdot,\pi-\cdot\}}}_{L^{k}(m^q ds)}\leq C\|\rho u+f\|_{L^{k}(m^q ds)}.
 $$
The above inequality with \eqref{eq:lm.est.us} implies that
 $$
\|u\|_{W^{2,{k}}(m^q ds)}\leq C\left(\rho\norm u_{L^{k}(m^q ds)}+\norm f_{L^{k}(m^q ds)}\right).
 $$
 The desired  result follows from \eqref{eq:measure-equivalence}.
\end{proof}


\section{Policy Iteration and Policy-dependent Bessel Weights}\label{sec:PIA}
Throughout Sections \ref{sec:PIA}--\ref{sec:convergence.entropy}, Assumptions~\ref{ass:b} and~\ref{ass:action} are in force unless explicitly stated otherwise. Results that hold under weaker conditions are specified in their respective statements.

Fix $\lambda>0$, the PIA is carried out over $\Ac_\lambda$. Recall the regularized value $V_\lambda^{\vpi}$ in Section 1 and the transformed value $W^{\vpi}_\lambda(s)=V_\lambda^{\vpi}(p(s))$ in Section 2. Under the arcsine transformation $S_t^{\vpi}=\Phi(P_t^{\vpi})$, applying It\^o's formula up to the first exit from any compact subinterval of $(0,\pi)$ gives
\be\label{eq:transformed.local.SDE}
dS_t^{\vpi}=Q^{\vpi}(S_t^{\vpi})\,dt+\sqrt{2}\,dB_t
\quad\text{before hitting $\{0,\pi\}$}.
\ee
For a trajectory that reaches an endpoint, \eqref{eq:transformed.local.SDE} is not a global SDE on $[0,\pi]$. Its extension at the endpoints is fixed by Lemma \ref{lm:closed.diffusion}. 

Define the two Gibbs maps
\be\label{eq:def.Gamma}
\Gamma_{\lambda,i}(s,y,a):=\frac{\exp\left(\lambda^{-1}[Q_i(s,a)y+\widetilde r_i(s,a)]\right)}{\int_A\exp\left(\lambda^{-1}[Q_i(s,a')y+\widetilde r_i(s,a')]\right)da'},\quad i=1,2.
\ee
We consider the following {\it Policy Iteration Algorithm}:\\
{\bf Initialization:} Set $\vpi_1^0\equiv\vpi_2^0\equiv\nu$ with $\nu(a)=1/\operatorname{Leb}(A)$ and let $u^0(s)=W^{\vpi^0}_\lambda(s)$ on $[0,\pi]$.\\
{\bf Iteration:} For each $n\geq1$, take
$$
\vpi_i^n(s,a):=\Gamma_{\lambda,i}(s,u_s^{n-1}(s),a)\;\text{for } i=1,2, \quad \text{and compute $u^n(s):=W^{\vpi^n}_\lambda(s)$ on $[0,\pi]$.}
$$
Accordingly, the PIA under the original coordinate starts from $v^0=V_\lambda^{\vpi^0}$ with the same uniform policy. For $n\geq1$, the iteration rules are given by, for $p\in(0,1)$,
\bee
\begin{aligned}
\vpi_1^n(p,a)
=\frac{e^{\frac1\lambda[-2q^{1,2}(p,a)v_p^{n-1}(p)+r(1,p,a)]}}
{\int_Ae^{\frac1\lambda[-2q^{1,2}(p,a')v_p^{n-1}(p)+r(1,p,a')]da'}},\quad
\vpi_2^n(p,a)
=\frac{e^{\frac1\lambda[2q^{2,1}(p,a)v_p^{n-1}(p)+r(2,p,a)]}}
{\displaystyle\int_Ae^{\frac1\lambda[2q^{2,1}(p,a')v_p^{n-1}(p)+r(2,p,a')]da'}}.
\end{aligned}
\eee
The policy evaluation is $v^n:=V_\lambda^{\vpi^n}$ on $[0,1]$. Note that these are two coordinate representations of the same algorithm. Indeed, for $s=\Phi(p)$ and $y=\sqrt{p(1-p)}z$, \eqref{eq:def.Q} gives
$Q_1(s,a)y=-2q^{1,2}(p,a)z+z/2$ and $Q_2(s,a)y=2q^{2,1}(p,a)z-z/2$.
The action-independent terms $\pm z/2$ cancel in the normalized Gibbs densities. Together with $W_\lambda^{\vpi}(s)=V_\lambda^{\vpi}(p(s))$, this identifies both the update and evaluation steps and yields, inductively whenever the iterates are well-defined, that
\begin{equation}\label{eq:PI.coordinates}
u^n(s)=v^n(p(s)),\quad s\in[0,\pi],\qquad
u_s^n(s)=\sqrt{p(s)(1-p(s))}\,v_p^n(p(s))\quad \forall s\in(0,\pi).
\end{equation}
We will establish well-posedness, improvement and convergence in the transformed coordinate below, \eqref{eq:PI.coordinates} transfers these conclusions to the original MFC problem.

The Gibbs formula \eqref{eq:def.Gamma} is used for $s\in(0,\pi)$. At the boundary $s=0,\pi$, the distribution of the policy can be chosen arbitrarily in $\Pc(A)$. The same choices are used at the corresponding original endpoints $p=0,1$. This does not affect the state law or the expected payoff because \eqref{eq:boundary.zero.occupation} ensures zero occupation time at the two boundary points on every finite horizon in Lemma \ref{lm:closed.diffusion}. The well-posedness of both coordinate formulations is provided in Corollary~\ref{cor:PIA.welldefined}.

\subsection{Policy-dependent Bessel structure and uniform endpoint control}
The transformed drift depends on the policy pair through its averaged endpoint weights. It is useful to record the corresponding Bessel structure once for an arbitrary admissible policy pair. All estimates below are uniform over the policy pair and therefore apply, in particular, to every iteration in PIA.

Given $\vpi=(\vpi_1,\vpi_2)\in\Ac$, Assumption~\ref{ass:b} gives
\be\label{eq:Theta.bounds}
\thl\leq\Lambda_i^\vpi(s)\leq\thu,\quad i=1,2.
\ee
Recall the policy-dependent weights $m^\vpi$ and $\widehat m^\vpi$ defined in \eqref{eq:pi.weights}.
Then by \eqref{eq:q.bound},
\be\label{eq:m.equivalence}
e^{-M\pi}m^\vpi(s)\leq\widehat m^\vpi(s)\leq e^{M\pi}m^\vpi(s),\quad s\in(0,\pi),
\ee
and, in the distributional sense, for any $u\in W^{2,k}(m^\theta ds)$ that
$$
u_{ss}+Q^\vpi u_s=\frac1{\widehat m^\vpi}(\widehat m^\vpi u_s)_s.
$$
Thanks to \eqref{eq:Theta.bounds} and \eqref{eq:q.bound}, for any $\theta\in \Ac$, applying Lemma \ref{lm:ratio} with $\Lambda_i=\Lambda_i^\vpi$ and $q=q^\vpi$, we obtain that \eqref{eq:lm.bessel.ratio}--\eqref{eq:general-masslog} hold uniformly over $\vpi\in\Ac$ for $m^\vpi$ and $\widehat m^\vpi$.

For $n\geq1$, set
$$
\Bc^n u:=u_{ss}+P^{\vpi^n}(s)u_s=\frac1{m^{\theta^n}(s)}\frac{d}{ds}(m^{\theta^n}(s)u_s).
$$
Given the well-posedness of PIA, $u^n = W^{\theta^n}_\lambda$ shall satisfies
\be\label{eq:un.PDE}
\begin{aligned}
0=&-\rho u^n+\Bc^n u^n+q^{\vpi^n}(s)u_s^n\\
&+ \sin^2(\frac s2)\left[ \rtild^{\theta^n_1}_1(s)+\lambda \Hc(\vpi^n_1(s))  \right]+\cos^2(\frac s2) \left[ \rtild^{\theta^n_2}_2(s)+\lambda \Hc(\vpi^n_2(s))  \right].
\end{aligned}
\ee
In light of Assumption \ref{ass:b}, $\|q^{\vpi^n}\|_{L^\infty((0,\pi))}\leq M$ for every $n$ as well as \eqref{eq:general-masslog} in Lemma \ref{lm:ratio}, we deduce that
\be\label{eq:mass.uniform}
\sup_{n\geq1}\int_0^\pi m^{\theta^n}(s)ds<\infty.
\ee

We have the next result of the sub-linear growth estimate on the entropy, which generalizes the estimate in \cite[Corollary 4.14]{HuangWangZhou2025}.

\begin{Lemma}[Sublinear estimate on entropy]\label{lm:entropy.est}
Suppose Assumption~\ref{ass:action} holds. Then there exists a finite constant $C>0$ independent of $\lambda\in(0,\infty)$, $s\in(0,\pi)$ and $y\in\mathbb R$ such that
\be\label{eq:entropy-pointwise}
\left|\Hc\left(\Gamma_{\lambda,i}(s,y,\cdot)\right)\right|\leq C\left[1+|\ln \lambda|+ \ln(1+|y|)+|\ln(\sin s)|\right] \quad \text{for $i=1,2$}.
\ee
\end{Lemma}

\begin{proof}
Recall \eqref{eq:def.Q} and \eqref{eq:def.Gamma}, and for $i=1,2$, define
\be\label{eq:lm.Gidef}
G_i(s,y,a):=Q_i(s,a)y+\widetilde r_i(s,a),
\ee
and it follows that
$
\Gamma_{\lambda,i}(s,y,a)
=\frac{\exp(\lambda^{-1}G_i(s,y,a))} {\int_A\exp(\lambda^{-1}G_i(s,y,a')) da'}.
$
Fix $a,\bar a\in A$. We first estimate the action dependence of $G_1$. By \eqref{eq:def.Q},
\begin{align}\label{eq:entropy.est1}
|G_1(s,y,a)-G_1(s,y,\bar a)|\leq& \frac{2|y|}{\sqrt{p(s)(1-p(s))}} |q^{1,2}(p(s), a)-q^{1,2}(p(s), \bar a)| +L_r|a-\bar a|.
\end{align}
Recall $q_1(s,a)= -2\sqrt{\frac{p(s)}{1-p(s)}}[q^{1,2}(p(s),a)-q^{1,2}(1,a)]$. It the holds that
\begin{align*}
q^{1,2}(p, a)-q^{1,2}(p, \bar a) =&q^{1,2}(1, a)-q^{1,2}(1, \bar a) +\big[q^{1,2}(p, a)-q^{1,2}(1, a)\big]-\big[q^{1,2}(p, \bar a)-q^{1,2}(1, \bar a)\big]\\
=& \left[ q^{1,2}(1, a)-q^{1,2}(1, \bar a) \right] -\frac{1}{2}\sqrt{\frac{1-p}{p}} \left[q_1(s, a) -q_1(s, \bar a) \right],
\end{align*}
It follows from Assumption~\ref{ass:action} (i)--(ii) that
\be\label{eq:entropy.est3}
|q^{1,2}(p, a)-q^{1,2}(p, \bar a)|\leq\frac12 \left(L_b+L_q\sqrt{\frac{1-p}{p}}\right)|a-\bar a|.
\ee
Combining \eqref{eq:entropy.est1} with \eqref{eq:entropy.est3} yields
\be\label{eq:entropy.est.G1}
 |G_1(s,y,a)-G_1(s,y,\bar a)|\leq \left[   L_r+|y| \left( \frac{L_b}{\sqrt{p(s)(1-p(s))}}+\frac{L_q}{p(s)}\right) \right]   |a-\bar a|.
\ee
 A similar calculation gives
\be\label{eq:entropy.est.G2}
  |G_2(s,y,a)-G_2(s,y,\bar a)| \leq  \left[   L_r+|y| \left( \frac{L_b}{\sqrt{p(s)(1-p(s))}}+\frac{L_q}{1-p(s)}\right) \right]   |a-\bar a|,
\ee
 Noting $p=\sin^2(s/2)$ and $\sqrt{p(1-p)}=\frac12\sin s$, we have $0<\frac1p\leq\frac1{p(1-p)}=\frac4{\sin^2s}$ and $0<\frac1{1-p}\leq\frac1{p(1-p)}=\frac4{\sin^2s}$.
 
As a result,
\be\label{eq:lm.entropy.est0}
\max\left\{
L_r+|y|\left(\frac{L_b}{\sqrt{p(1-p)}}+\frac{L_q}{p}\right),\;
L_r+|y|\left(\frac{L_b}{\sqrt{p(1-p)}}+\frac{L_q}{1-p}\right)
\right\}
\leq C_0\left(1+\frac{|y|}{\sin^2s}\right),
\ee
where $C_0\geq1$ is independent of $\lambda$. Fix an arbitrary $(s,y)\in (0,\pi)\times \R$, and take
set
\bee
L(s,y):=C_0\left(1+\frac{|y|}{\sin^2s}\right),
\qquad
r_0:=\min\left\{\frac\gamma2,\frac\lambda{C_0\left(1+\frac{|y|}{\sin^2s}\right)}\right\}.
\eee
For $i\in\{1,2\}$ and $\bar a\in A$, by \eqref{eq:entropy.est.G1}--\eqref{eq:lm.entropy.est0}, we arrive at
\bee
|G_i(s,y,a)-G_i(s,y,\bar a)|\leq C_0\left(1+\frac{|y|}{\sin^2s}\right) |a-\bar a| \quad \forall a\in A.
\eee
By $0<r_0<\gamma$ and $ C_0\left(1+\frac{|y|}{\sin^2s}\right) r_0\leq\lambda$, the Gibbs formula \eqref{eq:def.Gamma} and Assumption \ref{ass:action}(v), it holds that
\bee
\begin{aligned}
\Gamma_{\lambda,i}(s,y,\bar a)
&=\left[\int_A\exp\left(\frac{G_i(s,y,a)-G_i(s,y,\bar a)}\lambda\right)da\right]^{-1}\\
&\leq\left[\int_{A\cap B_{r_0}(\bar a)}\exp\left(-\frac{C_0\left(1+|y|/\sin^2s\right) |a-\bar a|}\lambda\right)da\right]^{-1}\\
&\leq\frac{e}{\iota\Leb(B_1(0)) r_0^\ell}
=\frac{e}{\iota\Leb(B_1(0))}
\max\left\{\left(\frac2\gamma\right)^\ell,
\left(\frac{C_0\left(1+|y|/\sin^2s\right)}\lambda\right)^\ell\right\}.
\end{aligned}
\eee
Then taking logarithms and combining with $C_0\left(1+|y|/\sin^2s\right)\geq1$ leads to
\bee
\ln\Gamma_{\lambda,i}(s,y,\bar a)
\leq C+\ell|\ln\lambda|+\ell\ln\left(1+\frac{|y|}{\sin^2s}\right),
\eee
where $C$ depends only on $C_0,\gamma,\iota,\ell$ and is independent of $\lambda$. This, together with
$1+\frac{|y|}{\sin^2s}\leq\frac{1+|y|}{\sin^2s}$, ensures that
\bee
\ln\Gamma_{\lambda,i}(s,y,\bar a)
\leq C+\ell|\ln\lambda|+\ell\ln(1+|y|)+2\ell|\ln(\sin s)|\quad \forall \bar a \in A.
\eee
Finally, as $\Gamma_{\lambda,i}(s,y,\cdot)$ is a probability density, it holds that $\int_A\Gamma_{\lambda,i}\ln\Gamma_{\lambda,i}da$ is bounded from above by the same quantity on the right-hand side above. Meanwhile,
$$
\int_A \ln\Gamma_{\lambda,i}(s,y,a)  \Gamma_{\lambda,i}(s,y,a)  da \geq -\ln(\Leb(A)).
$$
Hence \eqref{eq:entropy-pointwise} holds for both $i=1,2$.
\end{proof}

\subsection{Well-posedness of PIA}\label{subsec:wellposed}

For an admissible policy $\vpi\in\Ac_\lambda$, set
\begin{equation}\label{eq:Fpi}
F^{\vpi}(s):=p(s)\big(\rtild_1^{\vpi_1}(s)+\lambda\Hc(\vpi_1(s))\big)+(1-p(s))\big(\rtild_2^{\vpi_2}(s)+\lambda\Hc(\vpi_2(s))\big).
\end{equation}
Then the policy evaluation PDE becomes
\begin{equation}\label{eq:PDEpi}
-\rho u+u_{ss}+Q^{\vpi}u_s+F^{\vpi}=0,
\quad\text{or equivalently,}\quad
\rho u-\frac1{\widehat m^{\vpi}}(\widehat m^{\vpi}u_s)_s=F^{\vpi}.
\end{equation}

The following lemma constructs the transformed controlled diffusion on $[0,\pi]$, identifies its symmetric Dirichlet form and specifies its boundary behavior. Its proof is conducted with the help of the classical scale--speed theory of one-dimensional diffusions (see \cite[Chapter~4]{ItoMcKean1974}).
\begin{Lemma}\label{lm:closed.diffusion}
	Let Assumption~\ref{ass:b} hold and let $\vpi\in\Ac$ be Borel measurable. There exits a continuous-path Markov process $(S_t^{\vpi})_{t\geq0}$ with infinite lifetime and state space $[0,\pi]$, defined in law for every initial point $s\in[0,\pi]$. Its transition semigroup maps $\Cc([0,\pi])$ into itself, and its $L^2(\widehat m^\vpi ds)$-generator has the	interior differential expression
	\bee
	\Lc^\vpi u=u_{ss}+Q^\vpi(s)u_s=\frac{1}{\widehat m^\vpi(s)}\bigl(\widehat m^\vpi(s)u_s\bigr)_s.
	\eee
	With the domain $\Sc(\vpi)$, its symmetric Dirichlet form on $L^2(\widehat m^\vpi ds)$ is
	\be\label{eq:closed.form}
	\Ec^\vpi(u,v)=\int_0^\pi u_s(s)v_s(s)\widehat m^\vpi(s) ds,\quad \forall u,v\in\Sc(\vpi).
	\ee
Each endpoint is either regular \footnote{Here a regular boundary means to be both accessible from the interior and admits entrance into the interior. This terminology is distinct from the regular solution class defined by the zero-flux condition in \eqref{eq:zero-flux0}.} or entrance,
and
	\be\label{eq:boundary.zero.occupation}
	\E_s^{\vpi}\left[
	\int_0^\infty e^{-t}\mathbf 1_{\{0,\pi\}}(S_t^{\vpi})dt \right]=0,\quad \forall s\in[0,\pi].
	\ee
The above properties uniquely determine the law of $S^\vpi$ for every initial state in $[0,\pi]$.
\end{Lemma}

\begin{proof}
In view of the policy bounds \eqref{eq:Theta.bounds} and \eqref{eq:q.bound}, \eqref{eq:mass-ratio} and \eqref{eq:general-masslog} in Lemma~\ref{lm:ratio} give
	\be\label{eq:closed.mass}
	\int_0^\pi \widehat m^{\vpi}(t) dt<\infty,
	\quad
	\int_0^s \widehat m^{\vpi}(t) dt \leq Cs\widehat m^{\vpi}(s),	\quad \forall 0<s\leq \frac{\pi}{8},
	\ee
	which gives
	\be\label{eq:entrance.test}
	\int_0^{\frac\pi8}\frac{1}{\widehat m^{\vpi}(s)}	\left(\int_0^s \widehat m^{\vpi}(t) dt\right)ds	\leq C\frac{\pi^2}{128}<\infty.
	\ee
	The corresponding integral at the right endpoint $\pi$ is finite by the same argument.
	
	Define the scale function by$\mathfrak s(s)=\int_{\pi/2}^{s}\widehat m^{\vpi}(t)^{-1}dt$, and take the speed measure $\widehat m^{\vpi}(s)ds$ with zero mass at $\{0,\pi\}$. Then the scale--speed boundary criteria in \cite[Section~4.6, pp.~130--131]{ItoMcKean1974} (see also \cite[Section~V.51]{RogersWilliams2000}) together with \eqref{eq:closed.mass}--\eqref{eq:entrance.test} show that each endpoint is either regular or entrance, and exit-only and natural boundaries are excluded. Indeed, the left-endpoint accessibility integral satisfies
	\bee
	\int_0^{\frac{\pi}{8}}\frac{1}{\widehat m^{\vpi}(s)}	\left(\int_s^{\frac{\pi}{8}}\widehat m^{\vpi}(t) dt\right)ds \text{ being finite is equivalent to } \int_0^{\frac{\pi}{8}}\frac{ds}{\widehat m^{\vpi}(s)} \text{ being finite}.
	\eee
	If $\int_0^{\pi/8}\widehat m^{\vpi}(s)^{-1}ds<\infty$, the left endpoint	is regular, and we take the instantaneously reflecting extension. If this scale integral is infinite, the endpoint cannot be reached from the interior but admits a continuous entrance law. The same alternatives hold at $\pi$. These choices, with speed measure $\widehat m^{\vpi}(s)ds$ and no killing, uniquely determine a conservative diffusion on $[0,\pi]$ whose transition semigroup preserves $\Cc([0,\pi])$ (see the generator boundary conditions in \cite[Section~4.7, pp.~135--136]{ItoMcKean1974} and the converse scale--speed construction described in \cite[Section~5.1 and Chapter~5]{ItoMcKean1974}). At a regular endpoint, zero speed-measure mass gives instantaneous reflection and zero occupation time by the occupation-time formula (see \cite[Section~5.4, pp.~174--175]{ItoMcKean1974}). An entrance endpoint is left immediately and cannot be reached from the interior. Consequently, \eqref{eq:boundary.zero.occupation} holds for every initial point.
	
In view that $d\mathfrak s=ds/\widehat m^{\vpi}$ and $du/d\mathfrak s=\widehat m^{\vpi}u_s$, the  symmetric energy on $\Cc^1([0,\pi])$ is
\bee
\int_0^\pi	\frac{du}{d\mathfrak s}\frac{dv}{d\mathfrak s} d\mathfrak s
=\int_0^\pi u_s(s)v_s(s)\widehat m^{\vpi}(s) ds.
\eee
Its closure has domain $\Sc(\vpi)$. The resulting closed form is associated with a symmetric diffusion by \cite[Sections~4.5 and~7.2]{FukushimaOshimaTakeda2011}.
Its scale, speed measure and boundary conditions agree with those of $S^{\vpi}$; hence uniqueness of the scale--speed construction identifies its Dirichlet form with \eqref{eq:closed.form}.
\end{proof}

\begin{Remark}\label{rem:original.diffusion}
For every $s\in[0,\pi]$, the process $P_t^{\vpi}=\sin^2(S_t^{\vpi}/2)$, started from $p=\sin^2(s/2)$, is a weak solution of \eqref{eq:model.X}. Recall the transformation identities in Section~\ref{sec:arcsin.transform}, we have that
\be\label{eq:closed.generator.original}
\begin{aligned}
\Lc^{\vpi}p(s)&=b^{\vpi}(p(s)), \; \text{where}\; b^{\vpi}(p)&=2(1-p)\int_A q^{2,1}(p, a)\vpi_2(p,da)-2p\int_A q^{1,2}(p, a)\vpi_1(p,da),
\end{aligned}
\ee
and $\Lc^{\vpi}(p^2)(s)-2p(s)\Lc^{\vpi}p(s)=2p(s)(1-p(s)).$
By \eqref{eq:main.boundaryrelaxed}--\eqref{eq:assum.remainder}, both generators are bounded. Moreover, \eqref{eq:general-envelope} implies that $\widehat m^{\vpi}p_s$ and $\widehat m^{\vpi}(p^2)_s$ vanish at both endpoints. Then an integration by parts in \eqref{eq:closed.form} implies that $p$ and $p^2$ belong to its $L^2$ generator domain.
	
The generator martingale formula \cite[Chapter~5]{FukushimaOshimaTakeda2011} can be applied outside a set of zero capacity. Since point evaluation is continuous in the energy norm locally in the interior and at endpoints with finite scale distance, interior and regular-boundary points have positive capacity. For an entrance starting point, the same formula follows by applying the Markov property at time $\varepsilon>0$ and letting $\varepsilon\to 0+$, and combining with path continuity and boundedness of the test functions and their generators.
In sum, for every initial point,
\bee
P_t^{\vpi}=P_0^{\vpi}+\int_0^t b^{\vpi}(P_r^{\vpi})dr+M_t,\quad\langle M\rangle_t=\int_0^t2P_r^{\vpi}(1-P_r^{\vpi})dr,
\eee
where $M$ is a continuous martingale. Endpoint values of the generator expressions are immaterial by \eqref{eq:boundary.zero.occupation}. Then the Brownian representation, after enlarging the probability space if necessary, leads to \eqref{eq:model.X}. The law of $P^\vpi$ is uniquely determined by the law of $S^\vpi$ in Lemma~\ref{lm:closed.diffusion} through $P_t^\vpi=\sin^2(S_t^\vpi/2)$.
\end{Remark}

\begin{Remark}\label{rem:measurable.policy.diffusion}
The intrinsic-coordinate transformation and Bessel-type models are motivated by the Kimura diffusion theory developed in \cite{EpsteinMazzeo2010,EpsteinMazzeo2013}. However, the coefficient regularity required in \cite{EpsteinMazzeo2010,EpsteinMazzeo2013} is not ensured by our assumptions: Even when the transition rates are smooth in the population variable, averaging against a Borel feedback policy may produce a discontinuous effective drift $b^{\vpi}$ in \eqref{eq:closed.generator.original}. Therefore, results in \cite{EpsteinMazzeo2010,EpsteinMazzeo2013} cannot be directly applied to every $\vpi\in\Ac$.
	
Lemma \ref{lm:closed.diffusion} verifies the classical scale--speed criteria under the present assumptions and identifies the diffusion associated with the policy-dependent energy space $\Sc(\vpi)$. Together with Remark~\ref{rem:original.diffusion}, it connects the weighted analytic formulation to the original controlled diffusion, including initial points at the boundary, and contributes a justification for Borel feedback policies and its compatibility with the policy-uniform weighted estimates used in the policy iteration analysis.
\end{Remark}

\begin{Proposition}\label{prop:PDE.wellposed}
Fix $k>\max\{2,2\thu\}$ and $0<\alpha<1-\max\{1,2\thu\}/k$. Let $\vpi\in\Ac$ and $F\in L^k((0,\pi),\widehat m^{\vpi}(s) ds)$. Then
\begin{equation}\label{eq:prop.pdegeneral}
\rho u-\frac1{\widehat m^{\vpi}}(\widehat m^{\vpi}u_s)_s=F
\end{equation}
admits a unique weak solution $u^{\vpi,F}\in\Sc(\vpi)$, which is regular in the sense of Definition~\ref{def:regular.linear}. It has a unique representative in $\Cc^{1,\alpha}([0,\pi])$, and
\begin{equation}\label{eq:closed.resolvent}
u^{\vpi,F}(s)=\E_s^{\vpi}
\left[\int_0^\infty e^{-\rho t}F(S_t^{\vpi}) dt\right] \quad \forall s\in[0,\pi],
\end{equation}
where the integral is absolutely integrable for each starting point and is unchanged by modifications of $F$ on $\widehat m^{\vpi}(s)ds$-null sets.
Furthermore,
\begin{equation}\label{eq:closed.resolvent.est}
\|u^{\vpi,F}\|_{L^k(\widehat m^{\vpi}ds)}
\leq\rho^{-1}\|F\|_{L^k(\widehat m^{\vpi}ds)},
\qquad
\|u^{\vpi,F}\|_{\Cc^{1,\alpha}([0,\pi])}
\leq C\|F\|_{L^k(\widehat m^{\vpi}ds)},
\end{equation}
with $C$ depending only on $k,\alpha,\rho,M,\thl,\thu$ but not
on $\vpi$. In particular, for $\vpi\in\Ac_\lambda$ and $F=F^{\vpi}$, we have $W_\lambda^{\vpi}(s)=u^{\vpi,F^{\vpi}}(s)$ for all $s\in [0,\pi]$.


\end{Proposition}

\begin{proof}
{\bf Step 1.} We first show that \eqref{eq:prop.pdegeneral}  admits a unique weak solution in $\Sc(\theta)$ which satisfies the zero-flux condition:
\begin{equation}\label{eq:zero-flux}
 \lim_{s\to0+}\widehat m^{\vpi}(s)u_s^{\vpi,F}(s)=0,\quad \lim_{s\to\pi-}\widehat m^{\vpi}(s)u_s^{\vpi,F}(s)=0.
\end{equation}
Noting $k>2$, $\int_0^\pi\widehat m^{\vpi}ds<\infty$ by Lemma~\ref{lm:ratio}, we have $F\in L^2((0,\pi),\widehat m^{\vpi}ds).$ Consider the bilinear form on $\Sc(\vpi)$
$$
a^{\vpi}(\phi,\psi):=\rho\int_0^\pi\phi\psi\widehat m^{\vpi} ds
+\int_0^\pi\phi_s\psi_s\widehat m^{\vpi} ds.
$$
By Cauchy--Schwarz inequality, we get $|a^\theta(\phi, \psi)|\leq C_\rho \|\phi\|_{\Sc(\theta)}\|\psi\|_{\Sc(\theta)}$, which implies the continuity of $a^\theta$. Noting $a^\theta(\phi,\phi)\geq \min\{1, \rho\} \|\phi\|^2_{\Sc(\theta)}$, we have that $a^\theta$ is coercive. Meanwhile, the mapping
$
\phi \mapsto\int_0^\pi F\phi\widehat m^{\vpi} ds
$
is also continuous on $\Sc(\vpi)$. Then by Lax--Milgram Theorem, there exists a unique $u^{\vpi,F}\in\Sc(\vpi)$ such that
\begin{equation}\label{eq:weak.pde}
\rho\int_0^\pi u^{\vpi,F}\phi\widehat m^{\vpi} ds
+\int_0^\pi u_s^{\vpi,F}\phi_s\widehat m^{\vpi} ds
=\int_0^\pi F\phi\widehat m^{\vpi} ds, \qquad \forall\phi\in\Sc(\vpi).
\end{equation}
Hence, the following holds in the distributional sense:
\begin{equation}\label{eq:mup.derivative}
(\widehat m^{\vpi}u_s^{\vpi,F})_s=\widehat m^{\vpi}(\rho u^{\vpi,F}-F).
\end{equation}
In light of $u^{\vpi,F},F\in L^2(\widehat m^{\vpi}ds)$ and $\int_0^\pi\widehat m^{\vpi}ds<\infty$, it follows from Cauchy--Schwarz inequality again that the right-hand side belongs to $L^1((0,\pi))$. As a result, $(\widehat m^\theta u^{\theta, F}_s)_s \in L^1((0,\pi))$, which further implies that $\widehat m^\theta u^{\theta, F}_s \in W^{1,1}((0,\pi))$. By choosing its absolutely continuous representative, which extends continuously to
$[0,\pi]$, we conclude that
$$
\lim_{s\to 0+}\widehat m^{\vpi}(s)u_s^{\vpi,F}(s) \text{ and } \lim_{s\to \pi-}\widehat m^{\vpi}(s)u_s^{\vpi,F}(s) \text{exist and are finite.}
$$
For an arbitrary $\phi\in  \Cc^1([0,\pi])$, integrating by parts gives
$$
\begin{aligned}
(\widehat m^{\vpi}u_s^{\vpi,F}\phi)(\pi) -  (\widehat m^{\vpi}u_s^{\vpi,F}\phi)(0)
= &  \int_0^\pi u^{\theta, F}_s \phi_s \widehat m^{\theta}ds  +  \int_0^\pi (\widehat m^\theta u_s)_s \phi ds \\
= &\int_0^\pi u^{\theta, F}_s \phi_s \widehat m^{\theta}ds  +  \int_0^\pi \widehat m^{\vpi} \phi(\rho u^{\vpi,F}-F)  ds=0,
\end{aligned}
$$
where the second equality follows from \eqref{eq:mup.derivative}  and the last equality follows from \eqref{eq:weak.pde}.
Then, by the arbitrariness of $\phi\in \Cc^1([0,\pi])$, we conclude \eqref{eq:zero-flux}.

{\bf Step 2.} We now establish the $L^k$ estimate needed to apply Theorem \ref{thm:main.besselest} in Step 3. The process in Lemma \ref{lm:closed.diffusion} has Dirichlet form \eqref{eq:closed.form} on $\Sc(\vpi)$. By the resolvent characterization in \cite[Theorem~1.3.1]{FukushimaOshimaTakeda2011}, its $L^2(\widehat m^{\vpi}ds)$ resolvent $R_\rho^{\vpi}F$ is the unique element in $\Sc(\vpi)$ satisfying the variational identity \eqref{eq:weak.pde}. The uniqueness obtained in Step~1 then results in
\be\label{eq:closed.resolvent.L2}
u^{\vpi,F}=R_\rho^{\vpi}F
\qquad\text{in }L^2(\widehat m^{\vpi}ds).
\ee
In view that \eqref{eq:closed.form} is symmetric and the process has infinite lifetime, its transition semigroup is a contraction on $L^k(\widehat m^{\vpi}ds)$ for $k\geq1$. Applying Minkowski's inequality to the discounted semigroup then yields
\be\label{eq:closed.Lk.contraction}
\|u^{\vpi,F}\|_{L^k(\widehat m^{\vpi}ds)}
\leq\int_0^\infty e^{-\rho t} 
\|P_t^{\vpi}F\|_{L^k(\widehat m^{\vpi}ds)} dt
\leq\frac1\rho\|F\|_{L^k(\widehat m^{\vpi}ds)}.
\ee
Here $P_t^{\vpi}$ denotes the transition semigroup of $S^{\vpi}$. The identification of the pointwise stochastic resolvent, including the endpoints, is proved in Step 4 below.

{\bf Step 3.} To goal is to establish the regularity of $u^{\vpi,F}$. The $L^k$ estimate from Step 2 gives $u^{\vpi,F}\in L^k(\widehat m^{\vpi}ds)$. Fix
 $0<\delta<\pi/2$ and set $I_\delta=[\delta,\pi-\delta]$. On $I_\delta$, $Q^{\vpi}$ is bounded, then
\be\label{eq:prop.step3.0}
 0<c_\delta\leq \widehat m^{\vpi}(s)\leq C_\delta<\infty,
 \quad \text{and }\; \widehat m^{\vpi}, (\widehat m^{\vpi})^{-1} \in W^{1,\infty}(I_\delta).
\ee
Also, \eqref{eq:mup.derivative}, $\|u^{\vpi,F}\|_{L^k((0,\pi),\widehat m^{\vpi}ds)}\leq\frac1\rho \|F\|_{L^k((0,\pi),\widehat m^{\vpi}ds)}<\infty$, and the equivalence of $\widehat m^{\vpi}$ with the Lebesgue measure together imply
 $
(\widehat m^{\vpi}u_s^{\vpi,F})_s\in L^k(I_\delta).
 $
 Moreover, it follows from $u^{\vpi,F}\in\Sc(\vpi)$ that $\widehat m^{\vpi}u_s^{\vpi,F}\in L^1(I_\delta)$. It then holds that $\widehat m^{\vpi}u_s^{\vpi,F}\in W^{1,k}(I_\delta).$ Thanks to \eqref{eq:prop.step3.0}, we have 
 $$
 u_s^{\vpi,F}
 = (\widehat m^{\vpi})^{-1} \big(\widehat m^{\vpi}u_s^{\vpi,F}\big) \in W^{1,k}(I_\delta),
 $$
which tells $u^{\vpi,F}\in W^{2,k}(I_\delta).$ As $\delta>0$ is arbitrary, it holds that
 $
 u^{\vpi,F}\in W^{2,k}_{\rm loc}((0,\pi)).
 $
 For $k>1$, the classic one-dimensional Sobolev embedding further gives $u^{\vpi,F}\in \Cc^1((0,\pi))$. Together with the zero-flux condition \eqref{eq:zero-flux}, this shows that $u^{\vpi,F}$ is a regular solution in the sense of Definition~\ref{def:regular.linear}. 
The zero-flux condition from Step 1 and the regularity just proved allow us to apply Theorem \ref{thm:main.besselest} Part (ii) with the uniform equivalence of $\widehat m^{\vpi}ds$ and
$m^{\vpi}ds$. Then it follows from Theorem \ref{thm:main.besselest} Part (i) that
\begin{equation}\label{eq:prop.regular.closed.estimate}
\|u^{\vpi,F}\|_{\Cc^{1,\alpha}([0,\pi])}
\leq C\|F\|_{L^k(\widehat m^{\vpi}ds)}.
\end{equation}

{\bf Step 4.} We now prove \eqref{eq:closed.resolvent} for all starting points, without assuming that $F$ is continuous. First suppose $F\in\Cc([0,\pi])$. As the transition semigroup in Lemma \ref{lm:closed.diffusion} preserves $\Cc([0,\pi])$, the stochastic resolvent
\bee
\mathcal R_\rho^{\vpi}F(s):=\E_s^{\vpi}\left[\int_0^\infty e^{-\rho t}F(S_t^{\vpi}) dt\right],\quad\forall  s\in[0,\pi],
\eee
is continuous on $[0,\pi]$. By \eqref{eq:closed.resolvent.L2} and the $\Cc^{1,\alpha}([0,\pi])$ regularity established in Step~3, $\mathcal R_\rho^{\vpi}F$ and $u^{\vpi,F}$ are continuous functions that agree $\widehat m^{\vpi}ds$-almost everywhere. Since this measure charges every nonempty open subinterval of $(0,\pi)$, they agree on $[0,\pi]$. Hence \eqref{eq:prop.regular.closed.estimate} gives
\be\label{eq:closed.kernel.bound}
\sup_{s\in[0,\pi]}|\mathcal R_\rho^{\vpi}h(s)|
\leq C\|h\|_{L^k(\widehat m^{\vpi}ds)},\quad \forall h\in\Cc([0,\pi]).
\ee
For a fixed $s\in[0,\pi]$, by defining the finite positive measure $\mathcal K_s$ as 
\bee
\mathcal K_s(E):=\E_s^{\vpi}\left[\int_0^\infty e^{-\rho t}
\mathbf 1_{E}(S_t^{\vpi}) dt\right]\quad \text{for all Borel } E\subset[0,\pi],
\eee
the map $h\mapsto\mathcal R_\rho^{\vpi}h(s)$ is integration against $\mathcal K_s$. The estimate \eqref{eq:closed.kernel.bound} extends this positive linear functional continuously to $L^k(\widehat m^{\vpi}ds)$. By the $L^k$ duality, there exists a density $g_s\in L^{k'}(\widehat m^{\vpi}ds)$ with $k'=k/(k-1)$ such that
$
\mathcal K_s(E)=\int_E g_s(t)\widehat m^{\vpi}(t) dt$ and $\|g_s\|_{L^{k'}(\widehat m^{\vpi}ds)}\leq C$,
which implies $\mathcal K_s(\{0,\pi\})=0$ and
\be\label{eq:closed.integrability}
\E_s^{\vpi}\left[\int_0^\infty e^{-\rho t}|F(S_t^{\vpi})| dt\right]
\leq C\|F\|_{L^k(\widehat m^{\vpi}ds)}\quad \text{for every Borel $F\in L^k(\widehat m^{\vpi}ds)$}.
\ee
Therefore, the expectation is finite for every $s\in[0,\pi]$ and depends only on the weighted-$L^k$ equivalence class of $F$.

Finally, choose $F_j\in\Cc([0,\pi])$ converging to $F$ in $L^k(\widehat m^{\vpi}ds)$. By \eqref{eq:prop.regular.closed.estimate}, the associated PDE solutions converge in $\Cc([0,\pi])$ to $u^{\vpi,F}$. Meanwhile, \eqref{eq:closed.integrability} gives
\bee
\sup_{s\in[0,\pi]}\left|\mathcal R_\rho^{\vpi}F_j(s)-\E_s^{\vpi}\left[\int_0^\infty e^{-\rho t}F(S_t^{\vpi}) dt\right]\right|
\leq C\|F_j-F\|_{L^k(\widehat m^{\vpi}ds)}\to 0.
\eee
Then taking limit for each continuous $F_j$ proves \eqref{eq:closed.resolvent} for all $s\in[0,\pi]$. The estimates in \eqref{eq:closed.resolvent.est} follow from \eqref{eq:closed.Lk.contraction} and \eqref{eq:prop.regular.closed.estimate}. By Remark~\ref{rem:original.diffusion}, $P_t^{\vpi}=p(S_t^{\vpi})$ is the original-coordinate diffusion satisfying \eqref{eq:model.X}. For $F=F^{\vpi}$ in \eqref{eq:Fpi}, the change of variables $p=\sin^2(s/2)$ identifies the stochastic integral in \eqref{eq:closed.resolvent} with the value defined in \eqref{eq:value.regularized.x}. Hence $u^{\vpi,F^{\vpi}}=W_\lambda^{\vpi}$ on $[0,\pi]$.

\end{proof}

With the help of Proposition \ref{prop:PDE.wellposed}, we are ready to show the well-posedness of the PIA.
\begin{Corollary}\label{cor:PIA.welldefined}
The PIA under both coorfinate formulations are well-defined. For every $n\geq0$, $\vpi^n\in\Ac_\lambda$, $u^n=W_\lambda^{\vpi^n}$ and $v^n=V_\lambda^{\vpi^n}$ satisfy \eqref{eq:PI.coordinates}, the value identity holds also at the endpoints. Moreover,
$$
F^{\vpi^n}\in L^{k}((0,\pi),\widehat m^{\theta^n} ds), \quad \forall 1\leq k<\infty,
$$
and $u^n$ is the unique weak solution in $\Sc(\vpi^n)$ of the $n$th policy-evaluation equation:
\bee
  -\rho u+u_{ss}+Q^{\theta^{n}} u_s+ F^{\vpi^n}=0,
\eee
which is regular in the sense of Definition~\ref{def:regular.linear}. In addition, for any $ n\geq0$ and  $0<\alpha<1$, we have $u^n\in  \Cc^{1,\alpha}([0,\pi])$ and $v^n\in\Cc_K^{1,\alpha}([0,1])$.

\end{Corollary}

\begin{proof}
 By Lemma \ref{lm:closed.diffusion} and Remark~\ref{rem:original.diffusion}, once a Borel policy $\vpi^n$ is constructed, $S^{\vpi^n}$ and the corresponding Wright--Fisher diffusion $P^{\vpi^n}=p(S^{\vpi^n})$ exists for every initial state, with uniquely specified laws under the prescribed boundary conditions. The former satisfies \eqref{eq:boundary.zero.occupation}, and the latter spends zero Lebesgue time at the endpoints $\{0,1\}$.

Now we argue by induction. Fix ${k_0}>\max\{2,2\thu\}$ sufficiently large such that $0<\alpha<1-\max\{1,2\thu\}/{k_0}$. Invoking $\vpi_1^0\equiv\vpi_2^0\equiv\nu$ gives $\Hc(\vpi_i^0)\equiv \ln(\operatorname{Leb}(A))$ for $i=1,2$. In addition, \eqref{eq:general-masslog} in Lemma \ref{lm:ratio} yields $ \int_0^\pi \widehat m^{\theta^0}(s) ds<\infty.$ Then by the boundedness of $\widetilde r_1,\widetilde r_2$,
$$
F^{\vpi^0}\in L^{k}((0,\pi),\widehat m^{\theta^0} ds)\quad \forall 1\leq {k}<\infty.
$$
In particular, $ F^{\vpi^0}\in L^{{k_0}}((0,\pi),\widehat m^{\theta^0} ds)$. Thanks to Proposition \ref{prop:PDE.wellposed}, there exists a unique weak energy solution $u^0\in\Sc(\vpi^0)$ to the PDE
$$
-\rho u+u_{ss}+Q^{\vpi^0} u_s+F^{\vpi^0}=0,
$$
which is regular in the sense of Definition~\ref{def:regular.linear}. Moreover, Proposition~\ref{prop:PDE.wellposed} yields $u^0\in\Cc^{1,\alpha}([0,\pi])$, and $u^0(s)=W^{\vpi^0}_\lambda(s)$ for every $s\in[0,\pi]$, with an absolutely integrable stochastic payoff even for boundary initial conditions. 

Now suppose, for some $n\geq0$, that $\theta^n ,u^n$ are well-defined and $u^n\in \Cc^{1,\alpha}([0,\pi]).$ Then $\|u_s^n\|_{L^\infty([0,\pi])}<\infty,$
which gives that $\vpi_i^{n+1}(s,a)=\Gamma_{\lambda,i}(s,u_s^n(s),a)$ is well-defined for $i=1,2$. Lemma \ref{lm:entropy.est} then gives
\be\label{eq:cor.est1}
|\Hc(\vpi_i^{n+1}(s))|
\leq  C_\lambda\left[ 1+\ln\left(1+\|u_s^n\|_{L^\infty((0,\pi))}\right)+|\ln(\sin s)| \right] \quad \forall s\in (0,\pi).
\ee
Invoking \eqref{eq:general-masslog} in Lemma \ref{lm:ratio}, we get that, for any $1\leq {k}<\infty$,
$$
\int_0^\pi \left(1+|\ln(\sin s)|\right)^{k}\widehat{m}^{\theta^{n+1}} (s) ds<\infty.
$$
This together with \eqref{eq:cor.est1}, the weights $0\leq p(s),1-p(s)\leq1$, and Assumption \ref{ass:action} (iv) gives, for every $1\leq {k}<\infty$ that
$
F^{\vpi^{n+1}} \in L^{k}((0,\pi),\widehat{m}^{\theta^{n+1}}  ds).
$
In particular, $F^{\vpi^{n+1}} \in L^{{k_0}} ((0,\pi),\widehat{m}^{\theta^{n+1}} ds)$ and $\theta^{n+1}\in \Ac_\lambda$. Then Proposition \ref{prop:PDE.wellposed} yields the unique weak energy solution $u^{n+1}\in \Sc(\theta^{n+1})$ to the PDE:
$$
-\rho u+u_{ss}+Q^{\vpi^{n+1}} u_s+F^{\vpi^{n+1}}=0,
$$
which is regular in the sense of Definition \ref{def:regular.linear}, and $u^{n+1}(s)=W^{\theta^{n+1}}_\lambda(s)$ for any $s\in[0,\pi]$. Proposition \ref{prop:PDE.wellposed} further gives
$u^{n+1}\in \Cc^{1,\alpha}([0,\pi])$, which closes the induction. As a result, the assertions under the original coordinates follow from \eqref{eq:PI.coordinates}, the stochastic representation in Proposition~\ref{prop:PDE.wellposed} as well as \eqref{eq:CK.original}.
\end{proof}

\subsection{Policy improvement}\label{subsec:improvement}

To prove the {\it policy improvement} property of the PIA, we first upgrade the zero-flux condition in Definition \ref{def:regular.linear} to the following.

 \begin{Lemma}\label{lm:endpoint.crossflux}
Let $P$ have the form \eqref{eq:Pdef} with $\Lambda_1,\Lambda_2$ satisfying \eqref{eq:theta.bounds}, and suppose $\|q\|_{L^\infty((0,\pi))}\leq M$. Let $u$ be a regular solution to \eqref{eq:pde.bessel} in the sense of Definition~\ref{def:regular.linear}. 
	If
	\begin{equation}\label{eq:log-source.bound}
		|\rho u+f|\leq C_0\bigl(1+|\ln(\sin s)|\bigr)\quad \text{for a.e. } s\in(0,\pi),
	\end{equation}
	with $C_0<\infty$, then there exists a constant $C_1<\infty$ depending only on $C_0$ such that:
	\begin{equation}\label{eq:hard.regular.endpoint}
		|u_s(s)|\leq C_1 s(1+|\ln s|),\;\; 0<s\leq \frac{\pi}{8};
		\quad
		|u_s(s)|\leq C_1(\pi-s)(1+|\ln(\pi-s)|), \;\; \frac{7\pi}{8}\leq s<\pi.
	\end{equation}
	Moreover, if $\widetilde m(s)=\exp\left(\int_{\pi/2}^{s}(\widetilde P(t)+\widetilde q(t)) dt\right)$ is any full Bessel weight whose coefficient $\widetilde P$ has the form
	\eqref{eq:Pdef} satisfying \eqref{eq:theta.bounds} and whose remainder
	$\widetilde q$ is bounded, then
	\bee
		\lim_{s\to 0+}\widetilde m(s)u_s(s) = \lim_{s\to \pi-}\widetilde m(s)u_s(s)=0.
	\eee
\end{Lemma}

\begin{proof}
Notice that \eqref{eq:pde.bessel} and \eqref{eq:P.transform}  give
$
\frac1{m^q(s)}(m^q(s) u_s)_s = \rho u +f.
$
Since, $u$ satisfies the zero-flux condition at $0$, integrating from 0 to $s$ gives
$
u_s(s)=\frac1{m^q(s)}\int_0^s m^q(t)[ \rho u +f]dt.
$
Then by the first line of \eqref{eq:lm.bessel.ratio} in Lemma \ref{lm:ratio} and \eqref{eq:log-source.bound}, for $0<s\leq \pi/8$, we have
$$
\begin{aligned}
|u_s(s)|\leq & C C_0\int_0^s\left(\frac ts\right)^{2\thl-1}(1+|\ln t|) dt  = CC_0 s\int_0^1 r^{2\thl -1} (1+|\ln(sr)|)dr\\
\leq &CC_0 s\int_0^1 r^{2\thl -1} (1+|\ln s|+|\ln r|)dr 
 \leq  CC_0 \left[s(1+|\ln s|)\int_0^1 r^{2\thl -1} dr + s\int_0^1 r^{2\thl -1} |\ln r| dr\right]\\
 \leq & C_1 s(1+|\ln s|),
\end{aligned}
$$
where the first equality follows from the transformation $t=sr$ and the last inequality follows from $\thl>0$.
The right-endpoint estimate follows identically from the zero-flux condition at $\pi$ and the second line of \eqref{eq:lm.bessel.ratio}. Hence \eqref{eq:hard.regular.endpoint} holds.
	
Now let $\widetilde m$ be any full Bessel weight satisfying the conditions. By \eqref{eq:general-envelope} in Lemma \ref{lm:ratio}, we have
$\widetilde m(s) \leq Cs^{2\thl-1}$ near $0$. Therefore, by \eqref{eq:hard.regular.endpoint} and $\thl>0$,
$
|\widetilde m(s)u_s(s)|\leq C s^{2\thl}(1+|\ln s|)\to0
$ as $s\to 0+$.
The proof of the other endpoint $\pi$ follows the same argument.
\end{proof}

We refine in the next result the zero-flux property of $u^n$,  as stated in Proposition \ref{prop:PDE.wellposed}.

\begin{Lemma}
 \label{lm:individual.endpoint}
 For every fixed $n\geq0$, there exists $C_n<\infty$ such that
 \be\label{eq:ind.leftright}
 \begin{cases}
  |u_s^n(s)|\leq C_ns(1+|\ln s|),&0<s\leq \frac{\pi}{8},\\
  |u_s^n(s)|\leq C_n(\pi-s)(1+|\ln(\pi-s)|), &\frac{7\pi}{8}\leq s<\pi.
 \end{cases}
 \ee
 In particular,$ u_s^n(0)=u_s^n(\pi)=0.$ Moreover, for every feedback policy pair $\widetilde{\vpi}$,
 $$
 \lim_{s\to 0+} \widehat m^{\widetilde{\vpi}}(s)u_s^n(s) = \lim_{s\to \pi-} \widehat m^{\widetilde{\vpi}}(s)u_s^n(s)=0.
 $$
\end{Lemma}

\begin{proof}
 From the $n$th policy-evaluation equation, we have
 $$
 \frac1{\widehat m^{\theta^n}} (\widehat m^{\theta^n} u_s^n)_s = \rho u^n-F^{\vpi^n}.
 $$
 For fixed $n\geq1$, Corollary~\ref{cor:PIA.welldefined} gives that $u^n\in\Cc^{1,\alpha}([0,\pi])$ and satisfies the zero-flux condition. Moreover,  in view of $\vpi_i^n(s,a)=\Gamma_{\lambda,i}(s,u_s^{n-1}(s),a)$, Lemma~\ref{lm:entropy.est}, and the boundedness of $\widetilde r_1,\widetilde r_2$, we deduce that
\be\label{eq:lm.Fnbound}
 |\rho u^n(s)-F^{\vpi^n}(s)| \leq C_n\bigl(1+|\ln(\sin s)|\bigr).
\ee
For $n=0$, the same bound holds because $u^0$ is bounded and the uniform initial policy has bounded running reward $F^{\vpi^0}$.
Then applying Lemma~\ref{lm:endpoint.crossflux} with $m^q=\widehat m^{\theta^n}$ yields that \eqref{eq:ind.leftright} and for any feedback policy $\theta$ 
$$
 \lim_{s\to 0+} \widehat m^{\widetilde{\vpi}}(s)u_s^n(s)=\lim_{s\to \pi-} \widehat m^{\widetilde{\vpi}}(s)u_s^n(s)=0.
 $$
Thanks to
 $u^n\in\Cc^{1,\alpha}([0,\pi])$, \eqref{eq:ind.leftright} also gives
 $
 u_s^n(0)=u_s^n(\pi)=0.
 $
\end{proof}

\begin{Proposition}[Policy improvement]\label{prop:policy.improvement}
 For every $n\geq 0$, it holds that
\bee
u^{n+1}(s)\geq u^n(s),\quad \forall s\in[0,\pi],\qquad v^{n+1}(p)\geq v^n(p),\quad \forall p\in[0,1].
\eee
\end{Proposition}

\begin{proof}
 By the Gibbs variational property and the definition
 $\vpi_i^{n+1}=\Gamma_{\lambda,i}(s,u_s^n,\cdot)$, $i=1,2$, we have
\bee
Q^{\vpi^{n+1}}u_s^n+F^{\vpi^{n+1}}\geq Q^{\vpi^n}u_s^n+F^{\vpi^n}.
\eee
This together with
 $
 -\rho u^n+u_{ss}^n+Q^{\vpi^n}u_s^n+F^{\vpi^n}=0
 $
gives
 \begin{equation}\label{eq:sub.improved}
  -\rho u^n+u_{ss}^n+Q^{\vpi^{n+1}}u_s^n
  +F^{\vpi^{n+1}}\geq 0.
 \end{equation}
 On the other hand, $u^{n+1}$ satisfies
 \begin{equation}\label{eq:sol.improved}
  -\rho u^{n+1}+u_{ss}^{n+1}+Q^{\vpi^{n+1}}u_s^{n+1}
  +F^{\vpi^{n+1}}=0.
 \end{equation}
 Set $h=u^n-u^{n+1}$. Subtracting \eqref{eq:sol.improved} from \eqref{eq:sub.improved} gives
 \begin{equation}\label{eq:hsub}
  -\rho h+\frac1{\widehat{m}^{\theta^{n+1}} }(\widehat{m}^{\theta^{n+1}} h_s)_s\geq 0.
 \end{equation}
  By Proposition~\ref{prop:PDE.wellposed}, $u^{n+1}$ has zero flux for $\widehat{m}^{\theta^{n+1}} $. By Lemma~\ref{lm:individual.endpoint} with ${\widetilde{\vpi}}=\vpi^{n+1}$, $u^n$ also has zero flux for $\widehat{m}^{\theta^{n+1}} $. Hence $h$ has zero flux for $\widehat{m}^{\theta^{n+1}} $. That is,
 \begin{equation}\label{eq:h.zero.flux}
\widehat{m}^{\theta^{n+1}} (s)h_s(s)\to 0, \text{ as }  s\to0+ \text{ and } s\to\pi-.
 \end{equation}
 Since $h$ belongs to the corresponding weighted Sobolev space, so does $h^+$, and the Sobolev chain rule gives
 $ (h^+)_s=\mathbf 1_{\{h>0\}}h_s$  a.e. on $(0,\pi)$.
 In particular,
 $$
 h h^+=(h^+)^2, \quad  h_s(h^+)_s=|(h^+)_s|^2  \quad\text{a.e. on $(0,\pi)$.}
 $$
 Take an arbitrary $\varepsilon\in(0,\pi/2)$.  By multiplying \eqref{eq:hsub} by $h^+\widehat{m}^{\theta^{n+1}} $ and integrating over $(\varepsilon,\pi-\varepsilon)$, we have that
 $$
 -\rho\int_\varepsilon^{\pi-\varepsilon} (h^+)^2(s)\widehat{m}^{\theta^{n+1}} (s)ds+\int_\varepsilon^{\pi-\varepsilon}(\widehat{m}^{\theta^{n+1}} h_s)_s(s) h^+(s)ds \geq0.
 $$
Then integrating by parts yields
 $$
 \begin{aligned}
 \int_\varepsilon^{\pi-\varepsilon}(\widehat{m}^{\theta^{n+1}}  h_s)_s(s)  h^+(s)  ds
 &=\big[\widehat{m}^{\theta^{n+1}} (s) h_s(s) h^+(s)\big]_{\varepsilon}^{\pi-\varepsilon}
-\int_\varepsilon^{\pi-\varepsilon}\widehat{m}^{\theta^{n+1}} (s)  h_s(s) (h^+)_s(s) ds\\
 &=\big[\widehat{m}^{\theta^{n+1}} (s)  h_s(s)  h^+(s) \big]_{\varepsilon}^{\pi-\varepsilon}
 -\int_\varepsilon^{\pi-\varepsilon}|(h^+)_s|^2(s) \widehat{m}^{\theta^{n+1}} (s) ds.
 \end{aligned}
 $$
 Therefore,
 \begin{equation}
 \label{eq:positive.part.eps}
 \begin{aligned}
 0\leq&-\rho\int_\varepsilon^{\pi-\varepsilon}(h^+)^2 (s) \widehat{m}^{\theta^{n+1}} (s) ds
 -\int_\varepsilon^{\pi-\varepsilon}|(h^+)_s(s) |^2\widehat{m}^{\theta^{n+1}} (s)ds+\big[\widehat{m}^{\theta^{n+1}} (s)  h_s(s)  h^+(s) \big]_{\varepsilon}^{\pi-\varepsilon}.
 \end{aligned}
 \end{equation}
 Since $u^n,u^{n+1}\in  \Cc^{1,\alpha}([0,\pi])$, $h^+$ is bounded, then \eqref{eq:h.zero.flux} implies
 $
 \big[\widehat{m}^{\theta^{n+1}}  h_s h^+\big]_{\varepsilon}^{\pi-\varepsilon}\to 0\;\text{as }\varepsilon\to 0+.
 $
Hence, by letting $\varepsilon\to 0+$ in \eqref{eq:positive.part.eps}, we obtain
 $$
 -\rho\int_0^\pi (h^+)^2(s)\widehat{m}^{\theta^{n+1}}  (s) ds - \int_0^\pi |(h^+)_s(s)|^2\widehat{m}^{\theta^{n+1}} (s)ds \geq0.
 $$
 Both terms on the left-hand side are non-positive, and thus vanish. As a result, $h^+=0$, which, together with the continuity of $h^+$, implies that $u^n\leq u^{n+1}$ on $[0,\pi]$. Moreover, the inequality for $v^n$ follows from \eqref{eq:PI.coordinates}.
\end{proof}

\begin{Corollary}\label{cor:limit}
 There exists $K<\infty$ such that
 \be\label{eq:un.bound}
 \sup_{n\geq0}\norm{u^n}_{L^\infty([0,\pi])}\leq K.
 \ee
 Hence,
 $u^\infty(s):=\lim_{n\to\infty}u^n(s)$ is well defined and finite for every $s\in[0,\pi]$. Then \eqref{eq:mass.uniform} gives
 \be\label{eq:un.lpbound}
 \norm{u^n}_{L^k(m^{\theta^n} ds)} \leq K \left( \int_0^\pi m^{\theta^n}(s)ds \right)^{1/k} \leq C K,\; \text{uniformly in $n$.}
 \ee
\end{Corollary}

\begin{proof}
 Notice that
 $$
 u^n(s)\leq \frac{1}{\rho} \left( \sup_{i,p, a}|r(i,p, a)|+\lambda |\ln(\Leb(A))| \right) \text{ on $(0,\pi)$},\quad \forall n\geq 1.
 $$
Noting that $u^0$ is bigger or equal to the negative of the right-hand side above, we conclude \eqref{eq:un.bound} with $K=(\|r\|_\infty+\lambda|\ln\Leb(A)|)/\rho$ by the policy improvement and \eqref{eq:PI.coordinates}. 
\end{proof}

\section{Convergence of PIA to Optimal Policy}\label{sec:PIA.convergence}

Now we are ready to provide the main compactness result, which leads to the convergence of PIA.

\begin{Theorem}[Uniform weighted H\"older bound for PIA]\label{thm:compact.un}
Suppose Assumptions~\ref{ass:b} and \ref{ass:action} hold.  Then there exists $A^*<\infty$, depending only on $u^0, \alpha, \lambda$ and the fixed model constants, but not on $n$, such that
\bee
\sup_{n\geq0}\norm{u^n}_{\Cc^{1,\alpha}([0,\pi])} \leq A^*.
\eee
Invoking the inverse arcsine transformation, we get an equivalent uniform bound of the corresponding weighted H\"older norms of the original value sequence:
\bee
\sup_{n\geq0} \norm{v^n}_{C_{K}^{1,\alpha}([0,1])} \leq C A^*.
\eee
Moreover, the constant $A^*$ can be chosen uniformly over $\lambda\in(0,1]$.
\end{Theorem}

\begin{proof}
Substituting $y=u_s^{n-1}(s)$ into Lemma~\ref{lm:entropy.est} and using $0\leq p(s),1-p(s)\leq1$, we get
$$
|\lambda p(s)\Hc(\vpi_1^n(s))+\lambda(1-p(s))\Hc(\vpi_2^n(s))|\leq C\left[\lambda+\lambda|\ln\lambda|+1+\ln(1+|u_s^{n-1}(s)|)+|\ln(\sin s)|\right].
$$
Let $p$ be as in Theorem \ref{thm:main.besselest}.  Combining the above inequality with \eqref{eq:general-masslog} in Lemma \ref{lm:ratio}, we obtain the existence of constant $C_\lambda>0$, independent of $n$, such that
\be\label{eq:entropy}
\norm{\lambda p\Hc(\vpi_1^n)+\lambda(1-p)\Hc(\vpi_2^n)}_{L^p(m^{\theta^n} ds)}
\leq C\lambda \left[  1+ |\ln\lambda|+1+ \ln\left(1+\norm{u^{n-1}}_{\Cc^{1,\alpha}([0,\pi])}\right)\right].
\ee
Then applying Theorem \ref{thm:main.besselest} Part (i) and Theorem \ref{thm:main.besselest} Part (ii) to \eqref{eq:un.PDE}, together with \eqref{eq:un.lpbound}, \eqref{eq:entropy} and $\|r\|_\infty<\infty$, we get that
\begin{align}
 \|u^n\|_{\Cc^{1,\alpha}([0,\pi])}
 &\leq C_E\|u^n\|_{W^{2,p}(m^{\theta^n} ds)}\notag\\
 &\leq C_EC_S\bigg(\|u^n\|_{L^p(m^{\theta^n} ds)}+\norm{p\rtild_1^{\vpi_1^n}+(1-p)\rtild_2^{\vpi_2^n}}_{L^p(m^{\theta^n} ds)}\notag\\
 &\qquad \qquad \quad +\norm{\lambda p\Hc(\vpi_1^n)+\lambda(1-p)\Hc(\vpi_2^n)}_{L^p(m^{\theta^n} ds)}\bigg)\notag\\
 &\leq C_0+\lambda C_1\left[1 +|\ln \lambda|+\ln\left(1+ \|u^{n-1}\|_{\Cc^{1,\alpha}([0,\pi])} \right) \right],
 \label{eq:recursion}
\end{align}
where $C_0,C_1$ are independent of $n$ and $\lambda$. For fixed $\lambda>0$, logarithmic growth allows us to choose $A^*\geq\|u^0\|_{\Cc^{1,\alpha}([0,\pi])}$ such that
$
C_0+\lambda C_1[1+|\ln\lambda|+\ln(1+A^*)]\leq A^*.
$
Then an induction in \eqref{eq:recursion} gives $\sup_{n\geq0}\|u^n\|_{\Cc^{1,\alpha}([0,\pi])}\leq A^*$. The desired bound on $v^n$ follows from \eqref{eq:CK.original} and \eqref{eq:PI.coordinates}. Moreover, in view that $\lim_{\lambda\to 0+}\lambda|\ln\lambda|=0$, the constant bound $A^*$ can be chosen independently of $\lambda\in(0,1]$.
\end{proof}

Define the two soft Hamiltonians
$\mathfrak H_{\lambda,i}(s,y):=\sup_{\theta\in\Pc(A)}\left\{\int_A[Q_i(s,a)y+\rtild_i(s,a)]\theta(a)da+\lambda\Hc(\theta)\right\}$ for $i=1,2$,
and the weighted Hamiltonian
$$
\mathfrak H_\lambda(s,y):=p(s)\mathfrak H_{\lambda,1}(s,y)+(1-p(s))\mathfrak H_{\lambda,2}(s,y).
$$
The Gibbs variational formula gives
$$
\mathfrak H_{\lambda,i}(s,y)=\lambda\ln\int_A\exp\left(\frac{Q_i(s,a)y+\rtild_i(s,a)}{\lambda}\right)da,
$$
and its unique maximizer is $\Gamma_{\lambda,i}(s,y,\cdot)$.

\begin{Theorem}[Convergence of PIA and optimality]\label{thm:convergence}
 Let the assumptions in Theorem \ref{thm:compact.un} hold. Recall $u^\infty$ defined in Corollary \ref{cor:limit}. Then,
$u^\infty$ solves the entropy-regularized HJB equation
  \begin{equation}\label{eq:HJB.ulim}
   -\rho u^\infty+u_{ss}^\infty+\mathfrak H_\lambda(s,u_s^\infty)=0,\quad \text{a.e. on $(0,\pi)$.}
  \end{equation}
The Gibbs policy pair $\vpi_i^\infty(s,a):=\Gamma_{\lambda,i}(s,u_s^\infty(s),a)$, $i=1,2$, is optimal, and
  \begin{equation}\label{eq:optimal.value}
   u^\infty=W^{\vpi^\infty}_\lambda=\sup_{\vpi\in \Ac_\lambda}W^\vpi_\lambda.
  \end{equation}
In the original coordinate, the same feedback policy attains the optimal MFC value:
\begin{equation}\label{eq:PI.original.optimal}
v^\infty(p)=u^\infty(\Phi(p))=V_\lambda^{\vpi^\infty}(p)=v^{*,\lambda}(p)\quad \forall p\in[0,1].
\end{equation}
Moreover, $v^\infty$ solves \eqref{eq:HJB.x} a.e. on $(0,1)$, and for every $0<\beta<\alpha$,
\begin{equation}\label{eq:PI.both.convergence}
\lim_{n\to\infty}\|v^n-v^{*,\lambda}\|_{\Cc_K^{1,\beta}([0,1])}
=\lim_{n\to\infty}\|u^n-u^\infty\|_{\Cc^{1,\beta}([0,\pi])}=0.
\end{equation}
\end{Theorem}

\begin{proof}
By the compact embedding $ \Cc^{1,\alpha}([0,\pi])\Subset  \Cc^{1,\beta}([0,\pi])$, every subsequence has a further subsequence converging in $ \Cc^{1,\beta}$. Any such limit must equal the pointwise limit $u^\infty$. Therefore $\|u^n-u^\infty\|_{\Cc^{1,\beta}([0,\pi])}\to 0$ for any $\beta\in (0,\alpha)$. In particular,
 \begin{equation}\label{eq:derivative.converge}
  \|u_s^n-u_s^\infty\|_\infty\to0.
 \end{equation}

 By virtue of $\vpi_i^{n+1}=\Gamma_{\lambda,i}(s,u_s^n,\cdot)$, $i=1,2$, the Gibbs identity gives
$
 \mathfrak H_\lambda(s,u_s^n)=Q^{\vpi^{n+1}}u_s^n+F^{\vpi^{n+1}}.
$
 Then the PDE \eqref{eq:PDEpi} for $u^{n+1}$ becomes
 \begin{equation}\label{eq:PDE.unrewrite}
  -\rho u^{n+1}+u_{ss}^{n+1} +\mathfrak H_\lambda(s,u_s^n)+Q^{\vpi^{n+1}}(u_s^{n+1}-u_s^n)=0.
 \end{equation}

 Fix $\delta\in(0,\pi/2)$ and set $I_\delta=[\delta,\pi-\delta]$. On $I_\delta$, $Q_i(s,a)$ is uniformly bounded over $a$ for $i=1,2$, and hence $Q^{\vpi^{n+1}}$ is uniformly bounded. Therefore, \eqref{eq:derivative.converge} gives
 $
 \|Q^{\vpi^{n+1}} (u_s^{n+1}-u_s^n)\|_{L^\infty(I_\delta)}\to  0.
 $
 Moreover, the Hamiltonian is uniformly Lipschitz in $y$ on $I_\delta$, and hence
 $
 \mathfrak H_\lambda(s,u_s^n)\to   \mathfrak H_\lambda(s,u_s^\infty)
 $
 uniformly on $I_\delta$. As $u^{n+1}\to u^\infty$ uniformly, \eqref{eq:PDE.unrewrite} gives
 $$
 u_{ss}^{n+1}\to  \rho u^\infty-\mathfrak H_\lambda(s,u_s^\infty) \quad\text{in }L^\infty(I_\delta).
 $$
 It remains to identify $\rho u^\infty-\mathfrak H_\lambda(s,u_s^\infty)$ as the second weak derivative
 of $u^\infty$. For every $\varphi\in C_c^\infty((\delta, \pi-\delta))$, $\int_{I_\delta}u_{ss}^{n+1}\varphi ds=-\int_{I_\delta}u_s^{n+1}\varphi_s ds.$ Passing to the limit, using the preceding $L^\infty$ convergence and \eqref{eq:derivative.converge}, yields
 $
 \int_{I_\delta}\rho u^\infty-\mathfrak H_\lambda(s,u_s^\infty)\varphi ds =-\int_{I_\delta}u_s^\infty\varphi_sds$.
 Thus $u_s^\infty$ has weak derivative $\rho u^\infty-\mathfrak H_\lambda(s,u_s^\infty)$ on $I_\delta$, and
 $
 u_{ss}^\infty =\rho u^\infty-\mathfrak H_\lambda(s,u_s^\infty)\;\text{a.e. on }I_\delta.
 $
 In particular, $u^\infty\in W^{2,{k}}(I_\delta)$ for every finite ${k}$. As $\delta>0$ is arbitrary, $u^\infty\in W_{\rm loc}^{2,{k}}((0,\pi))$ and \eqref{eq:HJB.ulim} holds a.e. on $(0,\pi)$.

We next verify the endpoint behavior needed for the zero-flux condition of $u^\infty$ with respect to $\vpi^\infty$. Repeat the proof of Lemma \ref{lm:individual.endpoint}. The difference is that Corollary \ref{cor:PIA.welldefined}, Theorem~\ref{thm:compact.un}, and Lemma \ref{lm:entropy.est} now make the bound on $\rho u^n-F^{\vpi^n}$ uniform in $n$, i.e., the constant $C_n$ in \eqref{eq:lm.Fnbound} is now replaced by a constant $C$ independent of $n$. The ratio constants in \eqref{eq:lm.bessel.ratio}  are already uniform. We then conclude that
 \be\label{eq:uniform.endpoint}
 \begin{cases}
  |u_s^n(s)|\leq Cs(1+|\ln s|),&0<s\leq \frac{\pi}{8},\\
|u_s^n(s)|\leq C(\pi-s)(1+|\ln(\pi-s)|),&\frac{7\pi}{8}\leq s<\pi.
 \end{cases}
\ee
By Lemma~\ref{lm:endpoint.crossflux} and \eqref{eq:uniform.endpoint}, for an arbitrary feedback policy $\widetilde{\vpi}$,
\bee
  \widehat m^{\widetilde{\vpi}} u_s^n\to0  \text{ as } s\to 0+\text{ and }s\to \pi-, \text{ uniformly in }n.
\eee
Then \eqref{eq:derivative.converge} and  \eqref{eq:uniform.endpoint} imply the same endpoint decay for $u_s^\infty$, that is,
\be\label{eq:uinfty.zeroflux}
\lim_{s\to 0+}\widehat m^{\widetilde\vpi}(s)u_s^\infty(s)=\lim_{s\to\pi-}\widehat m^{\widetilde\vpi}(s)u_s^\infty(s)=0,\quad \forall \widetilde\vpi\in\Ac.
\ee
Lemma~\ref{lm:entropy.est}, the bound on $u_s^\infty$, and \eqref{eq:general-masslog} give $F^{\vpi^\infty}\in L^k(\widehat m^{\vpi^\infty}ds)$ for every finite $k$, so $\vpi^\infty\in\Ac_\lambda$. By the Gibbs identity, \eqref{eq:HJB.ulim} is exactly the policy-evaluation equation for $\vpi^\infty$. Proposition \ref{prop:PDE.wellposed} 
then give $ u^\infty=W^{\vpi^\infty}_\lambda.$

We then show the optimality. For any admissible policy pair $\widetilde{\vpi}\in\Ac_\lambda$, the supremum in the HJB equation implies
 \begin{equation}\label{eq:verification-sub}
  -\rho u^\infty+u_{ss}^\infty+Q^{\widetilde{\vpi}} u_s^\infty  +F^{\widetilde{\vpi}}\le0.
 \end{equation}
 Let $W^{\widetilde{\vpi}}_\lambda$ be the unique policy-evaluation solution of
 $$
 -\rho W^{\widetilde{\vpi}}_\lambda+(W^{\widetilde{\vpi}}_\lambda)_{ss}+Q^{\widetilde{\vpi}} (W^{\widetilde{\vpi}}_\lambda)_s+F^{\widetilde{\vpi}}=0.
 $$
 By \eqref{eq:uinfty.zeroflux}, $u^\infty$ satisfies the zero-flux condition also with respect to $\widehat m^{\widetilde{\vpi}}$. Set $h=W^{\widetilde{\vpi}}_\lambda-u^\infty$. By subtracting \eqref{eq:verification-sub} from the equation for $W^{\widetilde{\vpi}}_\lambda$, we have
 $
 -\rho h+\frac1{\widehat m^{\widetilde{\vpi}}}(\widehat m^{\widetilde{\vpi}} h_s)_s\geq 0.
 $
 Both $W^{\widetilde{\vpi}}_\lambda$ and $u^\infty$ have zero flux for $\widehat m^{\widetilde{\vpi}}$. Multiplying the above inequality with $h^+\widehat m^{\widetilde{\vpi}}$ then following a similar argument from \eqref{eq:h.zero.flux} to the end of the proof of  Proposition~\ref{prop:policy.improvement} yields
 $
 W^{\widetilde{\vpi}}_\lambda\leq u^\infty.
 $
 Since equality holds for ${\widetilde{\vpi}}=\vpi^\infty$, we reach to \eqref{eq:optimal.value}. The correspondence between the policy and value in Section~\ref{sec:arcsin.transform} gives \eqref{eq:PI.original.optimal} and the original HJB equation. Finally, \eqref{eq:CK.original} and \eqref{eq:PI.coordinates} translate the established convergence to \eqref{eq:PI.both.convergence}.
\end{proof}


\section{Convergence of Vanishing Entropy Regularization}\label{sec:convergence.entropy}

We next study the vanishing entropy convergence. The original unregularized HJB equation for the two-state MFC problem is
\be\label{eq:HJB.origin.x}
\begin{aligned}
	0=&-\rho v(p)+p(1-p)v_{pp}(p)\\
	&+\sup_{\vpi_1,\vpi_2\in\Pc(A)}\Bigg\{p\int_A\Big[-2q^{1,2}(p, a)v_p(p)+r(1,p, a)\Big]\vpi_1(da)\\
	&\qquad\qquad\qquad\qquad +(1-p)\int_A\Big[2q^{2,1}(p, a)v_p(p)+r(2,p, a)\Big]\vpi_2(da) \Bigg\},\quad  p\in (0,1).
\end{aligned}
\ee
For $i=1,2$, define the hard Hamiltonians
$$
\mathfrak H_{0,i}(s,y):=\sup_{\theta\in\Pc(A)}\int_A[Q_i(s,a)y+\rtild_i(s,a)]\theta(da)
=\sup_{a\in A}\{Q_i(s,a)y+\rtild_i(s,a)\},
$$
and set
\be\label{eq:hard.Hamiltonian}
\mathfrak H_0(s,y):=p(s)\mathfrak H_{0,1}(s,y)+(1-p(s))\mathfrak H_{0,2}(s,y).
\ee
Thus the arcsine transform of \eqref{eq:HJB.origin.x} is
\be\label{eq:HJB.origin.s}
-\rho u+u_{ss}+\mathfrak H_0(s,u_s)=0,\quad  s\in (0,\pi).
\ee
By \eqref{eq:hard.Hamiltonian} and the definitions of $Q^\vpi$ and $\rtild_i^{\vpi_i}$, for every $(s,y)\in(0,\pi)\times\R$,
\be\label{eq:hard.policy.sup}
\begin{aligned}
	\mathfrak H_0(s,y)=\sup_{\theta_1,\theta_2\in\Pc(A)}\Bigg\{&p(s)\int_A[Q_1(s,a)y+\rtild_1(s,a)]\theta_1(a) da\\
	&+(1-p(s))\int_A[Q_2(s,a)y+\rtild_2(s,a)]\theta_2(a) da\Bigg\}.
\end{aligned}
\ee

To make the dependence on the entropy parameter explicit, throughout this section we write $u^{*,\lambda}:=u^\infty, \vpi^{*,\lambda}:=\vpi^\infty\in\Ac_\lambda$ for the optimal regularized solution and its Gibbs policy pair in Theorem \ref{thm:convergence}, and $u^{n,\lambda}$ and $\vpi^{n,\lambda}\in\Ac_\lambda$ for the corresponding PIA iterates. We also write $v^{n,\lambda}=V_\lambda^{\vpi^{n,\lambda}}$, then $v^{n,\lambda}=u^{n,\lambda}\circ\Phi$. The transformed policy values $W_\lambda^{\vpi}$ and $W_0^{\vpi}$ are those introduced in Section~2 from the values in \eqref{eq:value.regularized.x} and \eqref{eq:value.unregularized.x}. We use $m^{\theta^{n, \lambda}},\widehat{m}^{\theta^{n, \lambda}}$ for the corresponding policy-dependent weights, and reserve $u^{*,0}$ for the solution of \eqref{eq:HJB.origin.s}. For an admissible policy pair $\vpi\in\Ac$, let $W_0^{\vpi}$ denote the unique zero-flux policy-evaluation solution of
\be\label{eq:policy.eval.lambda0}
-\rho W_0^\vpi+(W_0^\vpi)_{ss}+Q^\vpi(W_0^\vpi)_s+p(s)\rtild_1^{\vpi_1}(s)+(1-p(s))\rtild_2^{\vpi_2}(s)=0.
\ee
By Proposition~\ref{prop:PDE.wellposed}, such a solution is well defined. 

Recall the soft Hamiltonian $\mathfrak H_\lambda$ defined in Theorem \ref{thm:convergence}. We first bound the difference between $\mathfrak H_\lambda$ and $\mathfrak H_0$.

\begin{Lemma}\label{lm:regularize.Hconvergence}
There exists a finite constant $C>0$, independent of $\lambda\in(0,1]$, such that
	\be\label{eq:lambdaconverge.weighted}
	|\mathfrak H_\lambda(s,y)-\mathfrak H_0(s,y)|
	\leq C\lambda\left[1+|\ln\lambda|+\ln(1+|y|)+|\ln(\sin s)|\right],\quad \forall s\in(0,\pi),\ y\in\R.
	\ee
Moreover, for every $R>0$ and $1\leq k<\infty$, there exists a constant $C_{R,k}>0$ such that
\begin{align}
\|\mathfrak H_\lambda(\cdot,y(\cdot))-\mathfrak H_0(\cdot,y(\cdot))\|_{L^k(mds)}
&\leq C_{R,k}\lambda(1+|\ln\lambda|),\label{eq:lambdaconverge.lp}\\
\|\mathfrak H_\lambda(\cdot,y(\cdot))-\mathfrak H_0(\cdot,y(\cdot))\|_{L^k(\widehat m ds)}
&\leq C_{R,k}\lambda(1+|\ln\lambda|),	\label{eq:lambdaconverge.full}
\end{align}
for every $\lambda\in(0,1]$, every $y\in L^\infty((0,\pi))$ with $\|y\|_{L^\infty((0,\pi))}\leq R$, and every Bessel weight $m$ generated by coefficients satisfying \eqref{eq:theta.bounds}, where $\widehat m(s)=m(s)\exp\left(\int_{\pi/2}^s q(t)\,dt\right)$, $q\in L^\infty((0,\pi))$ with $\|q\|_{L^\infty((0,\pi))}\leq M.$ The constant $C_{R,k}$ depends only on $R$, $k$, $M$ and  is independent of $\lambda$ and the choices of $y$, $m$ and $q$.
\end{Lemma}

\begin{proof}
	Recall \eqref{eq:lm.Gidef}, it holds that $\mathfrak H_{0,i}(s,y)=\sup_{a\in A}G_i(s,y,a)$ for $i=1,2$. Note that, for any $\theta\in\Pc(A)$,
	$-\int_A\ln\theta(a)\theta(a)da\leq\ln(\Leb(A))$, we have
	$$
	\mathfrak H_{\lambda,i}(s,y)=\lambda\ln\int_Ae^{G_i(s,y,a)/\lambda}da
	\leq \mathfrak H_{0,i}(s,y)+\lambda\ln(\Leb(A)),\quad (s,y)\in(0,\pi)\times\R, ,\;\text{for $i=1,2$}..
	$$
	For the converse bound, fix $(s,y)\in(0,\pi)\times\R$ and choose $\bar a\in A$ such that
	$G_i(s,y,\bar a)\geq\mathfrak H_{0,i}(s,y)-\lambda$. By virtue of \eqref{eq:lm.entropy.est0}, we get
	$$
	|G_i(s,y,a)-G_i(s,y,\bar a)|\leq C\left(1+\frac{|y|}{\sin^2s}\right)|a-\bar a|,\;\text{for $i=1,2$}.
	$$
	For any $(s,y)\in (0,\pi)\times \R$, take
		$
		R(s,y):=\min\left\{\frac\gamma2,\frac{\lambda}{1+C\left(1+|y|/\sin^2s\right)}\right\}.
		$
		For every $a\in A\cap B_{R(s,y)}(\bar a)$, the above Lipschitz estimate and the choice of $\bar a$ give
		\bee
		G_i(s,y,a)\geq G_i(s,y,\bar a)-\lambda
		\geq\mathfrak H_{0,i}(s,y)-2\lambda.
		\eee
		In light of $0<R(s,y)<\gamma$, Assumption~\ref{ass:action}(v) implies 
		$
		\Leb(A\cap B_{R(s,y)}(\bar a))\geq\iota\Leb(B_1(0)) R(s,y)^\ell.
		$
		Consequently, it holds that
		\bee
		\begin{aligned}
			\mathfrak H_{\lambda,i}(s,y)
			&=\lambda\ln\left(\int_Ae^{G_i(s,y,a)/\lambda}da\right)
			\geq\lambda\ln\left(\int_{A\cap B_{R(s,y)}(\bar a)}e^{G_i(s,y,a)/\lambda}da\right)\\
			&\geq\mathfrak H_{0,i}(s,y)-2\lambda
			+\lambda\ln(\iota\Leb(B_1(0)))+\ell\lambda\ln R(s,y),\quad i=1,2.
		\end{aligned}
		\eee
	As $\lambda\in(0,1]$ and
	$
	|\ln R(s,y)|\leq C+|\ln\lambda|+\ln\left(1+\frac{|y|}{\sin^2s}\right),
	$
	the above inequality implies that
	$$
	\begin{aligned}
		\mathfrak H_{0,i}(s,y)-\mathfrak H_{\lambda,i}(s,y)
		&\leq C\lambda\left[1+|\ln\lambda|+\ln\left(1+\frac{|y|}{\sin^2s}\right)\right]\\
		&\leq C\lambda\left[1+|\ln\lambda|+\ln(1+|y|)+2|\ln(\sin s)|\right], \;\text{for $i=1,2$}..
	\end{aligned}
	$$
	Combining this lower estimate with the preceding upper estimate gives
	\be\label{eq:softmax.individual}
	|\mathfrak H_{\lambda,i}(s,y)-\mathfrak H_{0,i}(s,y)|
	\leq C\lambda\left[1+|\ln\lambda|+\ln(1+|y|)+|\ln(\sin s)|\right] ,\;\text{for $i=1,2$}..
	\ee
	By $0\leq p(s),1-p(s)\leq1$, \eqref{eq:lambdaconverge.weighted} follows from \eqref{eq:softmax.individual}.
	
	Moreover, for $\|y\|_\infty\leq A$, the uniform logarithmic integrability in \eqref{eq:general-masslog}, together with the corresponding endpoint estimates for an arbitrary Bessel weight in the class \eqref{eq:theta.bounds}, gives \eqref{eq:lambdaconverge.lp}.
	Finally, \eqref{eq:measure-equivalence} gives
	$\widehat m\leq e^{M\pi}m$, therefore \eqref{eq:lambdaconverge.full} follows directly from \eqref{eq:lambdaconverge.lp} after enlarging $C_{R,k}$.
\end{proof}

As $\lambda\to0+$,
we next show that \eqref{eq:HJB.origin.s} admits a unique solution satisfying the endpoint derivative estimates \eqref{eq:hard.regular.endpoint}.

\begin{Theorem}\label{thm:vanishing.entropy}
	Let Assumptions \ref{ass:b} and \ref{ass:action} hold, \eqref{eq:HJB.origin.s} admits a unique solution $u^{*,0}\in\Cc^1([0,\pi])\cap W_{\rm loc}^{2,k}((0,\pi))$, for some $k>1$, satisfying the endpoint derivative estimates \eqref{eq:hard.regular.endpoint} (with the constant possibly depending on the solution). Moreover, $u^{*,0}\in\Cc^{1,\alpha}([0,\pi])$, and
	\be\label{eq:W.verification}
	u^{*,0}(s)=\sup_{\vpi\in\Ac}W_0^\vpi(s),\quad \forall s\in[0,\pi],
	\ee
	with $v^{*,0}=u^{*,0}\circ\Phi$ solving \eqref{eq:HJB.origin.x} a.e. on $(0,1)$. For every $0<\beta<\alpha$, it holds that
	\be\label{eq:vanishing.entropy.C1beta}
	\lim_{\lambda\to0+}\|u^{*,\lambda}-u^{*,0}\|_{\Cc^{1,\beta}([0,\pi])}=0,
	\quad\text{equivalently,}\quad
	\lim_{\lambda\to0+}\|v^{*,\lambda}-v^{*,0}\|_{\Cc_K^{1,\beta}([0,1])}=0.
	\ee
\end{Theorem}

\begin{proof}
	{\bf Step 1.} We first show the existence of a solution to \eqref{eq:HJB.origin.s} satisfying \eqref{eq:hard.regular.endpoint}. By Theorem~\ref{thm:compact.un}, the optimal solutions inherit the uniform estimate
	$
	\sup_{0<\lambda\leq1}\|u^{*,\lambda}\|_{\Cc^{1,\alpha}([0,\pi])}\leq A^*.
	$
	Let $\lambda_{j}\to0+$. By the compact embedding $\Cc^{1,\alpha}([0,\pi])\Subset \Cc^{1,\beta}([0,\pi])$, after passing to a subsequence,  there exists some $\bar u$ such that
	$u^{*,\lambda_{j}}\to\bar u$ in $\Cc^{1,\beta}([0,\pi])$ for every $\beta<\alpha$. In particular, it holds that
	\be\label{eq:lambda.derivative.compact}
	\|u_s^{*,\lambda_{j}}-\bar u_s\|_{L^\infty([0,\pi])}\to0.
	\ee
	Fix $\delta\in(0,\pi/2)$ and let $I_\delta=[\delta,\pi-\delta]$. By \eqref{eq:def.Q}, Assumption~\ref{ass:b}, and the fact that $p(s)$ and $1-p(s)$ are bounded away from zero on $I_\delta$, there is $C_\delta<\infty$ such that, for $i=1,2$,
	$
	\sup_{s\in I_\delta, a\in A}|Q_i(s,a)|\leq C_\delta.
	$
	Then, for all $s\in I_\delta$ and $y,z\in\R$,
	$
	|\mathfrak H_0(s,y)-\mathfrak H_0(s,z)|\leq C_\delta|y-z|.
	$
	This together with \eqref{eq:lambdaconverge.weighted} in Lemma \ref{lm:regularize.Hconvergence}, \eqref{eq:lambda.derivative.compact}, $\|u_s^{*,\lambda_{j}}\|_\infty\leq A^*$ and $\inf_{I_\delta}\sin s>0$ yields that
	\be\label{eq:hard.H.limit}
	\begin{aligned}
		&\|\mathfrak H_{\lambda_{j}}(s,u_s^{*,\lambda_{j}})-\mathfrak H_0(s,\bar u_s)\|_{L^\infty(I_\delta)}\\
		&\leq \|\mathfrak H_{\lambda_{j}}(s,u_s^{*,\lambda_{j}})-\mathfrak H_0(s,u_s^{*,\lambda_{j}})\|_{L^\infty(I_\delta)}
		+\|\mathfrak H_0(s,u_s^{*,\lambda_{j}})-\mathfrak H_0(s,\bar u_s)\|_{L^\infty(I_\delta)} \\
		&\leq C_\delta\lambda_{j}\left[1+|\ln\lambda_{j}|+\ln(1+A^*)\right]
		+C_\delta\|u_s^{*,\lambda_{j}}-\bar u_s\|_\infty\to0,
	\end{aligned}
	\ee
	For each ${j}$, Theorem \ref{thm:convergence} tells that
	$
	-\rho u^{*,\lambda_{j}}+u_{ss}^{*,\lambda_{j}}
	+\mathfrak H_{\lambda_{j}}(s,u_s^{*,\lambda_{j}})=0.
	$
	Combining this equality with the uniform convergence of $u^{*,\lambda_{j}}$ and \eqref{eq:hard.H.limit} yields
	\bee
	u_{ss}^{*,\lambda_{j}}\to \rho\bar u-\mathfrak H_0(s,\bar u_s)\quad\text{in }L^\infty(I_\delta).
	\eee
	Then a similar argument as in the proof of Theorem \ref{thm:convergence} using test function $\varphi\in \Cc_c^\infty((\delta, \pi-\delta))$ identifies this limit as the weak derivative of $\bar u_s$ on $I_\delta$.
	Since $\delta>0$ is arbitrary, $\bar u\in W_{\rm loc}^{2,{\color{red}k}}((0,\pi))$ for every finite ${k}$ and $\bar u$ solves \eqref{eq:HJB.origin.s} a.e.
	
	It remains to pass the regular boundary behavior to the limit. Given $\lambda\in(0,1]$, let $\vpi^{*,\lambda}$ be the Gibbs optimizer from Theorem~\ref{thm:convergence}. Then Lemma \ref{lm:entropy.est}, the uniform bound $\|u_s^{*,\lambda}\|_\infty\leq A^*$, Assumption~\ref{ass:action}(iv) and
	$\sup_{0<\lambda\leq1}\lambda(1+|\ln\lambda|)<\infty$ together imply
	\bee
	\begin{aligned}
		&\bigg|\rho u^{*,\lambda}(s)-\left[ p(s)\big[\rtild_1^{\vpi_1^{*,\lambda}}(s)+\lambda\Hc(\vpi_1^{*,\lambda}(s))\big]+(1-p(s))\big[\rtild_2^{\vpi_2^{*,\lambda}}(s)+\lambda\Hc(\vpi_2^{*,\lambda}(s))\big] \right] \bigg|\\
		&\leq C[1+|\ln(\sin s)|]
	\end{aligned}
	\eee
	with $C$ independent of $\lambda\in(0,1]$.
	Then Lemma~\ref{lm:endpoint.crossflux}, applied with the above uniform logarithmic bound, gives the following estimates with a constant $C$ independent of $\lambda\in(0,1]$:
	\be\label{eq:lambda.uniform.endpoint}
	\begin{cases}
		|u_s^{*,\lambda}(s)|\leq Cs(1+|\ln s|),&0<s\leq \frac{\pi}{8},\\
		|u_s^{*,\lambda}(s)|\leq C(\pi-s)(1+|\ln(\pi-s)|),&\frac{7\pi}{8}\leq s<\pi.
	\end{cases}
	\ee
	Passing to the limit in \eqref{eq:lambda.uniform.endpoint} together with \eqref{eq:lambda.derivative.compact} shows that $\bar u$ satisfies \eqref{eq:hard.regular.endpoint}.  
	Moreover, by \eqref{eq:general-envelope} in Lemma \ref{lm:ratio}, these endpoint estimates imply that $\widehat m^\vpi(s)\bar u_s(s)\to0$ as $s\to0+$ and $s\to\pi-$ for every $\vpi\in\Ac$.
	
	{\bf Step 2.}  We next prove the uniqueness of solution with the regularity and endpoint behavior stated in the theorem. Let $u,w$ be two such solutions of \eqref{eq:HJB.origin.s} and set $h=u-w$. Fix $\varepsilon>0$. By \eqref{eq:hard.policy.sup} and a standard measurable $\varepsilon$-selection, there exists an admissible feedback policy pair $\vpi^\varepsilon\in\Ac$ such that
	\be\label{eq:eps.policy.u}
	\begin{aligned}
		Q^{\vpi^\varepsilon}(s)u_s(s)&+p(s)\rtild_1^{\vpi_1^\varepsilon}(s)+(1-p(s))\rtild_2^{\vpi_2^\varepsilon}(s) \geq \mathfrak H_0(s,u_s(s))-\varepsilon
		\quad\text{a.e. on } (0,\pi).
	\end{aligned}
	\ee
	Since $u$ satisfies \eqref{eq:HJB.origin.s}, \eqref{eq:eps.policy.u} implies
	$$
	-\rho u+u_{ss}+Q^{\vpi^\varepsilon}u_s+p\rtild_1^{\vpi_1^\varepsilon}+(1-p)\rtild_2^{\vpi_2^\varepsilon}\geq-\varepsilon, \quad\text{a.e. on }(0,\pi).
	$$
	On the other hand, by \eqref{eq:hard.policy.sup} and the equation for $w$,
	$$
	-\rho w+w_{ss}+Q^{\vpi^\varepsilon}w_s+p\rtild_1^{\vpi_1^\varepsilon}+(1-p)\rtild_2^{\vpi_2^\varepsilon}\leq0, \quad\text{a.e. on } (0,\pi).
	$$
	Taking difference of the above two inequalities yields
	$-\rho h+h_{ss}+Q^{\vpi^\varepsilon}h_s\geq-\varepsilon$ a.e. on $(0,\pi)$.
	Set $h^\eps:=h-\varepsilon/\rho$, then
	\bee
	-\rho h^\eps+\frac{1}{\widehat m^{\vpi^\eps}} \left(\widehat m^{\vpi^\eps} h^\eps_s  \right)_s=-\rho h^\eps+h^\eps_{ss}+Q^{\vpi^\varepsilon}h^\eps_s\geq0, \quad\text{a.e. on } (0,\pi).
	\eee
	By \eqref{eq:hard.regular.endpoint} in Lemma \ref{lm:endpoint.crossflux} and \eqref{eq:general-envelope} in Lemma \ref{lm:ratio}, both $u$ and $w$, and hence $h^\eps$, have zero flux for the full policy weight $\widehat m^{\vpi^\varepsilon}$. Moreover, in view that $u,w\in\Cc^1([0,\pi])$ and the weight has finite mass by \eqref{eq:general-masslog} in Lemma \ref{lm:ratio}, it holds that $h^\eps\in\Sc(\vpi^\varepsilon)$. Applying exactly the positive-part argument after \eqref{eq:hsub} in the proof of Proposition~\ref{prop:policy.improvement} by multiplying the last inequality by $(h^\eps)^+\widehat m^{\vpi^\varepsilon}$ gives $(h^\eps)^+=0$. Therefore
	$u-w\leq\varepsilon/\rho$. Then taking $\varepsilon\to0+$ gives $u\leq w$. By interchanging $u$ and $w$, we have $w\leq u$, which readily implies that the solution is unique in the stated class.
	
	Note that every sequence $\lambda_{j}\to0+$ admits a subsequence converging to a solution satisfying \eqref{eq:hard.regular.endpoint}, uniqueness then shows that every such subsequence limit equals the same function, denoted by $u^{*,0}$. Therefore the whole family converges and \eqref{eq:vanishing.entropy.C1beta} follows. By \eqref{eq:CK.original}, the corresponding convergence in the original coordinate follows once $u^{*,0}=v^{*,0}\circ p$ is established, which is confirmed in Step 3 by \eqref{eq:W.verification} and the relation $W_0^\vpi=V_0^\vpi\circ p$.
	
	{\bf Step 3.}  Finally, we finish the proof by verifying \eqref{eq:W.verification}.
	For any admissible policy pair $\vpi\in\Ac$. By \eqref{eq:hard.policy.sup} and \eqref{eq:HJB.origin.s},
	\be\label{eq:hard.verification.sub}
	-\rho u^{*,0}+u_{ss}^{*,0}+Q^\vpi u_s^{*,0}
	+p\rtild_1^{\vpi_1}+(1-p)\rtild_2^{\vpi_2}\leq0.
	\ee
	The function $W_0^\vpi$ solves \eqref{eq:policy.eval.lambda0}. By virtue of \eqref{eq:m.equivalence}, \eqref{eq:general-envelope} in Lemma \ref{lm:ratio} and \eqref{eq:hard.regular.endpoint}, $u^{*,0}$ has zero flux for $\widehat m^\vpi$, and $W_0^\vpi$ has zero flux for the same weight by its definition as the zero-flux policy-evaluation solution of \eqref{eq:policy.eval.lambda0}. Set $h=W_0^\vpi-u^{*,0}$. Subtracting \eqref{eq:hard.verification.sub} from \eqref{eq:policy.eval.lambda0} gives
	$
	-\rho h+\frac1{\widehat m^\vpi}(\widehat m^\vpi h_s)_s\geq0.
	$
	Then by applying the same positive-part argument in Step 2,
	we have $h\leq0$. That is,
	\be\label{eq:hard.verification.upper}
	W_0^\vpi\leq u^{*,0},\quad \forall\vpi\in\Ac.
	\ee
	We will construct a policy attaining equality. For $i=1,2$,
	we can choose an admissible measurable feedback policy pair $\vpi^*=(\vpi_1^*,\vpi_2^*)$ such that
	\be\label{eq:hard.support.argmax}
	\operatorname{supp}\vpi_i^*(s)\subset\argmax_{a\in A}\left\{Q_i(s,a)u_s^{*,0}(s)+\rtild_i(s,a)\right\},\quad \text{for a.e. $s\in(0,\pi)$. }
	\ee
	In particular, one may choose measurable maximizers $a_i^*(s)\in \argmax_{a\in A}\left\{Q_i(s,a)u_s^{*,0}(s)+\rtild_i(s,a)\right\}$ and take the classical controls $\vpi_i^*(s)=\delta_{a_i^*(s)}$. Noting that the integrands are constant on their respective maximizer sets, we get by \eqref{eq:hard.support.argmax} that
	$$
	\int_A\left[Q_i(s,a)u_s^{*,0}(s)+\rtild_i(s,a)\right]\vpi_i^*(s,da)=\mathfrak H_{0,i}(s,u_s^{*,0}(s)),\quad \text{for } i=1,2.
	$$
	This, together with \eqref{eq:hard.Hamiltonian}, gives
	\bee
	Q^{\vpi^*}(s)u_s^{*,0}(s)+p(s)\rtild_1^{\vpi_1^*}(s)+(1-p(s))\rtild_2^{\vpi_2^*}(s) =\mathfrak H_0(s,u_s^{*,0}(s)).
	\eee
	Therefore, by \eqref{eq:HJB.origin.s}, $u^{*,0}$ satisfies exactly the policy-evaluation equation corresponding to $\vpi^*$:
	\be\label{eq:hard.optimal.policy.PDE}
	-\rho u^{*,0}+u_{ss}^{*,0}+Q^{\vpi^*}u_s^{*,0}+p\rtild_1^{\vpi_1^*}+(1-p)\rtild_2^{\vpi_2^*}=0, \quad \text{a.e. on } (0,\pi).
	\ee
	We next verify that $u^{*,0}$ is the zero-flux policy-evaluation solution of \eqref{eq:hard.optimal.policy.PDE}. By Assumption~\ref{ass:b}, averaging the endpoint coefficients against $\vpi_1^*$ and $\vpi_2^*$ preserves the bounds \eqref{eq:main.boundaryrelaxed}. Also \eqref{eq:assum.remainder} in Assumption \ref{ass:b} ensures that the corresponding transformed remainder is uniformly bounded. Hence $Q^{\vpi^*}$ belongs to the same uniformly controlled Bessel class as the policy drifts considered in Section \ref{sec:Bessel.estimates}. As $u^{*,0}$ satisfies the endpoint decay \eqref{eq:hard.regular.endpoint}, it follows from \eqref{eq:m.equivalence} and \eqref{eq:general-envelope} in Lemma \ref{lm:ratio} that
	$$
	\widehat m^{\vpi^*}(s)u_s^{*,0}(s)\to0
	\quad\text{as }s\to0+\text{ and }s\to\pi-.
	$$
	Thus $u^{*,0}$ satisfies the zero-flux condition for \eqref{eq:hard.optimal.policy.PDE}. Consequently, the uniqueness argument in Proposition~\ref{prop:PDE.wellposed} can be invoked to this unregularized policy evaluation equation, which yields that  
	$
	u^{*,0}=W_0^{\vpi^*}
	$
	is the unique zero-flux policy-evaluation solution to \eqref{eq:hard.optimal.policy.PDE}. This, together with \eqref{eq:hard.verification.upper}, gives
	$
	u^{*,0}=W_0^{\vpi^*}=\sup_{\vpi\in\Ac}W_0^\vpi,
	$
	which proves \eqref{eq:W.verification}.
\end{proof}

\begin{Remark}\label{rmk:hard.attainment}
	The verification argument only requires that the support of $\vpi_i^*(s)$ be contained in $\argmax_{a\in A}\left\{Q_i(s,a)u_s^{*,0}(s)+\rtild_i(s,a)\right\}$. Thus no density is needed when $\lambda=0$. A classical feedback maximizer $a_i^*(s)$ corresponds to the Dirac relaxed control $\delta_{a_i^*(s)}$, while any probability measure supported on $\argmax_{a\in A}\left\{Q_i(s,a)u_s^{*,0}(s)+\rtild_i(s,a)\right\}$ gives the same value of the unregularized Hamiltonian. In particular, when the pointwise maximizers admit measurable selectors, the unregularized optimum is attained by a classical feedback control pair.
\end{Remark}

With the above preparation, we are ready to provide a quantitative regret analysis for the entropy regularization.
\begin{Theorem}\label{thm:lambda.rate}
	Under Assumptions \ref{ass:b} and \ref{ass:action}, there exists $C$ independent of $\lambda\in(0,1]$ such that
	\be\label{eq:lambda.rate}
	\|u^{*,\lambda}-u^{*,0}\|_{\Cc^{1,\alpha}([0,\pi])}\leq C\lambda(1+|\ln\lambda|),\;\text{equivalently, } \|v^{*,\lambda}-v^{*,0}\|_{\Cc_K^{1,\alpha}([0,1])} \leq C\lambda(1+|\ln\lambda|).
	\ee
\end{Theorem}

\begin{proof}
	Set $h^\lambda:=u^{*,\lambda}-u^{*,0}$ and define
	$$
	b_i^\lambda(s):=
	\begin{cases}
		\frac{\mathfrak H_{0,i}(s,u_s^{*,\lambda}(s))-\mathfrak H_{0,i}(s,u_s^{*,0}(s))}{u_s^{*,\lambda}(s)-u_s^{*,0}(s)},\quad h_s^\lambda(s)\neq0,\\
\frac{1}{\operatorname{Leb}(A)}\int_A Q_i(s,a)da,\qquad\quad \; h_s^\lambda(s)=0.
	\end{cases}
	$$
	We first show that
	\be\label{eq:secant.slope.interval}
	\inf_{a\in A}Q_i(s,a)\leq b_i^\lambda(s)\leq\sup_{a\in A}Q_i(s,a),\qquad i=1,2.
	\ee
	The above inequalities hold obviously on $\{h_s^\lambda=0\}$. Now take $y_0, y_i\in \R$ such that $y_0<y_1$. For every \(a\in A\),
	$$
	\begin{aligned}
		Q_i(s,a)y_1+\widetilde r_i(s,a)
		&=Q_i(s,a)y_0+\widetilde r_i(s,a)+Q_i(s,a)(y_1-y_0)\\
		&\leq\mathfrak H_{0,i}(s,y_0)+\left(\sup_{\bar a\in A}Q_i(s,\bar a)\right)(y_1-y_0).
	\end{aligned}
	$$
	Taking the supremum over $a$ on the left-hand side of the above inequality, we have that
	$$
	\mathfrak H_{0,i}(s,y_1)\leq \mathfrak H_{0,i}(s,y_0)+ \left(\sup_{a\in A}Q_i(s,a)\right)(y_1-y_0).
	$$
	That is,
	$
	\frac{\mathfrak H_{0,i}(s,y_1)-\mathfrak H_{0,i}(s,y_0)}{y_1-y_0} \leq\sup_{a\in A}Q_i(s,a).
	$
	Now for any $a\in A$,
	$$
	Q_i(s,a)y_0+\widetilde r_i(s,a)=Q_i(s,a)y_1+\widetilde r_i(s,a)-Q_i(s,a)(y_1-y_0).
	$$
	Then, by $Q_i(s,a)\geq \inf_{\bar a\in A}Q_i(s,\bar a)$ and $y_1-y_0>0$,
	$$
	Q_i(s,a)y_0+\widetilde r_i(s,a)\leq\mathfrak H_{0,i}(s,y_1)-\left(\inf_{\bar a\in A}Q_i(s,\bar a)\right)(y_1-y_0).
	$$
	Taking the supremum over $a$ on the left-hand side above implies
	$
	\frac{\mathfrak H_{0,i}(s,y_1)-\mathfrak H_{0,i}(s,y_0)}{y_1-y_0} \geq \inf_{a\in A}Q_i(s,a).
	$
	In sum, for $y_0<y_1$,
	$$
	\inf_{a\in A}Q_i(s,a)\leq
	\frac{\mathfrak H_{0,i}(s,y_0)-\mathfrak H_{0,i}(s,y_1)}{y_0-y_1}
	\leq \sup_{a\in A}Q_i(s,a).
	$$
	The case $y_0>y_1$ follows by interchanging $y_0$ and $y_1$. Hence, \eqref{eq:secant.slope.interval} holds on $\{ h^\lambda_s\neq 0\}$, and therefore on all of $(0,\pi)$.

	Let us define $B^\lambda(s):=p(s)b_1^\lambda(s)+(1-p(s))b_2^\lambda(s)$, and verify that $B^\lambda$ belongs to a uniformly controlled Bessel class.  By  \eqref{eq:secant.slope.interval}, $b_i^\lambda(s)$ belongs to the closed convex hull of the scalar set $\{Q_i(s,a):a\in A\}$. Equivalently, by a measurable convex-selection argument, for $i=1,2,$ we can take a measurable function $\nu_i: \{ h^\lambda_s \neq 0\}\to \Pc(A)$ such that $Q^{\nu_i(s)}(s)=b_i^\lambda(s)$. 
	And we extend $\nu^i$ to $(0,\pi)$ by simply setting $\nu_i=\vpi^0$ on $\{ h^\lambda_s =0 \}$. From \eqref{eq:def.Q} and the definitions of $q_1,q_2$ in \eqref{eq:def.q}
	By setting
	$$
	\begin{aligned}
		&\Lambda_1^\lambda(s):=2\int_Aq^{1,2}(1, a)\nu_1^{\lambda,s}(da), \; \Lambda_2^\lambda(s):=2\int_Aq^{2,1}(0, a)\nu_2^{\lambda,s}(da),\\
		&q^\lambda(s):=\int_Aq_1(s,a)\nu_1^{\lambda,s}(da)+\int_Aq_2(s,a)\nu_2^{\lambda,s}(da),
	\end{aligned}\
	$$
	we have
	$
	B^\lambda(s)=\left(\Lambda_2^\lambda(s)-\frac12\right)\cot\frac{s}{2}-\left(\Lambda_1^\lambda(s)-\frac12\right)\tan\frac{s}{2}+q^\lambda(s).
	$
	Then \eqref{eq:main.boundaryrelaxed} in Assumption \ref{ass:b} gives $\thl\leq\Lambda_i^\lambda(s)\leq\thu$, and \eqref{eq:assum.remainder} together with the definitions of $q_1(s,a),q_2(s,a)$ give $|q_1(s,a)|\vee |q_2(s,a)|\leq M$, hence $|q^\lambda(s)|\leq2M$. Therefore, $B^\lambda$ is of the form required in Section \ref{sec:PIA}, with endpoint coefficients in the fixed interval $[\thl,\thu]$ and a remainder bounded uniformly over $\lambda\in(0,1]$.
	
	Taking the difference between \eqref{eq:HJB.origin.s} for $u^{*,0}$ and \eqref{eq:HJB.ulim} for $u^{*,\lambda}$, we obtain
	$
	\rho h^\lambda-h_{ss}^\lambda=\mathfrak H_\lambda(s,u_s^{*,\lambda})-\mathfrak H_0(s,u_s^{*,0}).
	$
	By the definition of $B^\lambda$, $\mathfrak H_0(s,u_s^{*,\lambda})-\mathfrak H_0(s,u_s^{*,0})=B^\lambda(s)h_s^\lambda(s),$ where the equality is also valid on $\{h_s^\lambda=0\}$ because both sides vanish. As a result,
	\be\label{eq:lambda.error.linear}
	\rho h^\lambda
	-\left(h_{ss}^\lambda+B^\lambda h_s^\lambda\right)=
	\mathfrak H_\lambda(s,u_s^{*,\lambda})-\mathfrak H_0(s,u_s^{*,\lambda}).
	\ee
	
	Let $m^{B^\lambda}$ and $\widehat m^{B^\lambda}$ be the singular and full integrating factors associated with the decomposition just established. Thanks to \eqref{eq:lambda.uniform.endpoint}, \eqref{eq:hard.regular.endpoint}, \eqref{eq:m.equivalence},  \eqref{eq:general-envelope} in Lemma~\ref{lm:ratio} applied to the remainder bound $M$, $h^\lambda$ satisfies the zero-flux condition for $\widehat m^{B^\lambda}$. Theorem \ref{thm:compact.un} gives $\|u_s^{*,\lambda}\|_\infty\leq A^*$ uniformly over $\lambda\in(0,1]$, then \eqref{eq:lambdaconverge.lp} in Lemma \ref{lm:regularize.Hconvergence}, applied to the Bessel weight $m^{B^\lambda}$, yields
	\be\label{eq:lambda.error.source}
	\|\mathfrak H_\lambda(s,u_s^{*,\lambda})-\mathfrak H_0(s,u_s^{*,\lambda})\|_{L^k(m^{B^\lambda}ds)}\leq C\lambda(1+|\ln\lambda|).
	\ee
	By \eqref{eq:lambdaconverge.full}, the same estimate holds with $m^{B^\lambda}$ replaced by $\widehat m^{B^\lambda}$, with a possibly different constant $C$ independent of $\lambda\in(0,1]$.
	
	We next obtain the $L^k$ estimate for $h^\lambda$ from \eqref{eq:lambda.error.linear}. In divergence form,
	$$
	\rho h^\lambda-\frac1{\widehat m^{B^\lambda}}
	(\widehat m^{B^\lambda}h_s^\lambda)_s=\mathfrak H_\lambda(s,u_s^{*,\lambda})-\mathfrak H_0(s,u_s^{*,\lambda}).
	$$
		Multiplying by $|h^\lambda|^{{k}-2}h^\lambda\widehat m^{B^\lambda}$ and integrating over $(0,\pi)$, the zero-flux condition eliminates the boundary terms and gives
	$$
	\rho\int_0^\pi|h^\lambda|^{k}\widehat m^{B^\lambda}ds+({k}-1)\int_0^\pi|h^\lambda|^{{k}-2}|h_s^\lambda|^2\widehat m^{B^\lambda}ds
	=\int_0^\pi \left[ \mathfrak H_\lambda(s,u_s^{*,\lambda})-\mathfrak H_0(s,u_s^{*,\lambda}) \right] |h^\lambda|^{{k}-2}h^\lambda\widehat m^{B^\lambda}ds.
	$$
	For $1<{k}<2$, the identity follows by the standard approximation
	$|h|^{{k}-2}h\mapsto(h^2+\varepsilon)^{({k}-2)/2}h$ and then letting $\varepsilon\to0+$. Utilizing  H\"older inequality yields
	$$
	\|h^\lambda\|_{L^{k}(\widehat m^{B^\lambda}ds)}\leq\frac1\rho\|\mathfrak H_\lambda(s,u_s^{*,\lambda})-\mathfrak H_0(s,u_s^{*,\lambda})\|_{L^{k}(\widehat m^{B^\lambda}ds)}\leq C\lambda(1+|\ln\lambda|).
	$$
	Using the uniform equivalence \eqref{eq:measure-equivalence} between $m^{B^\lambda}$ and $\widehat m^{B^\lambda}$, we also have
	\be\label{eq:lambda.error.Lp}
	\|h^\lambda\|_{L^{k}(m^{B^\lambda}ds)}
	\leq C\lambda(1+|\ln\lambda|).
	\ee
	Then applying the uniform Bessel Sobolev estimate in Theorem \ref{thm:main.besselest}(ii) to \eqref{eq:lambda.error.linear}, then combining with the remainder bound just obtained for $B^\lambda$, together with \eqref{eq:lambda.error.source} and \eqref{eq:lambda.error.Lp}, we arrive at
	$$
	\|h^\lambda\|_{W^{2,{k}}(m^{B^\lambda}ds)}\leq C\lambda(1+|\ln\lambda|).
	$$
	Applying the uniform embedding in Theorem~\ref{thm:main.besselest}(i) yields \eqref{eq:lambda.rate}. The estimate in the original variable follows from \eqref{eq:CK.original}.
\end{proof}

\begin{Corollary}\label{cor:entropy.optimal.regret}
	The entropy-optimal policy pair $\vpi^{*,\lambda}$ is $O(\lambda(1+|\ln\lambda|))$-optimal for the unregularized problem in the sense that
	$$
	0\leq u^{*,0}(s)-W_0^{\vpi^{*,\lambda}}(s) \quad \forall s\in(0,\pi),\quad\text{and }
	\|u^{*,0}-W_0^{\vpi^{*,\lambda}}\|_{\Cc^{1,\alpha}([0,\pi])}\leq C\lambda(1+|\ln\lambda|).
	$$
	The same estimate holds in $\Cc_K^{1,\alpha}([0,1])$ after the inverse arcsine transformation. 
\end{Corollary}

\begin{proof}
	The pointwise nonnegativity follows from \eqref{eq:W.verification}. By \eqref{eq:optimal.value}, $u^{*,\lambda}=W_\lambda^{\vpi^{*,\lambda}}$.
	We first compare the regularized and unregularized evaluations of the same policy $\vpi^{n,\lambda}$. By definition,
	$u^{n,\lambda}=W_\lambda^{\vpi^{n,\lambda}}$, while $W_0^{\vpi^{n,\lambda}}$ solves \eqref{eq:policy.eval.lambda0}. Hence
	$d^{n,\lambda}:=u^{n,\lambda}-W_0^{\vpi^{n,\lambda}}$ satisfies
	\be\label{eq:same.policy.difference.PDE}
	\rho d^{n,\lambda}-\left(d_{ss}^{n,\lambda}+Q^{\vpi^{n,\lambda}}d_s^{n,\lambda}\right)
	=\lambda\left[p\Hc(\vpi_1^{n,\lambda})+(1-p)\Hc(\vpi_2^{n,\lambda})\right].
	\ee
	Both terms in the difference are regular zero-flux solutions for the same policy drift $Q^{\vpi^{n,\lambda}}$, so $d^{n,\lambda}$ satisfies the zero-flux condition for $\widehat{m}^{\theta^{n, \lambda}} $.
	
	For $n\geq1$, the policy update in Section~4 gives
	$\vpi_i^{n,\lambda}=\Gamma_{\lambda,i}(s,u_s^{n-1,\lambda},\cdot)$, $i=1,2$. The uniform estimate in Theorem~\ref{thm:compact.un} and the entropy estimate \eqref{eq:entropy-pointwise} in Lemma~\ref{lm:entropy.est} imply
	$$
	\lambda|\Hc(\vpi_i^{n,\lambda}(s))|
	\leq C\lambda\left[1+|\ln\lambda|+\ln(1+A^*)+|\ln(\sin s)|\right],\qquad i=1,2.
	$$
	For $n=0$, $\vpi_i^{0,\lambda}$ is the uniform density and
	$\Hc(\vpi_i^{0,\lambda})=\ln(\Leb(A))$. Therefore, applying \eqref{eq:general-masslog} in Lemma~\ref{lm:ratio} and the measure equivalence \eqref{eq:m.equivalence}, uniformly in $n\geq0$ and $\lambda\in(0,1]$, leads to
	\be\label{eq:same.policy.entropy.source}
	\left\|\lambda\left[p\Hc(\vpi_1^{n,\lambda})+(1-p)\Hc(\vpi_2^{n,\lambda})\right]\right\|_{L^p(\widehat{m}^{\theta^{n, \lambda}} ds)}
	\leq C\lambda(1+|\ln\lambda|).
	\ee
	
	Applying the $L^p$ resolvent estimate in Proposition~\ref{prop:PDE.wellposed} to \eqref{eq:same.policy.difference.PDE} gives
	$$
	\|d^{n,\lambda}\|_{L^p(\widehat{m}^{\theta^{n, \lambda}} ds)}
	\leq\frac1\rho
	\left\|\lambda\left[p\Hc(\vpi_1^{n,\lambda})+(1-p)\Hc(\vpi_2^{n,\lambda})\right]\right\|_{L^p(\widehat{m}^{\theta^{n, \lambda}} ds)}.
	$$
	By \eqref{eq:m.equivalence}, we pass from $\widehat{m}^{\theta^{n, \lambda}} $ to $m^{\theta^{n, \lambda}}$ and apply Theorem \ref{thm:main.besselest}(ii) to \eqref{eq:same.policy.difference.PDE}, then we obtain
	$
	\|d^{n,\lambda}\|_{W^{2,p}(m^{\theta^{n, \lambda}} ds)}\leq C\lambda(1+|\ln\lambda|),
	$
	where $C$ is independent of $n$ and $\lambda$. The embedding in Theorem~\ref{thm:main.besselest}(i) then yields
	\be\label{eq:same.policy.entropy.error}
	\|u^{n,\lambda}-W_0^{\vpi^{n,\lambda}}\|_{\Cc^{1,\alpha}([0,\pi])}
	\leq C\lambda(1+|\ln\lambda|).
	\ee
	Taking $n\to \infty$ further gives
	$$\|u^{*,\lambda}-W_0^{\vpi^{*,\lambda}}\|_{\Cc^{1,\alpha}([0,\pi])}\leq C\lambda(1+|\ln\lambda|).$$
	Combining this estimate with Theorem~\ref{thm:lambda.rate} proves the desired norm bound.
\end{proof}

\section{An Example of Epidemic Containment}\label{sec:epidemic}

In this section, we apply the exploratory two-state mean-field control problem to epidemic containment, considering  an SIS-type model  infectious disease in a large population of agents.

We consider two states: infected (State 1) and susceptible (State 2).  In the SIS model, recovered individuals return to the susceptible state and may become reinfected. Let \(p \in [0,1]\) denote the fraction of the infected population. The social planner implements distinct relaxed feedback policies, \(\theta_1(p,a)\) and \(\theta_2(p,a)\), to control the social activity levels \(a \in [\underline{a}, \bar{a}]\) of the infected and susceptible individuals, respectively. Relaxed feedback policies are well suited to capturing uncertainty in individual behavioral responses to epidemic interventions through randomized actions. Inspired by the discussion in \cite{Lefevre1981} regarding large populations, we specify the transition rates as follows:
\begin{equation}
    q^{2,1}(p, a) = \beta a p^\delta + \epsilon_0 , \quad
    q^{1,2}(p, a) = \gamma,
\end{equation}
where \(\delta > 0\) is the infection power, \(\beta > 0\) is the transmission rate, \(\gamma > 0\) is the constant recovery rate, and \(\epsilon_0 > 0\) is a small external infection rate. To account for aggregate uncertainty arising from factors such as seasonal variations and changes in social interactions, we incorporate a Kimura-type common noise into the epidemic dynamics. The diffusion coefficient  vanishes at the boundaries, making it a natural choice for modeling stochastic fluctuations of the infected fraction.

The social planner aims to maximize the expected discounted social welfare. Following \cite{Aurell2022}, the individual running rewards are modeled as negative quadratic penalties for deviating from the recommended activity levels \(a_i^{\mathrm{rec}}(p)\). Additionally, an aggregate social cost \(c_H P_t\) for the infected agents reduces the social welfare. Thus, the resulting objective functional is
\begin{equation*}
\begin{aligned}
V^{\boldsymbol{\pi}}(p) = \mathbb{E}_p^{\boldsymbol{\pi}} \Bigg[ \int_0^\infty e^{-\rho t} \Bigg\{ & - c_H P_t + P_t \int_A \Big[ - \frac{c_{a_1}}{2}\big(a_1^{\mathrm{rec}}(P_t) - a\big)^2  - \lambda \ln \theta_1(a|P_t) \Big] \theta_1(a|P_t) da \\
& + (1-P_t) \int_A \Big[ - \frac{c_{a_2}}{2}\big(a_2^{\mathrm{rec}}(P_t) - a\big)^2 - \lambda \ln \theta_2(a|P_t) \Big] \theta_2(a|P_t) da \Bigg\} dt \Bigg],
\end{aligned}
\end{equation*}
where \(c_{a_1}, c_{a_2}, c_H > 0\) are penalty coefficients, and \(\lambda > 0\) is the entropy regularization parameter.

We note that when \(\delta \ge \frac{1}{2}\) and \(\epsilon_0 > 0\), the proposed model naturally satisfies Assumptions~\ref{ass:b} and~\ref{ass:action}. Specifically, \(\epsilon_0 > 0\) guarantees a strictly positive inward drift at \(p=0\). Furthermore, \(\delta \ge \frac{1}{2}\) ensures ~\eqref{eq:assum.remainder} holds in Assumption~\ref{ass:b}. Since \([\underline a, \bar a]\) is bounded and closed, the conditions in Assumption~\ref{ass:action} hold trivially.

As all the conditions are fulfilled, we use the PIA to find the optimal value function and the optimal policy.

\paragraph{Numerical implementations}
We implement the PIA in both the
original and transformed spaces. We set the parameters \(\rho = 0.1\), \(\beta = 0.2\), \(\delta = 1.0\), \(\epsilon_0 = 0.01\), \(\gamma = 0.1\), \(c_H = 5\), \(c_{a1} = 2\), \(c_{a2} = 10\), \(a_1^{\mathrm{rec}}(p) = 0\), \(a_2^{\mathrm{rec}}(p) = 1 - 0.8p\), \(\underline{a} = 0\), \(\bar{a} = 1\), and \(\lambda = 0.1\).

In the original space, the state variable $p\in[0,1]$ is discretized on a uniform
grid with \(N_p = 2000\) intervals. The policy evaluation step is solved using a finite
difference scheme with directional discretization of the
drift terms.  In the transformed space, we discretize the transformed equation on a uniform
$s$-grid with \(N_s = 2000\) intervals using central finite differences at interior
nodes. At the endpoints, we use the boundary equations
inherited from the degenerate generator in the original
space. In both formulations, the policy improvement step is performed
using the Gibbs formula. The resulting tridiagonal systems are solved numerically
at each policy evaluation step.The action space \([0,1]\) is discretized using \(N_a = 200\) intervals, with trapezoidal quadrature for the integrals.

For comparison, we also solve the corresponding discrete
optimal HJB equations in both coordinates to obtain
reference value functions for PIA.

\begin{figure}[htbp]
  \centering
  \begin{subfigure}[t]{0.24\linewidth}
    \centering
    \includegraphics[width=\linewidth]{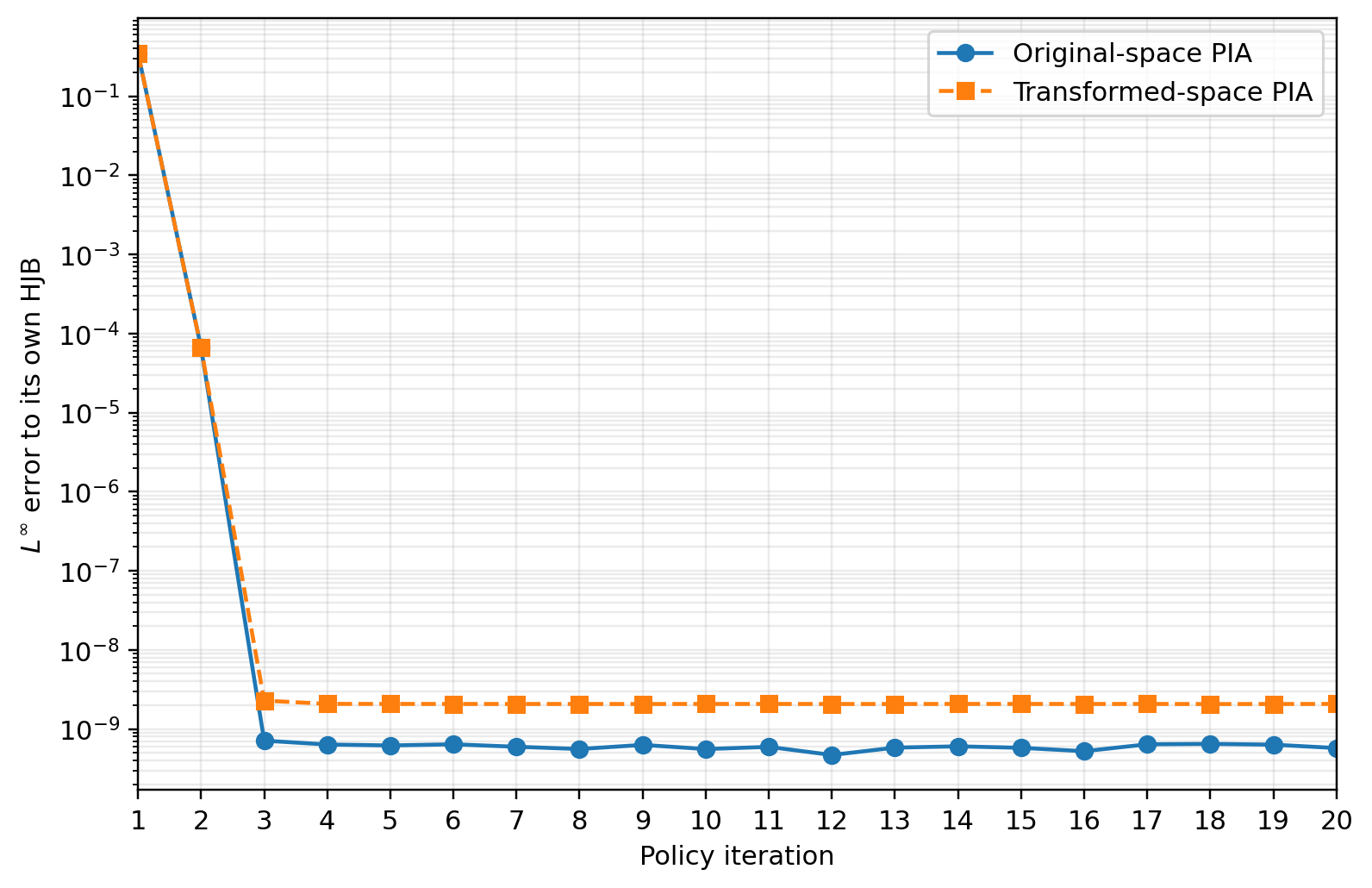}
    \caption{Value convergence}
    \label{fig:Fig1}
  \end{subfigure}\hfill
  \begin{subfigure}[t]{0.24\linewidth}
    \centering
    \includegraphics[width=\linewidth]{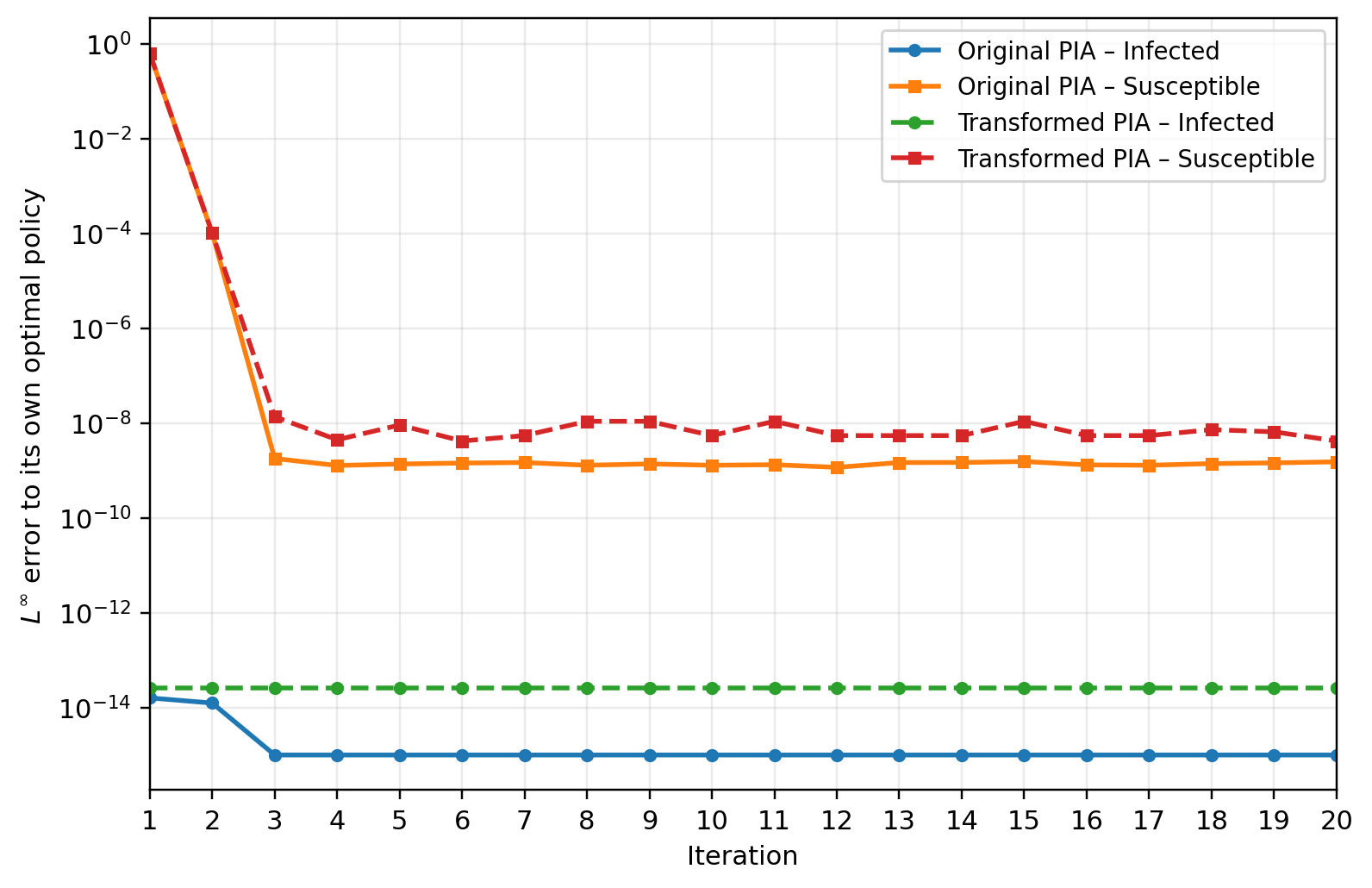}
    \caption{Policy convergence}
    \label{fig:Fig2}
  \end{subfigure}\hfill
  \begin{subfigure}[t]{0.24\linewidth}
    \centering
    \includegraphics[width=\linewidth]{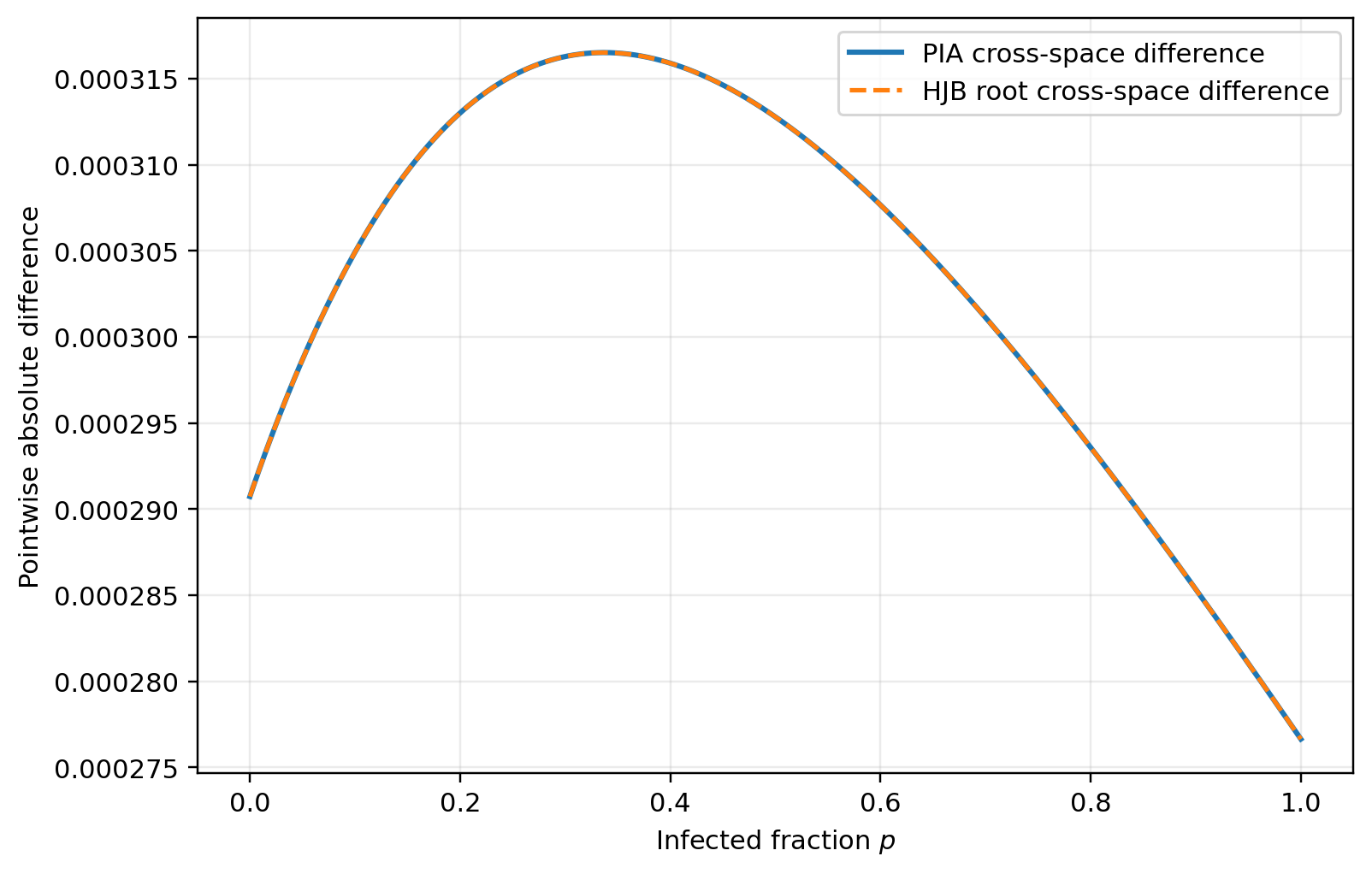}
    \caption{Cross discrepency}
    \label{fig:Fig3}
  \end{subfigure}\hfill
  \begin{subfigure}[t]{0.24\linewidth}
    \centering
    \includegraphics[width=\linewidth]{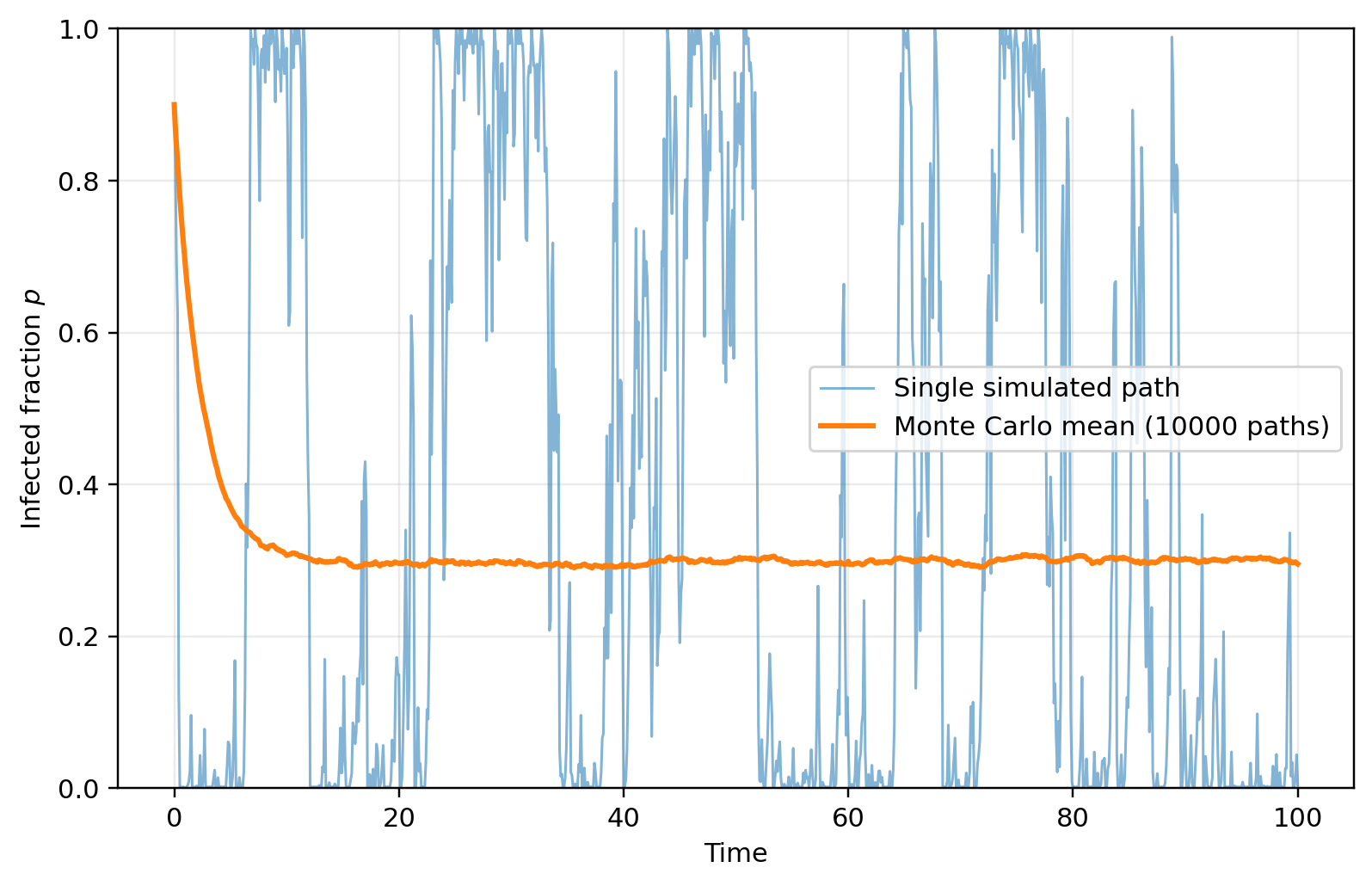}
    \caption{Population dynamics}
    \label{fig:Fig4}
  \end{subfigure}

  \caption{Performance of the policy iteration algorithm}
  \label{fig:Fig1-4}
\end{figure}

\paragraph{Numerical results} We first investigate the convergence of the PIA, as illustrated in Figure~\ref{fig:Fig1-4}. Figure~\ref{fig:Fig1} presents the $L^\infty$ errors
$\|V^n-V^*\|_\infty$
and 
$\|u^n-u^*\|_\infty$,
where $V^n$ and $u^n$ denote the value-function iterates
in the original and transformed spaces, respectively, and $V^*$ and $u^*$ are the corresponding HJB solutions. Figure~\ref{fig:Fig2} illustrates the $L^\infty$ convergence of the relaxed policies for both infected and susceptible populations towards their respective optimal policies. Both algorithms exhibit rapid convergence during the
first few iterations. In Figure~\ref{fig:Fig3}, we compare the optimal value functions obtained from the HJB equations in two coordinates systems by plotting the pointwise discrepancy  $|V_{\mathrm{trans}}^*(p) - V^*(p)|$ throughout $p \in [0, 1]$, where $V_{\mathrm{trans}}^*(p) = u^*(2\arcsin\sqrt{p})$. The cross-coordinate discrepancy in Figure~\ref{fig:Fig3} is larger than the PIA errors in Figure~\ref{fig:Fig1}, as it reflects differences in spatial grids and finite-difference discretizations rather than errors arising from policy iteration.  Consequently, a cross-coordinate discrepancy may remain even when both PIA algorithms have converged. Figure~\ref{fig:Fig4} shows the evolution of the optimally controlled infected fraction with a single sample path and the Monte Carlo mean. The latter approaches a positive level around $0.28$ rather than zero, reflecting the possibility of reinfection in the
SIS framework and the presence of a positive baseline infection rate $\epsilon_0 >0$.

\begin{figure}[htbp]
  \centering
  \begin{subfigure}[t]{0.24\linewidth}
    \centering
    \includegraphics[width=\linewidth]{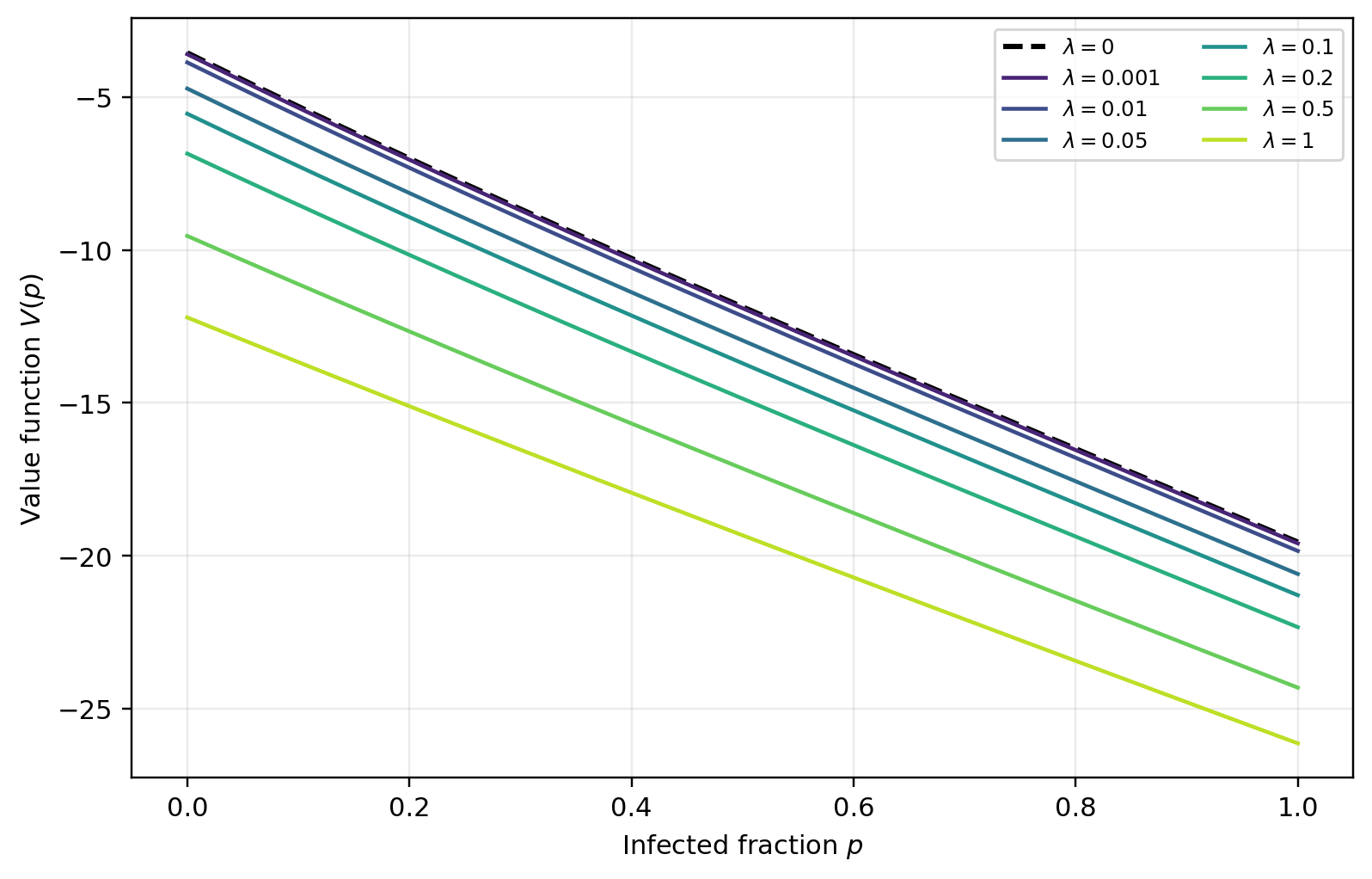}
    \caption{Original values}
    \label{fig:Fig5}
  \end{subfigure}\hfill
  \begin{subfigure}[t]{0.24\linewidth}
    \centering
    \includegraphics[width=\linewidth]{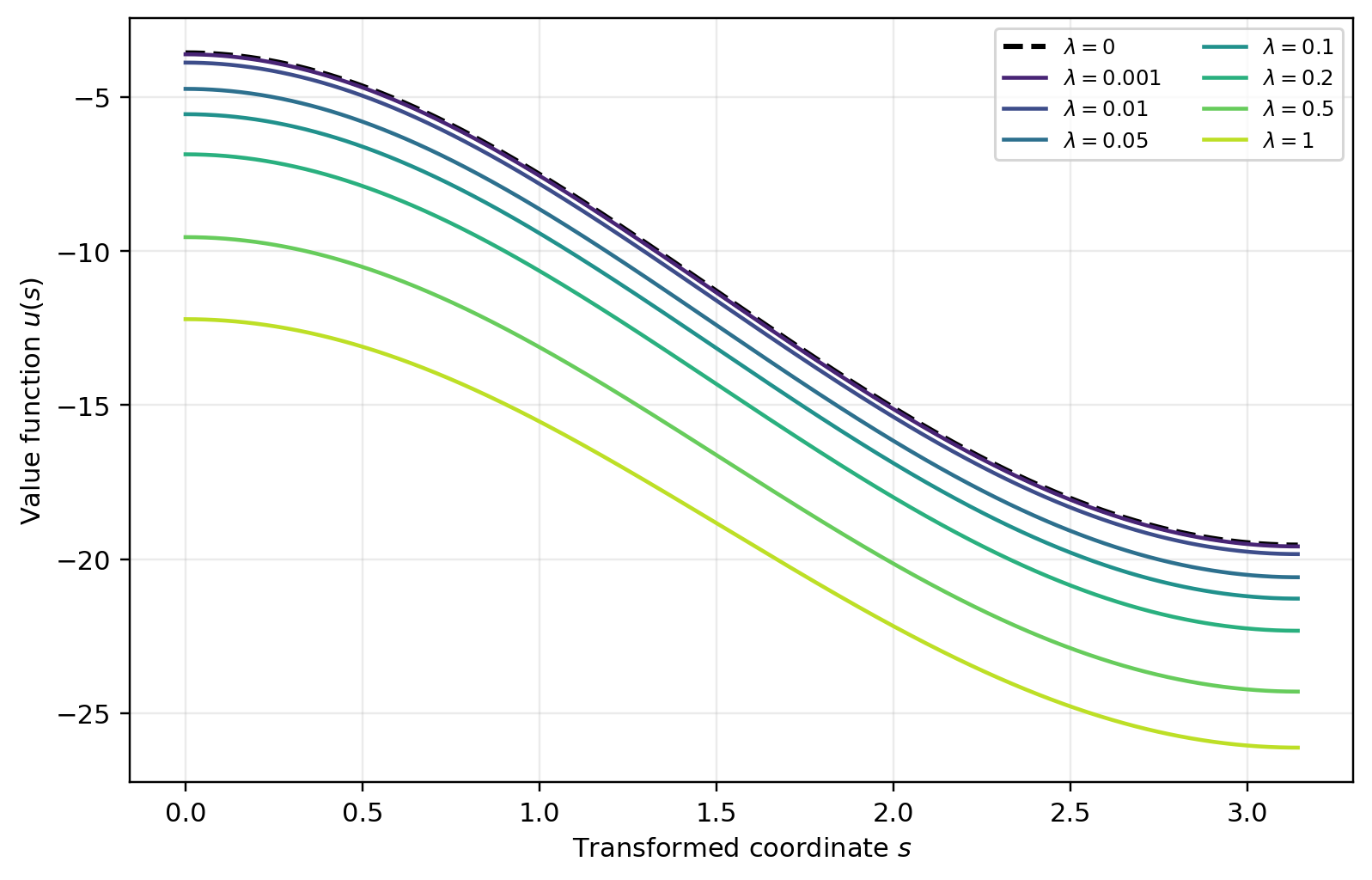}
    \caption{Transformed values}
    \label{fig:Fig6}
  \end{subfigure}\hfill
  \begin{subfigure}[t]{0.24\linewidth}
    \centering
    \includegraphics[width=\linewidth]{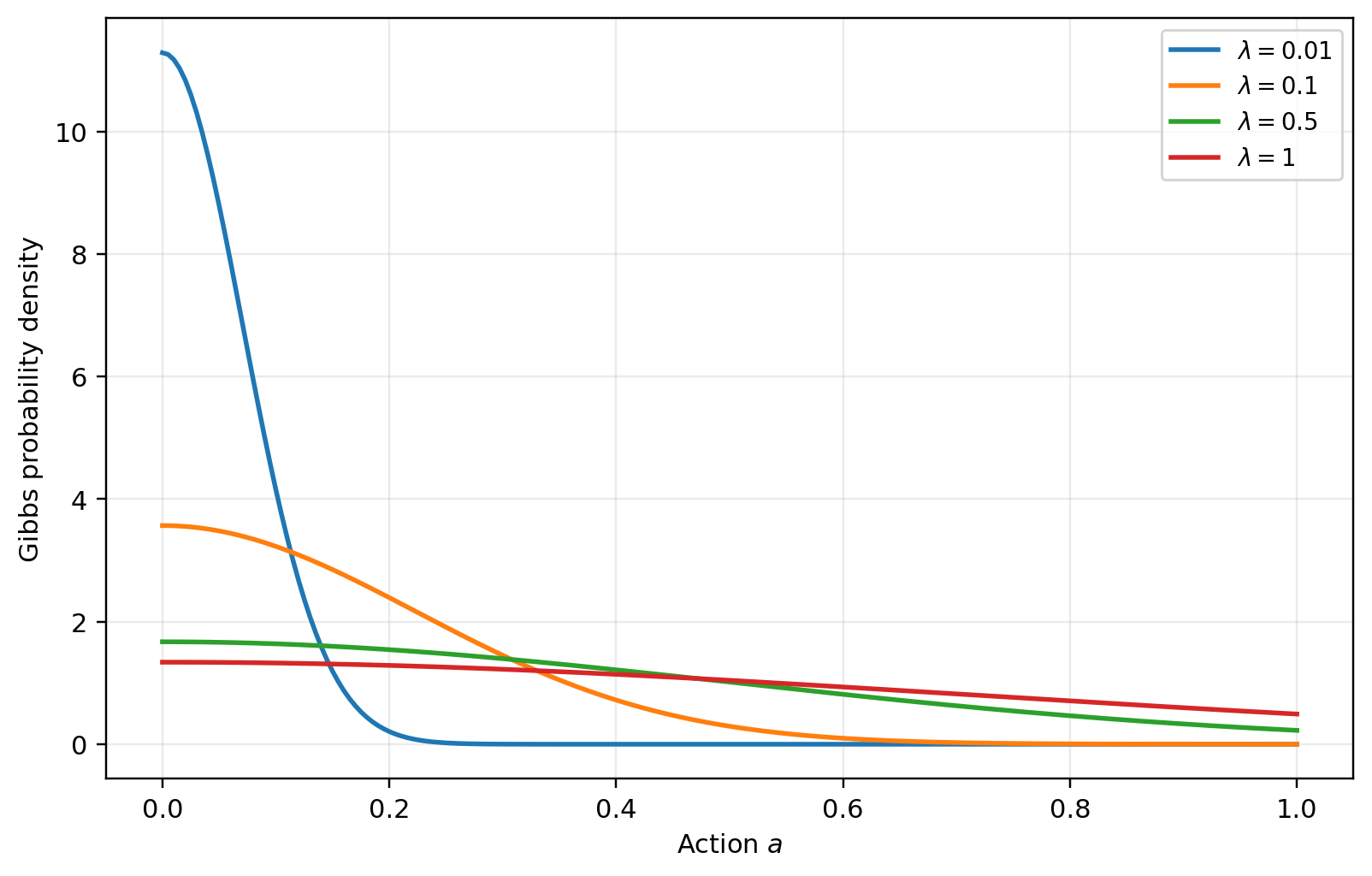}
    \caption{Infected policies}
    \label{fig:Fig7}
  \end{subfigure}\hfill
  \begin{subfigure}[t]{0.24\linewidth}
    \centering
    \includegraphics[width=\linewidth]{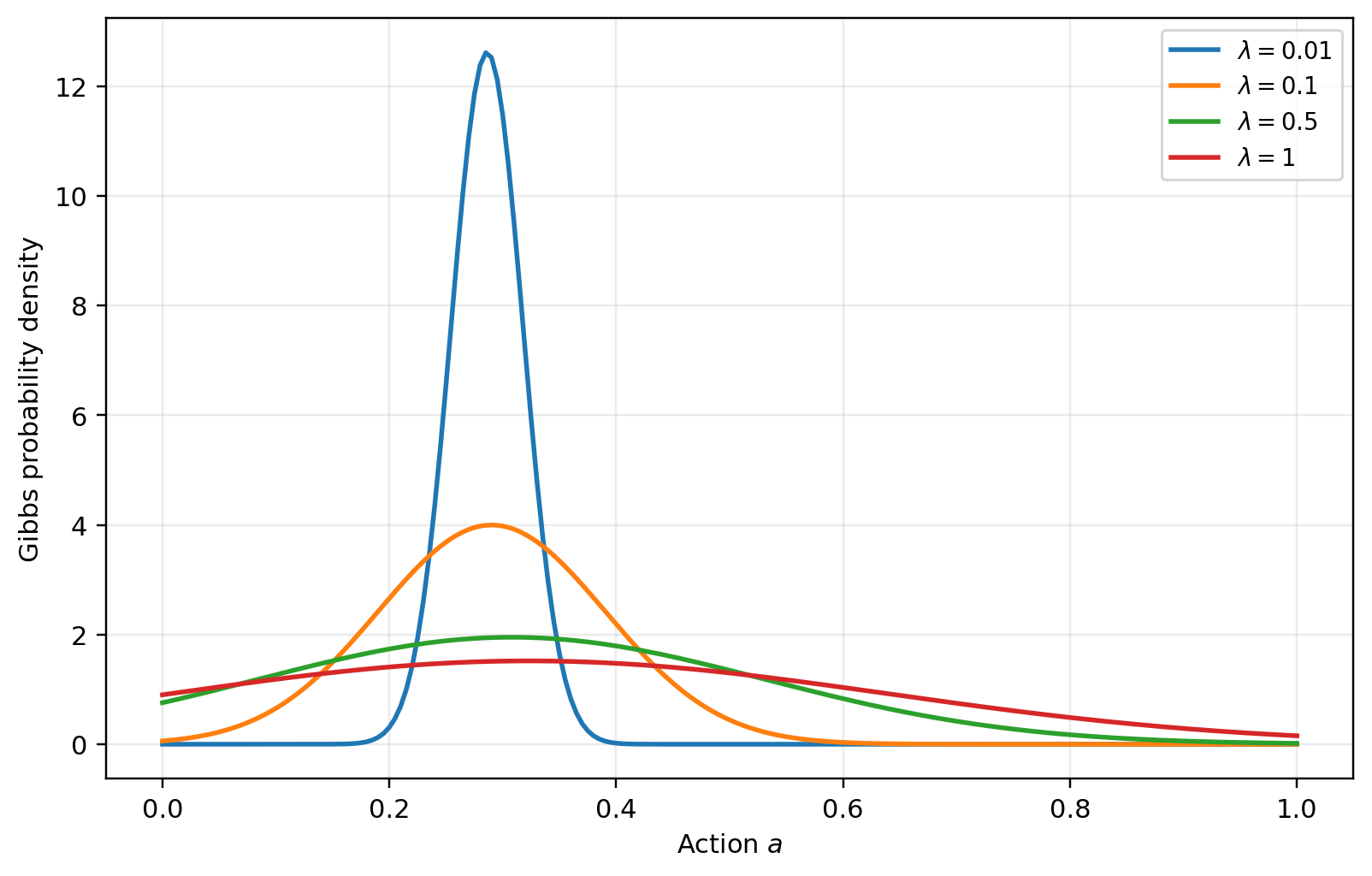}
    \caption{Susceptible policies}
    \label{fig:Fig8}
  \end{subfigure}

  \caption{Sensitivity analysis with respect to the temperature parameter $\lambda$}
  \label{fig:Fig5-8}
\end{figure}

We next investigate the sensitivity of the optimal value functions and relaxed policies with respect to the temperature rate $\lambda$, as illustrated in Figure~\ref{fig:Fig5-8}. Figure~\ref{fig:Fig5} and ~\ref{fig:Fig6} show the optimal value functions for different values of $\lambda \in \{0, 0.001, 0.01, 0.05, 0.1, 0.2, 0.5 , 1\}$ in the original and transformed coordinates, respectively. In both cases, the value functions decrease as $\lambda$ increases. 
Figures~\ref{fig:Fig7} and~\ref{fig:Fig8} illustrate the corresponding optimal policy densities for the infected and susceptible populations at $p=0.5$ for different values of $\lambda \in \{0.01, 0.1, 0.5, 1\}$. As $\lambda$ increases, the optimal policies become more dispersed. For sufficiently small $\lambda$, the policy densities become sharply concentrated around the optimal actions, approaching Dirac measures as $\lambda \to 0$; these cases  $\lambda=0$  and $\lambda=0.001$ are therefore omitted from the plots for clarity.

\ \\
\textbf{Acknowledgements}: Zhenhua Wang is supported by Shandong Excellent Young Scientists Fund Program (Overseas) under grant No. 2025HWYQ–022 and National Natural Science Foundation of China (NSFC) under grant No.12501659. Xiaoli Wei is supported by National Natural Science Foundation of China
grants under no.12201343 and no.12571509. 
Xiang Yu is supported by the Hong Kong RGC General Research Fund (GRF) under grants no. 15211524 and no. 15214125.

\end{document}